\documentclass[opre,nonblindrev]{informs3} 

\DoubleSpacedXI 
\usepackage[colorlinks=true,breaklinks=true,bookmarks=true,urlcolor=blue,
     citecolor=blue,linkcolor=blue,bookmarksopen=false,draft=false]{hyperref}

\usepackage{endnotes}
\let\footnote=\endnote
\let\enotesize=\normalsize
\def\notesname{Endnotes}%
\def\makeenmark{$^{\theenmark}$}
\def\enoteformat{\rightskip0pt\leftskip0pt\parindent=1.75em
  \leavevmode\llap{\theenmark.\enskip}}

\usepackage{natbib}
 \bibpunct[, ]{(}{)}{,}{a}{}{,}%
 \def\bibfont{\small}%
\TheoremsNumberedThrough     
\ECRepeatTheorems

\EquationsNumberedThrough    

\usepackage{amsmath}
\usepackage{amssymb}
\usepackage{empheq}
\usepackage[mathscr]{euscript}
\usepackage{stmaryrd}
\usepackage{mathabx}
\usepackage{dsfont}
\usepackage{graphicx}
\usepackage{multirow}
\usepackage{algorithm}
\usepackage[noend]{algpseudocode}
\usepackage[table]{xcolor}
\usepackage{multirow}
\usepackage{setspace}
\usepackage{inconsolata} 
\usepackage{listings}
\usepackage{arydshln}
\usepackage{tikz}
\usetikzlibrary{3d}
\usepackage{xfrac}
\usepackage{braket}
\usepackage{relsize}
\newcommand*\Let[2]{\State #1 $\gets$ #2}

\usepackage{listings}

\definecolor{Ablue}{RGB}{43,131,186}
\definecolor{Agreen}{RGB}{216,100,118}

\lstdefinelanguage{julia}
{
  keywordsprefix=,
  morekeywords={
    exit,whos,edit,load,is,isa,isequal,typeof,tuple,ntuple,uid,hash,finalizer,convert,promote,
    subtype,typemin,typemax,realmin,realmax,sizeof,promote_type,method_exists,applicable,
    invoke,dlopen,dlsym,system,error,throw,assert,new,Inf,Nan,pi,im,begin,while,for,in,return,
    break,continue,macro,quote,let,if,elseif,else,try,catch,end,bitstype,ccall,do,using,module,
    import,export,importall,baremodule,immutable,local,global,const,Bool,Int,Int8,Int16,Int32,
    Int64,Uint,Uint8,Uint16,Uint32,Uint64,Float32,Float64,Complex64,Complex128,Any,Nothing,None,
    function,type,typealias,abstract
  },
  sensitive=true,
  morecomment=[l]{\#},
  morestring=[b]',
  morestring=[b]"
}

\newcommand{\blue}[1]{{#1}}
\newcommand{\scrA}{\mathscr{A}}

\newcommand{\bbN}{\mathbb{N}}
\newcommand{\bbR}{\mathbb{R}}
\newcommand{\bbZ}{\mathbb{Z}}

\newcommand{\calT}{\mathcal{T}}

\newcommand{\aff}{\operatorname{aff}}

\newcommand{\Conv}{\operatorname{Conv}}
\newcommand{\dom}{\operatorname{dom}}
\newcommand{\Em}{\operatorname{Em}}
\newcommand{\gr}{\textbf{gr}}
\newcommand{\ext}{\operatorname{ext}}
\renewcommand{\int}{\operatorname{int}}

\newcommand{\Proj}{\operatorname{Proj}}
\newcommand{\relint}{\operatorname{relint}}

\newcommand{\supp}{\operatorname{supp}}

\newcommand{\Vol}{\operatorname{Vol}}

\newcommand\defeq{\mathrel{\overset{\makebox[0pt]{\mbox{\normalfont\tiny\sffamily def}}}{=}}}

\newcommand{\norm}[1]{\left\lVert#1\right\rVert}
\newcommand{\sidx}[1]{\left\llbracket     #1 \right\rrbracket}

\usepackage[TABBOTCAP]{subfigure}

\newcommand{\ZZB}{\texttt{ZZB}}
\newcommand{\ZZI}{\texttt{ZZI}}
\newcommand{\MC}{\texttt{MC}}
\newcommand{\CC}{\texttt{CC}}
\newcommand{\Log}{\texttt{LogE}}
\newcommand{\LogIB}{\texttt{LogIB}}
\newcommand{\DLog}{\texttt{DLog}}
\newcommand{\Inc}{\texttt{Inc}}
\newcommand{\SOS}{\texttt{SOS2}}

\begin{document}


\RUNAUTHOR{Huchette and Vielma}

\RUNTITLE{MIP formulations for piecewise linear functions}

\TITLE{Nonconvex piecewise linear functions: Advanced formulations and simple modeling tools}

\ARTICLEAUTHORS{%
\AUTHOR{Joey Huchette}
\AFF{Department of Computational and Applied Mathematics, Rice University, \EMAIL{joehuchette@rice.edu}} 
\AUTHOR{Juan Pablo Vielma}
\AFF{Sloan School of Management, Massachusetts Institute of Technology, \EMAIL{jvielma@mit.edu}}
} 

\ABSTRACT{%
We present novel mixed-integer programming (MIP) formulations for optimization over nonconvex piecewise linear functions. We exploit recent advances in the systematic construction of MIP formulations to derive new formulations for univariate functions using a geometric approach, and for bivariate functions using a combinatorial approach. All formulations are strong, small (so-called \emph{logarithmic} formulations), and have other desirable computational properties. We present extensive experiments in which they exhibit substantial computational performance improvements over existing approaches. To accompany these advanced formulations, we present \texttt{PiecewiseLinearOpt}, an extension of the JuMP modeling language in Julia that implements our models (alongside other formulations from the literature) through a high-level interface, hiding the complexity of the formulations from the end-user.
}%



\maketitle

%


\section{Introduction}
Consider a piecewise linear function $f: D \to \bbR$, where $D \subseteq \bbR^n$. That is, $f$ can be described by a partition of the domain $D$ into a finite family $\{C^i\}_{i=1}^d$ of polyhedral pieces, where for each piece $C^i$ there is an affine function $f^i : C^i \to \bbR$ such that $f(x) = f^i(x)$ for all $x \in C^i$. In this work, we will study methods to solve optimization problems containing piecewise linear functions. This encompasses cases where $f$ appears either in the objective function (e.g. $\min_{x} f(x)$), or in a constraint (e.g. the feasible domain for the optimization problem is partially defined by the inequality $f(x) \leq 0$).

The potential applications for this class of optimization problems are legion. Piecewise linear functions arise naturally throughout operations~\citep{Croxton:2003,Croxton:2007,Liu:2015a} and engineering~\citep{Fugenschuh:2014,Graf:1990,Silva:2012}. They are a natural choice for approximating nonlinear functions, as they often lead to optimization problems that are easier to solve than the original problem~\citep{Bergamini:2005,Bergamini:2008,Castro:2013,Geisler:2012,Kolodziej:2013,Misener:2011,Misener:2012}. For example, there has been recently been significant interest in using piecewise linear functions to approximate complex nonlinearities arising in gas network optimization~\citep{Codas:2012,Codas:2012a,Martin:2006,Mahlke:2010,Misener:2009,Silva:2014}; see \cite{Koch:2015} for a recent book on the subject.

If the function $f$ happens to be convex, it is possible to reformulate our optimization problem into an equivalent linear programming (LP) problem (provided that $D$ is polyhedral). However, if $f$ is nonconvex, this problem is NP-hard in general~\citep{Keha:2006}. A number of specialized algorithms for solving piecewise linear optimization problems have been proposed over the years~\citep{Beale:1970,Farias-Jr.:2008,Farias-Jr.:2013,Keha:2006,Tomlin:1981}. Another popular approach is to use mixed-integer programming (MIP) to encode the logical constraints $x \in C^i \implies f(x) = f^i(x)$ using auxiliary integer decision variables. There are many possible ways to do this, and the MIP approach to modeling optimization problems containing piecewise linear functions has been an active and fruitful area of research for decades~\citep{Balakrishnan:1989,Croxton:2003,DAmbrosio:2010,Dantzig:1960,Jeroslow:1984,Jeroslow:1985a,Keha:2004,Lee:2001,Magnanti:2004,Markowitz:1957,Padberg:2000,Sherali:2001,Vielma:2009a,Vielma:2010,Wilson:1998}. This line of work has produced a large number of MIP formulations that exploit the high performance and flexibility of modern MIP solvers~\citep{Bixby:2007,Junger:2010}, with varying degrees of success. The 2010 \emph{Operations Research} paper of \cite{Vielma:2010} compiled these formulations into a unified framework and provided extensive comparisons of their computational performance. Notably, they showcase the substantial computational advantage of \emph{logarithmic} formulations~\citep{Vielma:2009a}, so-called because their size scales logarithmically in the number of piecewise segments. This work has subsequently sparked attempts to construct logarithmic formulations for other nonconvex constraints~\citep{Huchette:2016a,Huchette:2015,Vielma:2016}. However, the complexity of the logarithmic formulations has resulted in a relatively low rate of adoption in practice, despite their computational efficacy.

In this paper, we study piecewise linear functions as a case study for recent developments in the systematic construction of advanced MIP formulations for nonconvex structures. We present novel logarithmic formulations for piecewise linear functions that improve on the state-of-the-art, and also provide accessible software modeling tools that hide the resulting complexity of these formulations from end users. Specifically, the main contributions of this paper are:
\begin{enumerate}
	\item \blue{\textbf{For univariate functions: A speed-up of up to 3x on harder instances.} In Section~\ref{sec:univariate} we present new formulations for univariate piecewise linear functions that preserve the size and strength of the existing logarithmic formulations, while significantly improving their branching behavior. We show how these formulations computationally outperform the crowded field of existing formulations. In particular, we focus on regimes that are known to be problematic for existing formulations, and observe an improvement of up to 3x on harder instances, though the magnitude of improvement varies among different solvers. To accomplish this, we adapt the geometric formulation construction technique of \cite{Vielma:2016} to develop an unorthodox MIP formulation that exploits general integer (rather than binary) variables. We believe that our results suggest that general integer formulations are a fruitful direction for future MIP formulation research.}

	\item \textbf{For bivariate functions: An order-of-magnitude speed-up.} In Section~\ref{sec:bivariate} we study bivariate piecewise linear functions with generic grid triangulated domains, extending and applying the combinatorial formulation construction technique of~\cite{Huchette:2016a} to develop several families of novel logarithmic formulations. Along the way, we show that for the disjunctive constraints considered in this work (the vast class of ``combinatorial disjunctive constraints''~\citep{Huchette:2016a}), the common loss of strength resulting from intersecting MIP formulations is entirely avoided (Theorem~\ref{thm:combining-formulations}). Finally, we show that the formulations we derive offer a significant computational advantage over existing techniques.

     \item {\bf An accessible modeling library for advanced formulations.} In Section~\ref{sec:computational-tools}, we present a \texttt{PiecewiseLinearOpt}, an extension of the JuMP algebraic modeling language~\citep{Dunning:2015a} that offers a high-level way to model piecewise linear functions in practice. The package supports all the MIP formulations for piecewise linear functions discussed in this work, and generates them automatically and transparently from the user. We believe that easy-to-use modeling interfaces such as \texttt{PiecewiseLinearOpt} are crucial for the practical adoption of advanced MIP formulations like those presented in this work.
\end{enumerate}

\section{Piecewise linear functions and combinatorial disjunctive constraints}
\begin{figure}[htpb]
\centering
    \begin{tikzpicture}
    \begin{scope}
      \draw [thick,->] (-1,0) -- (6,0);
      \draw [thick,->] (-1,0) -- (-1,5);
      \node [right] at (6,0) {$x$};
      \node [above] at (-1,5) {$z$};

      \draw (0,1) -- (1.25,3.5) -- (2.25,3.5) -- (4.5,2) -- (5,4);

      \node [right] at (5,4) {$\gr(f)$};
    \end{scope}
    \end{tikzpicture}
    \hfill
{
    \begin{tikzpicture}[scale=0.4]
\begin{scope}[canvas is xz plane at y=0]
  \draw [fill=gray!20] (0,0) -- (8,0) -- (8,8) -- (0,8);
\draw (0, 8) -- (1, 7) -- (1, 8) -- (0, 8);
\draw (0, 8) -- (1, 7) -- (0, 7) -- (0, 8);
\draw (1, 8) -- (2, 8) -- (1, 7) -- (1, 8);
\draw (2, 7) -- (2, 8) -- (1, 7) -- (2, 7);
\draw (2, 8) -- (3, 7) -- (3, 8) -- (2, 8);
\draw (2, 8) -- (3, 7) -- (2, 7) -- (2, 8);
\draw (3, 8) -- (4, 8) -- (3, 7) -- (3, 8);
\draw (4, 7) -- (4, 8) -- (3, 7) -- (4, 7);
\draw (4, 8) -- (5, 8) -- (4, 7) -- (4, 8);
\draw (5, 7) -- (5, 8) -- (4, 7) -- (5, 7);
\draw (5, 8) -- (6, 8) -- (5, 7) -- (5, 8);
\draw (6, 7) -- (6, 8) -- (5, 7) -- (6, 7);
\draw (6, 8) -- (7, 8) -- (6, 7) -- (6, 8);
\draw (7, 7) -- (7, 8) -- (6, 7) -- (7, 7);
\draw (7, 8) -- (8, 7) -- (8, 8) -- (7, 8);
\draw (7, 8) -- (8, 7) -- (7, 7) -- (7, 8);
\draw (0, 7) -- (1, 7) -- (0, 6) -- (0, 7);
\draw (1, 6) -- (1, 7) -- (0, 6) -- (1, 6);
\draw (1, 7) -- (2, 7) -- (1, 6) -- (1, 7);
\draw (2, 6) -- (2, 7) -- (1, 6) -- (2, 6);
\draw (2, 7) -- (3, 6) -- (3, 7) -- (2, 7);
\draw (2, 7) -- (3, 6) -- (2, 6) -- (2, 7);
\draw (3, 7) -- (4, 6) -- (4, 7) -- (3, 7);
\draw (3, 7) -- (4, 6) -- (3, 6) -- (3, 7);
\draw (4, 7) -- (5, 7) -- (4, 6) -- (4, 7);
\draw (5, 6) -- (5, 7) -- (4, 6) -- (5, 6);
\draw (5, 7) -- (6, 6) -- (6, 7) -- (5, 7);
\draw (5, 7) -- (6, 6) -- (5, 6) -- (5, 7);
\draw (6, 7) -- (7, 7) -- (6, 6) -- (6, 7);
\draw (7, 6) -- (7, 7) -- (6, 6) -- (7, 6);
\draw (7, 7) -- (8, 7) -- (7, 6) -- (7, 7);
\draw (8, 6) -- (8, 7) -- (7, 6) -- (8, 6);
\draw (0, 6) -- (1, 5) -- (1, 6) -- (0, 6);
\draw (0, 6) -- (1, 5) -- (0, 5) -- (0, 6);
\draw (1, 6) -- (2, 5) -- (2, 6) -- (1, 6);
\draw (1, 6) -- (2, 5) -- (1, 5) -- (1, 6);
\draw (2, 6) -- (3, 5) -- (3, 6) -- (2, 6);
\draw (2, 6) -- (3, 5) -- (2, 5) -- (2, 6);
\draw (3, 6) -- (4, 5) -- (4, 6) -- (3, 6);
\draw (3, 6) -- (4, 5) -- (3, 5) -- (3, 6);
\draw (4, 6) -- (5, 6) -- (4, 5) -- (4, 6);
\draw (5, 5) -- (5, 6) -- (4, 5) -- (5, 5);
\draw (5, 6) -- (6, 5) -- (6, 6) -- (5, 6);
\draw (5, 6) -- (6, 5) -- (5, 5) -- (5, 6);
\draw (6, 6) -- (7, 6) -- (6, 5) -- (6, 6);
\draw (7, 5) -- (7, 6) -- (6, 5) -- (7, 5);
\draw (7, 6) -- (8, 6) -- (7, 5) -- (7, 6);
\draw (8, 5) -- (8, 6) -- (7, 5) -- (8, 5);
\draw (0, 5) -- (1, 5) -- (0, 4) -- (0, 5);
\draw (1, 4) -- (1, 5) -- (0, 4) -- (1, 4);
\draw (1, 5) -- (2, 4) -- (2, 5) -- (1, 5);
\draw (1, 5) -- (2, 4) -- (1, 4) -- (1, 5);
\draw (2, 5) -- (3, 4) -- (3, 5) -- (2, 5);
\draw (2, 5) -- (3, 4) -- (2, 4) -- (2, 5);
\draw (3, 5) -- (4, 4) -- (4, 5) -- (3, 5);
\draw (3, 5) -- (4, 4) -- (3, 4) -- (3, 5);
\draw (4, 5) -- (5, 4) -- (5, 5) -- (4, 5);
\draw (4, 5) -- (5, 4) -- (4, 4) -- (4, 5);
\draw (5, 5) -- (6, 4) -- (6, 5) -- (5, 5);
\draw (5, 5) -- (6, 4) -- (5, 4) -- (5, 5);
\draw (6, 5) -- (7, 4) -- (7, 5) -- (6, 5);
\draw (6, 5) -- (7, 4) -- (6, 4) -- (6, 5);
\draw (7, 5) -- (8, 5) -- (7, 4) -- (7, 5);
\draw (8, 4) -- (8, 5) -- (7, 4) -- (8, 4);
\draw (0, 4) -- (1, 4) -- (0, 3) -- (0, 4);
\draw (1, 3) -- (1, 4) -- (0, 3) -- (1, 3);
\draw (1, 4) -- (2, 4) -- (1, 3) -- (1, 4);
\draw (2, 3) -- (2, 4) -- (1, 3) -- (2, 3);
\draw (2, 4) -- (3, 4) -- (2, 3) -- (2, 4);
\draw (3, 3) -- (3, 4) -- (2, 3) -- (3, 3);
\draw (3, 4) -- (4, 4) -- (3, 3) -- (3, 4);
\draw (4, 3) -- (4, 4) -- (3, 3) -- (4, 3);
\draw (4, 4) -- (5, 4) -- (4, 3) -- (4, 4);
\draw (5, 3) -- (5, 4) -- (4, 3) -- (5, 3);
\draw (5, 4) -- (6, 4) -- (5, 3) -- (5, 4);
\draw (6, 3) -- (6, 4) -- (5, 3) -- (6, 3);
\draw (6, 4) -- (7, 4) -- (6, 3) -- (6, 4);
\draw (7, 3) -- (7, 4) -- (6, 3) -- (7, 3);
\draw (7, 4) -- (8, 3) -- (8, 4) -- (7, 4);
\draw (7, 4) -- (8, 3) -- (7, 3) -- (7, 4);
\draw (0, 3) -- (1, 3) -- (0, 2) -- (0, 3);
\draw (1, 2) -- (1, 3) -- (0, 2) -- (1, 2);
\draw (1, 3) -- (2, 3) -- (1, 2) -- (1, 3);
\draw (2, 2) -- (2, 3) -- (1, 2) -- (2, 2);
\draw (2, 3) -- (3, 2) -- (3, 3) -- (2, 3);
\draw (2, 3) -- (3, 2) -- (2, 2) -- (2, 3);
\draw (3, 3) -- (4, 2) -- (4, 3) -- (3, 3);
\draw (3, 3) -- (4, 2) -- (3, 2) -- (3, 3);
\draw (4, 3) -- (5, 3) -- (4, 2) -- (4, 3);
\draw (5, 2) -- (5, 3) -- (4, 2) -- (5, 2);
\draw (5, 3) -- (6, 2) -- (6, 3) -- (5, 3);
\draw (5, 3) -- (6, 2) -- (5, 2) -- (5, 3);
\draw (6, 3) -- (7, 2) -- (7, 3) -- (6, 3);
\draw (6, 3) -- (7, 2) -- (6, 2) -- (6, 3);
\draw (7, 3) -- (8, 2) -- (8, 3) -- (7, 3);
\draw (7, 3) -- (8, 2) -- (7, 2) -- (7, 3);
\draw (0, 2) -- (1, 2) -- (0, 1) -- (0, 2);
\draw (1, 1) -- (1, 2) -- (0, 1) -- (1, 1);
\draw (1, 2) -- (2, 1) -- (2, 2) -- (1, 2);
\draw (1, 2) -- (2, 1) -- (1, 1) -- (1, 2);
\draw (2, 2) -- (3, 2) -- (2, 1) -- (2, 2);
\draw (3, 1) -- (3, 2) -- (2, 1) -- (3, 1);
\draw (3, 2) -- (4, 2) -- (3, 1) -- (3, 2);
\draw (4, 1) -- (4, 2) -- (3, 1) -- (4, 1);
\draw (4, 2) -- (5, 2) -- (4, 1) -- (4, 2);
\draw (5, 1) -- (5, 2) -- (4, 1) -- (5, 1);
\draw (5, 2) -- (6, 1) -- (6, 2) -- (5, 2);
\draw (5, 2) -- (6, 1) -- (5, 1) -- (5, 2);
\draw (6, 2) -- (7, 2) -- (6, 1) -- (6, 2);
\draw (7, 1) -- (7, 2) -- (6, 1) -- (7, 1);
\draw (7, 2) -- (8, 1) -- (8, 2) -- (7, 2);
\draw (7, 2) -- (8, 1) -- (7, 1) -- (7, 2);
\draw (0, 1) -- (1, 1) -- (0, 0) -- (0, 1);
\draw (1, 0) -- (1, 1) -- (0, 0) -- (1, 0);
\draw (1, 1) -- (2, 1) -- (1, 0) -- (1, 1);
\draw (2, 0) -- (2, 1) -- (1, 0) -- (2, 0);
\draw (2, 1) -- (3, 1) -- (2, 0) -- (2, 1);
\draw (3, 0) -- (3, 1) -- (2, 0) -- (3, 0);
\draw (3, 1) -- (4, 0) -- (4, 1) -- (3, 1);
\draw (3, 1) -- (4, 0) -- (3, 0) -- (3, 1);
\draw (4, 1) -- (5, 0) -- (5, 1) -- (4, 1);
\draw (4, 1) -- (5, 0) -- (4, 0) -- (4, 1);
\draw (5, 1) -- (6, 1) -- (5, 0) -- (5, 1);
\draw (6, 0) -- (6, 1) -- (5, 0) -- (6, 0);
\draw (6, 1) -- (7, 1) -- (6, 0) -- (6, 1);
\draw (7, 0) -- (7, 1) -- (6, 0) -- (7, 0);
\draw (7, 1) -- (8, 1) -- (7, 0) -- (7, 1);
\draw (8, 0) -- (8, 1) -- (7, 0) -- (8, 0);

  \draw [->, thick] (0,0) -- (0,9);
  \draw [->, thick] (0,0) -- (9,0);

  \node [below] at (0,9) {$x_1$};
  \node [right] at (9,0) {$x_2$};
\end{scope}

\begin{scope}[canvas is xy plane at z=0]
  \draw [->, thick] (0,0) -- (0,8);
  \node [above] at (0,8) {$z$};
  \node [right] at (8,6) {$\gr(f)$};
\end{scope}

\draw (0, 6.0, 8) -- (1, 7.060660171779821, 7) -- (1, 6.574025148547634, 8) -- (0, 6.0, 8);
\draw (0, 6.0, 8) -- (1, 7.060660171779821, 7) -- (0, 6.574025148547634, 7) -- (0, 6.0, 8);
\draw (1, 6.574025148547634, 8) -- (2, 7.060660171779821, 8) -- (1, 7.060660171779821, 7) -- (1, 6.574025148547634, 8);
\draw (2, 7.38581929876693, 7) -- (2, 7.060660171779821, 8) -- (1, 7.060660171779821, 7) -- (2, 7.38581929876693, 7);
\draw (2, 7.060660171779821, 8) -- (3, 7.5, 7) -- (3, 7.38581929876693, 8) -- (2, 7.060660171779821, 8);
\draw (2, 7.060660171779821, 8) -- (3, 7.5, 7) -- (2, 7.38581929876693, 7) -- (2, 7.060660171779821, 8);
\draw (3, 7.38581929876693, 8) -- (4, 7.5, 8) -- (3, 7.5, 7) -- (3, 7.38581929876693, 8);
\draw (4, 7.38581929876693, 7) -- (4, 7.5, 8) -- (3, 7.5, 7) -- (4, 7.38581929876693, 7);
\draw (4, 7.5, 8) -- (5, 7.38581929876693, 8) -- (4, 7.38581929876693, 7) -- (4, 7.5, 8);
\draw (5, 7.060660171779821, 7) -- (5, 7.38581929876693, 8) -- (4, 7.38581929876693, 7) -- (5, 7.060660171779821, 7);
\draw (5, 7.38581929876693, 8) -- (6, 7.060660171779821, 8) -- (5, 7.060660171779821, 7) -- (5, 7.38581929876693, 8);
\draw (6, 6.574025148547635, 7) -- (6, 7.060660171779821, 8) -- (5, 7.060660171779821, 7) -- (6, 6.574025148547635, 7);
\draw (6, 7.060660171779821, 8) -- (7, 6.574025148547635, 8) -- (6, 6.574025148547635, 7) -- (6, 7.060660171779821, 8);
\draw (7, 6.0, 7) -- (7, 6.574025148547635, 8) -- (6, 6.574025148547635, 7) -- (7, 6.0, 7);
\draw (7, 6.574025148547635, 8) -- (8, 5.425974851452366, 7) -- (8, 6.0, 8) -- (7, 6.574025148547635, 8);
\draw (7, 6.574025148547635, 8) -- (8, 5.425974851452366, 7) -- (7, 6.0, 7) -- (7, 6.574025148547635, 8);
\draw (0, 6.574025148547634, 7) -- (1, 7.060660171779821, 7) -- (0, 7.060660171779821, 6) -- (0, 6.574025148547634, 7);
\draw (1, 7.38581929876693, 6) -- (1, 7.060660171779821, 7) -- (0, 7.060660171779821, 6) -- (1, 7.38581929876693, 6);
\draw (1, 7.060660171779821, 7) -- (2, 7.38581929876693, 7) -- (1, 7.38581929876693, 6) -- (1, 7.060660171779821, 7);
\draw (2, 7.5, 6) -- (2, 7.38581929876693, 7) -- (1, 7.38581929876693, 6) -- (2, 7.5, 6);
\draw (2, 7.38581929876693, 7) -- (3, 7.38581929876693, 6) -- (3, 7.5, 7) -- (2, 7.38581929876693, 7);
\draw (2, 7.38581929876693, 7) -- (3, 7.38581929876693, 6) -- (2, 7.5, 6) -- (2, 7.38581929876693, 7);
\draw (3, 7.5, 7) -- (4, 7.060660171779821, 6) -- (4, 7.38581929876693, 7) -- (3, 7.5, 7);
\draw (3, 7.5, 7) -- (4, 7.060660171779821, 6) -- (3, 7.38581929876693, 6) -- (3, 7.5, 7);
\draw (4, 7.38581929876693, 7) -- (5, 7.060660171779821, 7) -- (4, 7.060660171779821, 6) -- (4, 7.38581929876693, 7);
\draw (5, 6.574025148547635, 6) -- (5, 7.060660171779821, 7) -- (4, 7.060660171779821, 6) -- (5, 6.574025148547635, 6);
\draw (5, 7.060660171779821, 7) -- (6, 6.0, 6) -- (6, 6.574025148547635, 7) -- (5, 7.060660171779821, 7);
\draw (5, 7.060660171779821, 7) -- (6, 6.0, 6) -- (5, 6.574025148547635, 6) -- (5, 7.060660171779821, 7);
\draw (6, 6.574025148547635, 7) -- (7, 6.0, 7) -- (6, 6.0, 6) -- (6, 6.574025148547635, 7);
\draw (7, 5.425974851452366, 6) -- (7, 6.0, 7) -- (6, 6.0, 6) -- (7, 5.425974851452366, 6);
\draw (7, 6.0, 7) -- (8, 5.425974851452366, 7) -- (7, 5.425974851452366, 6) -- (7, 6.0, 7);
\draw (8, 4.939339828220179, 6) -- (8, 5.425974851452366, 7) -- (7, 5.425974851452366, 6) -- (8, 4.939339828220179, 6);
\draw (0, 7.060660171779821, 6) -- (1, 7.5, 5) -- (1, 7.38581929876693, 6) -- (0, 7.060660171779821, 6);
\draw (0, 7.060660171779821, 6) -- (1, 7.5, 5) -- (0, 7.38581929876693, 5) -- (0, 7.060660171779821, 6);
\draw (1, 7.38581929876693, 6) -- (2, 7.38581929876693, 5) -- (2, 7.5, 6) -- (1, 7.38581929876693, 6);
\draw (1, 7.38581929876693, 6) -- (2, 7.38581929876693, 5) -- (1, 7.5, 5) -- (1, 7.38581929876693, 6);
\draw (2, 7.5, 6) -- (3, 7.060660171779821, 5) -- (3, 7.38581929876693, 6) -- (2, 7.5, 6);
\draw (2, 7.5, 6) -- (3, 7.060660171779821, 5) -- (2, 7.38581929876693, 5) -- (2, 7.5, 6);
\draw (3, 7.38581929876693, 6) -- (4, 6.574025148547635, 5) -- (4, 7.060660171779821, 6) -- (3, 7.38581929876693, 6);
\draw (3, 7.38581929876693, 6) -- (4, 6.574025148547635, 5) -- (3, 7.060660171779821, 5) -- (3, 7.38581929876693, 6);
\draw (4, 7.060660171779821, 6) -- (5, 6.574025148547635, 6) -- (4, 6.574025148547635, 5) -- (4, 7.060660171779821, 6);
\draw (5, 6.0, 5) -- (5, 6.574025148547635, 6) -- (4, 6.574025148547635, 5) -- (5, 6.0, 5);
\draw (5, 6.574025148547635, 6) -- (6, 5.425974851452366, 5) -- (6, 6.0, 6) -- (5, 6.574025148547635, 6);
\draw (5, 6.574025148547635, 6) -- (6, 5.425974851452366, 5) -- (5, 6.0, 5) -- (5, 6.574025148547635, 6);
\draw (6, 6.0, 6) -- (7, 5.425974851452366, 6) -- (6, 5.425974851452366, 5) -- (6, 6.0, 6);
\draw (7, 4.939339828220179, 5) -- (7, 5.425974851452366, 6) -- (6, 5.425974851452366, 5) -- (7, 4.939339828220179, 5);
\draw (7, 5.425974851452366, 6) -- (8, 4.939339828220179, 6) -- (7, 4.939339828220179, 5) -- (7, 5.425974851452366, 6);
\draw (8, 4.61418070123307, 5) -- (8, 4.939339828220179, 6) -- (7, 4.939339828220179, 5) -- (8, 4.61418070123307, 5);
\draw (0, 7.38581929876693, 5) -- (1, 7.5, 5) -- (0, 7.5, 4) -- (0, 7.38581929876693, 5);
\draw (1, 7.38581929876693, 4) -- (1, 7.5, 5) -- (0, 7.5, 4) -- (1, 7.38581929876693, 4);
\draw (1, 7.5, 5) -- (2, 7.060660171779821, 4) -- (2, 7.38581929876693, 5) -- (1, 7.5, 5);
\draw (1, 7.5, 5) -- (2, 7.060660171779821, 4) -- (1, 7.38581929876693, 4) -- (1, 7.5, 5);
\draw (2, 7.38581929876693, 5) -- (3, 6.574025148547635, 4) -- (3, 7.060660171779821, 5) -- (2, 7.38581929876693, 5);
\draw (2, 7.38581929876693, 5) -- (3, 6.574025148547635, 4) -- (2, 7.060660171779821, 4) -- (2, 7.38581929876693, 5);
\draw (3, 7.060660171779821, 5) -- (4, 6.0, 4) -- (4, 6.574025148547635, 5) -- (3, 7.060660171779821, 5);
\draw (3, 7.060660171779821, 5) -- (4, 6.0, 4) -- (3, 6.574025148547635, 4) -- (3, 7.060660171779821, 5);
\draw (4, 6.574025148547635, 5) -- (5, 5.425974851452366, 4) -- (5, 6.0, 5) -- (4, 6.574025148547635, 5);
\draw (4, 6.574025148547635, 5) -- (5, 5.425974851452366, 4) -- (4, 6.0, 4) -- (4, 6.574025148547635, 5);
\draw (5, 6.0, 5) -- (6, 4.939339828220179, 4) -- (6, 5.425974851452366, 5) -- (5, 6.0, 5);
\draw (5, 6.0, 5) -- (6, 4.939339828220179, 4) -- (5, 5.425974851452366, 4) -- (5, 6.0, 5);
\draw (6, 5.425974851452366, 5) -- (7, 4.61418070123307, 4) -- (7, 4.939339828220179, 5) -- (6, 5.425974851452366, 5);
\draw (6, 5.425974851452366, 5) -- (7, 4.61418070123307, 4) -- (6, 4.939339828220179, 4) -- (6, 5.425974851452366, 5);
\draw (7, 4.939339828220179, 5) -- (8, 4.61418070123307, 5) -- (7, 4.61418070123307, 4) -- (7, 4.939339828220179, 5);
\draw (8, 4.5, 4) -- (8, 4.61418070123307, 5) -- (7, 4.61418070123307, 4) -- (8, 4.5, 4);
\draw (0, 7.5, 4) -- (1, 7.38581929876693, 4) -- (0, 7.38581929876693, 3) -- (0, 7.5, 4);
\draw (1, 7.060660171779821, 3) -- (1, 7.38581929876693, 4) -- (0, 7.38581929876693, 3) -- (1, 7.060660171779821, 3);
\draw (1, 7.38581929876693, 4) -- (2, 7.060660171779821, 4) -- (1, 7.060660171779821, 3) -- (1, 7.38581929876693, 4);
\draw (2, 6.574025148547635, 3) -- (2, 7.060660171779821, 4) -- (1, 7.060660171779821, 3) -- (2, 6.574025148547635, 3);
\draw (2, 7.060660171779821, 4) -- (3, 6.574025148547635, 4) -- (2, 6.574025148547635, 3) -- (2, 7.060660171779821, 4);
\draw (3, 6.0, 3) -- (3, 6.574025148547635, 4) -- (2, 6.574025148547635, 3) -- (3, 6.0, 3);
\draw (3, 6.574025148547635, 4) -- (4, 6.0, 4) -- (3, 6.0, 3) -- (3, 6.574025148547635, 4);
\draw (4, 5.425974851452366, 3) -- (4, 6.0, 4) -- (3, 6.0, 3) -- (4, 5.425974851452366, 3);
\draw (4, 6.0, 4) -- (5, 5.425974851452366, 4) -- (4, 5.425974851452366, 3) -- (4, 6.0, 4);
\draw (5, 4.939339828220179, 3) -- (5, 5.425974851452366, 4) -- (4, 5.425974851452366, 3) -- (5, 4.939339828220179, 3);
\draw (5, 5.425974851452366, 4) -- (6, 4.939339828220179, 4) -- (5, 4.939339828220179, 3) -- (5, 5.425974851452366, 4);
\draw (6, 4.61418070123307, 3) -- (6, 4.939339828220179, 4) -- (5, 4.939339828220179, 3) -- (6, 4.61418070123307, 3);
\draw (6, 4.939339828220179, 4) -- (7, 4.61418070123307, 4) -- (6, 4.61418070123307, 3) -- (6, 4.939339828220179, 4);
\draw (7, 4.5, 3) -- (7, 4.61418070123307, 4) -- (6, 4.61418070123307, 3) -- (7, 4.5, 3);
\draw (7, 4.61418070123307, 4) -- (8, 4.61418070123307, 3) -- (8, 4.5, 4) -- (7, 4.61418070123307, 4);
\draw (7, 4.61418070123307, 4) -- (8, 4.61418070123307, 3) -- (7, 4.5, 3) -- (7, 4.61418070123307, 4);
\draw (0, 7.38581929876693, 3) -- (1, 7.060660171779821, 3) -- (0, 7.060660171779821, 2) -- (0, 7.38581929876693, 3);
\draw (1, 6.574025148547635, 2) -- (1, 7.060660171779821, 3) -- (0, 7.060660171779821, 2) -- (1, 6.574025148547635, 2);
\draw (1, 7.060660171779821, 3) -- (2, 6.574025148547635, 3) -- (1, 6.574025148547635, 2) -- (1, 7.060660171779821, 3);
\draw (2, 6.0, 2) -- (2, 6.574025148547635, 3) -- (1, 6.574025148547635, 2) -- (2, 6.0, 2);
\draw (2, 6.574025148547635, 3) -- (3, 5.425974851452366, 2) -- (3, 6.0, 3) -- (2, 6.574025148547635, 3);
\draw (2, 6.574025148547635, 3) -- (3, 5.425974851452366, 2) -- (2, 6.0, 2) -- (2, 6.574025148547635, 3);
\draw (3, 6.0, 3) -- (4, 4.939339828220179, 2) -- (4, 5.425974851452366, 3) -- (3, 6.0, 3);
\draw (3, 6.0, 3) -- (4, 4.939339828220179, 2) -- (3, 5.425974851452366, 2) -- (3, 6.0, 3);
\draw (4, 5.425974851452366, 3) -- (5, 4.939339828220179, 3) -- (4, 4.939339828220179, 2) -- (4, 5.425974851452366, 3);
\draw (5, 4.61418070123307, 2) -- (5, 4.939339828220179, 3) -- (4, 4.939339828220179, 2) -- (5, 4.61418070123307, 2);
\draw (5, 4.939339828220179, 3) -- (6, 4.5, 2) -- (6, 4.61418070123307, 3) -- (5, 4.939339828220179, 3);
\draw (5, 4.939339828220179, 3) -- (6, 4.5, 2) -- (5, 4.61418070123307, 2) -- (5, 4.939339828220179, 3);
\draw (6, 4.61418070123307, 3) -- (7, 4.61418070123307, 2) -- (7, 4.5, 3) -- (6, 4.61418070123307, 3);
\draw (6, 4.61418070123307, 3) -- (7, 4.61418070123307, 2) -- (6, 4.5, 2) -- (6, 4.61418070123307, 3);
\draw (7, 4.5, 3) -- (8, 4.939339828220179, 2) -- (8, 4.61418070123307, 3) -- (7, 4.5, 3);
\draw (7, 4.5, 3) -- (8, 4.939339828220179, 2) -- (7, 4.61418070123307, 2) -- (7, 4.5, 3);
\draw (0, 7.060660171779821, 2) -- (1, 6.574025148547635, 2) -- (0, 6.574025148547635, 1) -- (0, 7.060660171779821, 2);
\draw (1, 6.0, 1) -- (1, 6.574025148547635, 2) -- (0, 6.574025148547635, 1) -- (1, 6.0, 1);
\draw (1, 6.574025148547635, 2) -- (2, 5.425974851452366, 1) -- (2, 6.0, 2) -- (1, 6.574025148547635, 2);
\draw (1, 6.574025148547635, 2) -- (2, 5.425974851452366, 1) -- (1, 6.0, 1) -- (1, 6.574025148547635, 2);
\draw (2, 6.0, 2) -- (3, 5.425974851452366, 2) -- (2, 5.425974851452366, 1) -- (2, 6.0, 2);
\draw (3, 4.939339828220179, 1) -- (3, 5.425974851452366, 2) -- (2, 5.425974851452366, 1) -- (3, 4.939339828220179, 1);
\draw (3, 5.425974851452366, 2) -- (4, 4.939339828220179, 2) -- (3, 4.939339828220179, 1) -- (3, 5.425974851452366, 2);
\draw (4, 4.61418070123307, 1) -- (4, 4.939339828220179, 2) -- (3, 4.939339828220179, 1) -- (4, 4.61418070123307, 1);
\draw (4, 4.939339828220179, 2) -- (5, 4.61418070123307, 2) -- (4, 4.61418070123307, 1) -- (4, 4.939339828220179, 2);
\draw (5, 4.5, 1) -- (5, 4.61418070123307, 2) -- (4, 4.61418070123307, 1) -- (5, 4.5, 1);
\draw (5, 4.61418070123307, 2) -- (6, 4.61418070123307, 1) -- (6, 4.5, 2) -- (5, 4.61418070123307, 2);
\draw (5, 4.61418070123307, 2) -- (6, 4.61418070123307, 1) -- (5, 4.5, 1) -- (5, 4.61418070123307, 2);
\draw (6, 4.5, 2) -- (7, 4.61418070123307, 2) -- (6, 4.61418070123307, 1) -- (6, 4.5, 2);
\draw (7, 4.939339828220179, 1) -- (7, 4.61418070123307, 2) -- (6, 4.61418070123307, 1) -- (7, 4.939339828220179, 1);
\draw (7, 4.61418070123307, 2) -- (8, 5.425974851452365, 1) -- (8, 4.939339828220179, 2) -- (7, 4.61418070123307, 2);
\draw (7, 4.61418070123307, 2) -- (8, 5.425974851452365, 1) -- (7, 4.939339828220179, 1) -- (7, 4.61418070123307, 2);
\draw (0, 6.574025148547635, 1) -- (1, 6.0, 1) -- (0, 6.0, 0) -- (0, 6.574025148547635, 1);
\draw (1, 5.425974851452366, 0) -- (1, 6.0, 1) -- (0, 6.0, 0) -- (1, 5.425974851452366, 0);
\draw (1, 6.0, 1) -- (2, 5.425974851452366, 1) -- (1, 5.425974851452366, 0) -- (1, 6.0, 1);
\draw (2, 4.939339828220179, 0) -- (2, 5.425974851452366, 1) -- (1, 5.425974851452366, 0) -- (2, 4.939339828220179, 0);
\draw (2, 5.425974851452366, 1) -- (3, 4.939339828220179, 1) -- (2, 4.939339828220179, 0) -- (2, 5.425974851452366, 1);
\draw (3, 4.61418070123307, 0) -- (3, 4.939339828220179, 1) -- (2, 4.939339828220179, 0) -- (3, 4.61418070123307, 0);
\draw (3, 4.939339828220179, 1) -- (4, 4.5, 0) -- (4, 4.61418070123307, 1) -- (3, 4.939339828220179, 1);
\draw (3, 4.939339828220179, 1) -- (4, 4.5, 0) -- (3, 4.61418070123307, 0) -- (3, 4.939339828220179, 1);
\draw (4, 4.61418070123307, 1) -- (5, 4.61418070123307, 0) -- (5, 4.5, 1) -- (4, 4.61418070123307, 1);
\draw (4, 4.61418070123307, 1) -- (5, 4.61418070123307, 0) -- (4, 4.5, 0) -- (4, 4.61418070123307, 1);
\draw (5, 4.5, 1) -- (6, 4.61418070123307, 1) -- (5, 4.61418070123307, 0) -- (5, 4.5, 1);
\draw (6, 4.939339828220179, 0) -- (6, 4.61418070123307, 1) -- (5, 4.61418070123307, 0) -- (6, 4.939339828220179, 0);
\draw (6, 4.61418070123307, 1) -- (7, 4.939339828220179, 1) -- (6, 4.939339828220179, 0) -- (6, 4.61418070123307, 1);
\draw (7, 5.425974851452365, 0) -- (7, 4.939339828220179, 1) -- (6, 4.939339828220179, 0) -- (7, 5.425974851452365, 0);
\draw (7, 4.939339828220179, 1) -- (8, 5.425974851452365, 1) -- (7, 5.425974851452365, 0) -- (7, 4.939339828220179, 1);
\draw (8, 6.0, 0) -- (8, 5.425974851452365, 1) -- (7, 5.425974851452365, 0) -- (8, 6.0, 0);
\end{tikzpicture}
}
\caption{\textbf{(Left)} A univariate piecewise linear function, and \textbf{(Right)} a bivariate piecewise linear function with a grid triangulated domain.}
\label{fig:pwl-examples}
\end{figure}
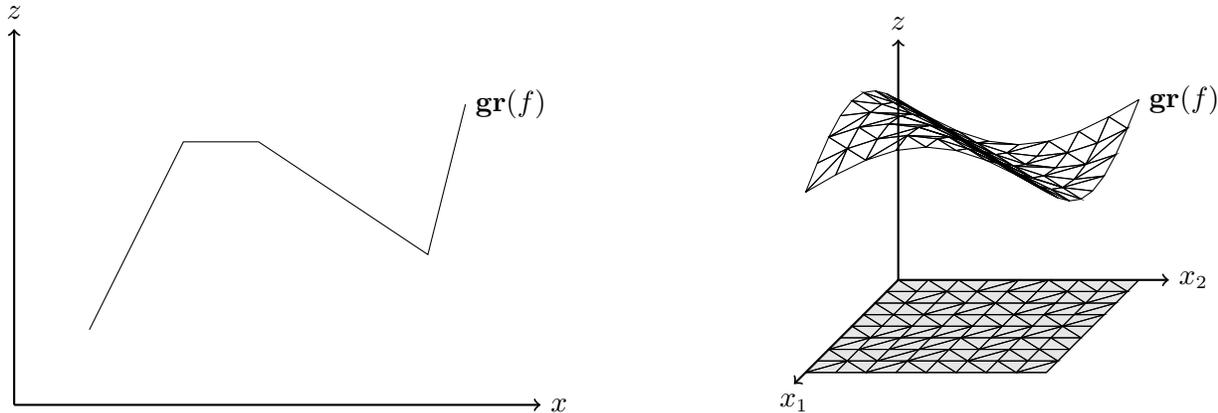

Consider a continuous\endnote{We refer the reader interested in  modeling discontinuous functions to \cite{Vielma:2010}.} piecewise linear function $f : D \to \bbR$, where $D \subset \bbR^n$ is bounded. We will describe $f$ in terms of the domain pieces $\{C^i \subseteq D\}_{i=1}^d$ and affine functions $\{f^i\}_{i=1}^d$ as above; we assume that the pieces cover the domain $D$ and that their interiors do not overlap. Furthermore, we assume that our function $f: \mathbb{R}^n \to \bbR$ is \emph{non-separable} and cannot be decomposed as the sum of lower-dimensional piecewise linear functions. This is without loss of generality, as if such a decomposition exists, we could apply our formulation techniques to the individual pieces separately. Finally, we will focus primarily on the regime where the dimension $n$ of the domain is relatively small: when $f$ is either univariate ($n=1$) or bivariate ($n=2$) with a grid triangulated domain; see Figure~\ref{fig:pwl-examples} for an illustrative example of each. Low dimensional piecewise linear functions are broadly applicable (especially with the non-separability assumption), and are sufficiently complex to warrant in-depth analysis. We tabulate notation we will use for the remainder in Table~\ref{tab:notation}.

\begin{table}[htpb]
\begin{center}
{\scriptsize
\begin{tabular}{r|l|l}
Notation & Formal Definition & Description \\\hline
$\llbracket d \rrbracket$ & $\{1,\ldots,d\}$ & All integers from 1 to $d$ \\
$\bbR^n_{\geq 0}$ & $\Set{x \in \bbR^n | x \geq 0 }$ & Nonnegative orthant in $n$-dimensional space \\
$\Delta^V$ & $\Set{\lambda \in \bbR^V_{\geq 0} | \sum_{v \in V} \lambda_v = 1 }$ & Unit simplex on ground set $V$ \\
$\supp(\lambda)$ & $\Set{v \in V | \lambda_v \neq 0 }$ & Nonzero values (\emph{support}) of $\lambda$ \\
$P(T)$ & $\Set{\lambda \in \Delta^V | \supp(\lambda) \subseteq T }$ & Face of the unit simplex given by components $T$ \\
$\ext(P)$ & - & Extreme points of polyhedra $P$ \\
$\gr(f)$ & $\Set{(x,f(x)) | x \in \dom(f) }$ & Graph of the function $f$ \\
$[V]^2$ & $\Set{\{u,v\} \in V \times V | u \neq v }$ & All unordered pairs of elements in $V$ \\
$\Em(\calT,H)$ & $\bigcup_{i=1}^d P(T^i) \times \{H_i\}$ &  Embedding of disjunctive constraint (where $H_i$ is the $i$-th row of $H$) \\
$\Conv(S)$ & - & Convex hull of $S$ \\
$Q(\calT,H)$ & $\Conv(\Em(\calT,H))$ & Convex hull of embedding \\
$\aff(H)$ & - & Affine hull of the rows of $H$ \\
$L(H)$ & $\Set{y-H_1 | y \in \aff(H) }$ & Linear space parallel to the affine hull $\aff(H)$ (where $H_1$ is first row of $H$) \\
$M(b)$ & $\Set{y \in L(H) | b \cdot y = 0 }$ & The hyperplane in $L(H)$ normal to $b$ \\
$\Vol(D)$ & - & Volume of set $D$ \\
$A * B$ & $\Set{\{u,v\} | u \in A, v \in B }$ & Unordered pairs of elements in $A$ and $B$
\end{tabular}
\caption{Notation used throughout the paper.}
\label{tab:notation}
}
\end{center}
\end{table}

In order to solve an optimization problem containing $f$, we will construct a formulation for its graph $\gr(f) \defeq \Set{(x,f(x)) | x \in D }$, which will couple the argument $x$ with the function output $f(x)$. We can view the graph disjunctively as the union $\gr(f) = \bigcup_{i=1}^d S^i$, where each $S^i = \Set{ (x,f^i(x)) | x \in C^i }$ is a \emph{segment} of the graph.

\begin{example}\label{startex}
  Consider the univariate piecewise linear function $f : [1,5] \to \bbR$ with the domain pieces $C^1=[1,2]$, $C^2=[2,3]$, $C^3=[3,4]$, and $C^4=[4,5]$, where
  \begin{subequations} \label{eqn:pwl-example}
  \begin{alignat}{3}
      x \in C^1 &\implies f(x) = 4x-4,   \quad\quad &x \in C^2 \implies& f(x) = 3x-2, \\
      x \in C^3 &\implies f(x) = 2x+1, \quad\quad &x \in C^4 \implies& f(x) = x+5.
  \end{alignat}
  \end{subequations}
The graph of the piecewise linear function is then
    \[
        \gr(f) = \Set{(x,4x-4) | x \in C^1 } \cup \Set{(x,3x-2) | x \in C^2 } \cup \Set{(x,2x+1) | x \in C^3 } \cup \Set{(x,x+5) | x \in C^4}.
    \]
    Similarly, we take the bivariate piecewise linear function $g : [0,1]^2 \to \bbR$ with the domain partition $C^1 = \Set{x \in [0,1]^2 | x_1 \leq x_2 }$ and $C^2 = \Set{x \in [0,1]^2 | x_1 \geq x_2 }$, and
  	\begin{equation}\label{eqn:bivariate-pwl-example}
      x \in C^1 \implies g(x) = -x_1+3x_2+1, \quad\quad x \in C^2 \implies g(x) = x_1+x_2+1.
  	\end{equation}
	The corresponding graph is
  \[
  	  \gr(g) = \Set{(x_1,x_2,-x_1+3x_2+1) | x \in C^1 } \cup \Set{(x_1,x_2,x_1+x_2+1) | x \in C^2 }.
  \]
\end{example}

We refer the reader to \cite{Vielma:2010} for an exhaustive taxonomy of existing MIP formulations for piecewise linear functions. In this work, we will build formulations for piecewise linear functions using the \emph{combinatorial disjunctive constraint} approach~\citep{Huchette:2016a}.

Given a piecewise linear function, take the family of sets $\calT = (T^i = \ext(C^i))_{i=1}^d$ corresponding to the extreme points of each piece of the domain $C^i$. This describes the underlying combinatorial structure among the segments of of the graph, induced by the shared breakpoints over the ground set $V = \bigcup_{i=1}^d T^i$. Define $\Delta^V \defeq \Set{\lambda \in \bbR^V_{\geq 0} | \sum_{v \in V} \lambda_v = 1 }$ as the standard simplex, $\supp(\lambda) \defeq \Set{v \in V | \lambda_v \neq 0 }$ as the nonzero values (\emph{support}) of $\lambda$, and $P(T) \defeq \Set{\lambda \in \Delta^V | \supp(\lambda) \subseteq T }$ as the face of the standard simplex with support restricted to $T$. Then we can express the graph in terms of $\calT$ as $\gr(f) = \Set{\sum_{v \in V} \lambda_v (v,f(v)) | \lambda \in \bigcup_{i=1}^d P(T^i) }$. In particular, we can build a formulation for $f$ through the \emph{combinatorial disjunctive constraint} $\lambda \in \bigcup_{i=1}^d P(T^i)$ \citep{Huchette:2016a}, which is a disjunctive constraint on convex multipliers $\lambda$ where each alternative $P(T^i)$ is some face of the unit simplex $\Delta^V$. 

\begin{example} \label{example-1}
	Take $f$ as given in Example~\ref{startex}. The graph of this function has $d=4$ segments, and the breakpoints between segments are given by the set $V=\llbracket d+1 \rrbracket$. We have that $(x,z) \in \gr(f)$ if and only if $(x,z) = \sum_{v \in V} (v,f(v))\lambda_v$ for some $\lambda\in \bigcup_{i=1}^4 P(\{i,i+1\})$.

	Similarly, for the function $g$ as in Example~\ref{startex}, we can take $V = \{0,1\}^2$ and observe that $(x,z) \in \gr(g)$ if and only if $(x,z) = \sum_{v \in V} (v,f(v)) \lambda_v$ for some $\lambda \in P(\{(0,0),(1,0),(1,1)\}) \cup P(\{(0,0),(0,1),(1,1)\})$.
\end{example}


For the remainder, we assume without loss of generality (w.l.o.g.) that $V=\llbracket d+1\rrbracket$ for univariate functions. We also note that the constraint $\lambda\in \bigcup_{i=1}^{d+1} P(\{i,i+1\})$ from Example~\ref{example-1} is the classical special ordered set of type 2 (SOS2) constraint~\citep{Beale:1970}. \blue{Additionally, we will restrict our attention to bivariate functions with \emph{grid triangulated} domains; that is, the domain of the function is a rectangle in $\bbR^2$ decomposed according to a regular grid, and that each subrectangle within that decomposition is partitioned into exactly two triangular domain pieces. In this setting, we may assume w.l.o.g. and that $V=\llbracket d_1+1 \rrbracket \times \llbracket d_2+1 \rrbracket$ for bivariate functions. 
}

The logarithmic formulations of \cite{Vielma:2009a} apply to several special classes of combinatorial disjunctive constraints, including SOS2 constraints. These logarithmic formulations have been observed to perform extremely well computationally; this can be largely attributed to their strength and size. With regards to strength, the formulations are \emph{ideal}: their LP relaxations offer the tightest possible convex relaxation for the underlying nonconvex set $\gr(f)$, and their extreme points naturally satisfy the desired integrality condition (see \cite{Vielma:2015} for more about formulation strength). Moreover, the formulations are small, as the number of auxiliary variables and (general inequality) constraints scale logarithmically in the number of segments of the piecewise linear functions. The novel formulations presented in this work will also possess these two properties. Moreover, we will also design them to have other desirable computational properties (univariate functions in Section~\ref{sec:univariate}), and such that they apply to a much larger class of piecewise linear functions than previously considered (bivariate functions in Section~\ref{sec:bivariate}).  To achieve this, we use and extend two recent generalizations of \cite{Vielma:2009a}: the geometric \emph{embedding formulation} technique of \cite{Vielma:2016}, and the combinatorial \emph{independent branching formulation} technique of  \cite{Huchette:2016a}. When applied to the SOS2 constraint, these techniques yield two formulations we denote the \emph{logarithmic embedding} (\Log{}) and the \emph{logarithmic independent branching} (\LogIB{}) formulations, respectively. Both formulations are quite similar; however, while \LogIB{} always exactly coincides with the original logarithmic formulation of \cite{Vielma:2009a}, \Log{} only does so when $d$ is a power-of-two \citep{Muldoon:2012} (we provide an example of this divergence in Appendix~\ref{app:log-logib}). These two formulations will serve as the reference benchmark formulation in our computational experiments.


%
%
%

\section{Formulations for univariate piecewise linear functions} \label{sec:univariate}

In this section we will adapt a geometric formulation construction method to build novel strong logarithmic formulations for univariate piecewise linear functions.

\subsection{The embedding approach} \label{ssec:embedding-approach}
The embedding approach of \cite{Vielma:2016} provides one way to construct strong formulations for disjunctive constraints. To formulate $\bigcup_{i=1}^d P(T^i)$, assign each alternative $P(T^i)$ a unique integer \emph{code} $H_i \in \bbZ^r$. We call the collection of all codes as rows in a matrix $H \in \bbZ^{d \times r}$ an \emph{encoding}, where $H_i$ is the $i$-th row of $H$. Then the disjunctive set is ``embedded'' in a higher-dimensional space as $\Em(\calT,H) \defeq \bigcup_{i=1}^d (P(T^i) \times \{H_i\})$. In the case studied by \cite{Vielma:2016} where $H \in \{0,1\}^{d \times r}$ is a binary encoding, this easily leads to a MIP formulation for $\bigcup_{i=1}^d P(T^i)$. However, we will be interested in constructing formulations using general integer encodings, which requires some care to ensure that the embedding leads to a valid formulation.

\begin{definition}
Take the matrix $H \in \bbZ^{d \times r}$, and the collection of its rows as $\Lambda = \{H_i\}_{i=1}^d$.
\begin{itemize}
	\item $H$ is in \emph{convex position} if $\ext(\Conv(\Lambda)) = \Lambda$.
	\item $H$ is \emph{hole-free} if $\Conv(\Lambda) \cap \bbZ^r = \Lambda$.
\end{itemize}
Take $\mathcal{H}_r(d)\defeq \Set{ H\in \mathbb{Z}^{d\times r} | H \text{ is hole-free and in convex position, and each } H_i \text{ is distinct}}$.
\end{definition}

The following straightforward extension of Proposition 1 and Corollary 1 in \cite{Vielma:2016} shows that encodings in $\mathcal{H}_r(d)$ always lead to valid formulations. 
\begin{proposition}\label{firstprop}
Take the family of sets $\calT = (T^i \subseteq V)_{i=1}^d$, along with $r\geq \left\lceil \log_2 (d)\right\rceil$ and $H \in \mathcal{H}_r(d)$. Then $Q(\calT,H) \defeq \Conv(\Em(\calT,H))$ is a rational polyhedron, and an ideal formulation for $\bigcup_{i=1}^d P(T^i)$ is $\Set{(\lambda,y) \in Q(\calT,H) | y \in \bbZ^r }$. We call this the \emph{embedding formulation} of $\calT$ associated to $H$.
\end{proposition}


In general, constructing a linear inequality description of $Q(\calT,H)$ is difficult, the resulting representation may be exponentially large, and its structure is highly dependent on the interplay between the sets $\calT$ and the encoding $H$.
Fortunately, \cite[Proposition 2]{Vielma:2016} gives an explicit description of $Q(\calT,H)$ for the SOS2 constraint with any choice of binary encoding $H$. This description is geometric, in terms of the difference directions $H_{i+1}-H_i$ between adjacent codes. In particular, we will need to compute all hyperplanes $M(b) \defeq \Set{ y \in L(H) | b \cdot y = 0}$ spanned by these difference directions in $L(H) \defeq \Set{y - H_1 | y \in \aff(H) }$, the linear space parallel to the affine hull of $H$. The following straightforward extension of Proposition 2 from \cite{Vielma:2016} shows that this description  also holds for any encoding in  $\mathcal{H}_r(d)$.

\begin{proposition} \label{prop:sos2}
Take $H \in \mathcal{H}_r(d)$, along with $H_0 \equiv H_1$ and $H_{d+1} \equiv H_d$ for notational convenience. Let $\mathcal{B} \subset L(H) \backslash \{{\bf 0}^r\}$ be normal directions such that $\{M(b)\}_{b \in \mathcal{B}}$ is the set of hyperplanes spanned by $\{H_{i+1}-H_i\}_{i=1}^{d-1}$ in $L(H)$. If $\calT = (\{i,i+1\})_{i=1}^d$ is the family of sets defining the SOS2 constraint on $d+1$ breakpoints, then $Q(\calT,H)$ is equal to all $(\lambda,y) \in \Delta^{d+1} \times \aff(H)$ such that
\begin{equation*}
        \sum\nolimits_{v=1}^{d+1} \min\{b \cdot H_{v-1}, b \cdot H_{v}\}\lambda_v \leq b \cdot y \leq \sum\nolimits_{v=1}^{d+1} \max\{b \cdot H_{v-1}, b \cdot H_{v}\}\lambda_v \quad \forall b \in \mathcal{B}.
    \end{equation*}
\end{proposition}

Consider the class of encodings $K^r \in \mathcal{H}_r(d)$ for $r = \lceil \log_2(d) \rceil$ known as Gray codes \citep{Savage:1997}, where adjacent codes differ in exactly one component (i.e. $\norm{K^r_{j+1}-K^r_j}_1=1$ for all $j\in \sidx{d-1}$). This class of encodings enjoys the desirable property that the spanning hyperplanes needed for Proposition~\ref{prop:sos2} are parsimonious and simple to describe in closed form. For the remainder, we will work with a particular Gray code known as the binary reflected Gray code (BRGC); see Appendix~\ref{app:encodings} for a formal definition.

We can construct the \emph{logarithmic embedding} (\Log{}) formulation for the SOS2 constraint due to \cite{Vielma:2016} by applying Proposition~\ref{prop:sos2} with the BRGC. This formulation is ideal, and its size scales logarithmically in the number of segments $d$. 

\begin{example}\label{example:sos2}
The \Log{} formulation for the SOS2 constraint with $d=4$ (arising in \eqref{eqn:pwl-example}) is:
\begin{equation}\label{eqn:log-example}
    \lambda_3 \leq y_1, \quad\quad \lambda_1 + \lambda_5 \leq 1 - y_1, \quad\quad \lambda_4 + \lambda_5 \leq y_2, \quad\quad \lambda_1 + \lambda_2 \leq 1-y_2, \quad\quad (\lambda,y) \in \Delta^{V} \times \{0,1\}^2.
\end{equation}
\end{example}

\subsection{Branching behavior of existing formulations}

As observed by \cite{Vielma:2010} and in our computational experiments, logarithmic formulations such as \Log{} can offer a considerable computational advantage over other approaches, particularly for univariate piecewise linear functions with many segments (i.e. large $d$). However, it has also been observed that variable branching with logarithmic formulations such as \Log{} can produce weak dual bounds (e.g. \cite{Martin:2006,Rebennack:2016,Yildiz:2013}). 

To quantitatively assess relaxation strength after branching, we consider two metrics. The first is the volume of the projection of the LP relaxation onto $(x,z)$-space (cf. \citep{Lee:2018} for a recent work using volume as a metric for formulation quality). The second is the proportion of the function domain where the LP relaxation after branching is stronger than the LP relaxation before branching. More formally, if $D$ is the domain of $f$, $F$ is the projection of the original LP relaxation onto $(x,z)$-space, and $F'$ is the same projection of the LP relaxation after branching, then we report $\frac{1}{\Vol(D)}\Vol\left(\Set{x \in D | \min_{(x,z) \in F} z < \min_{(x,z) \in F'} z }\right)$, which we dub the \emph{strengthened proportion}.

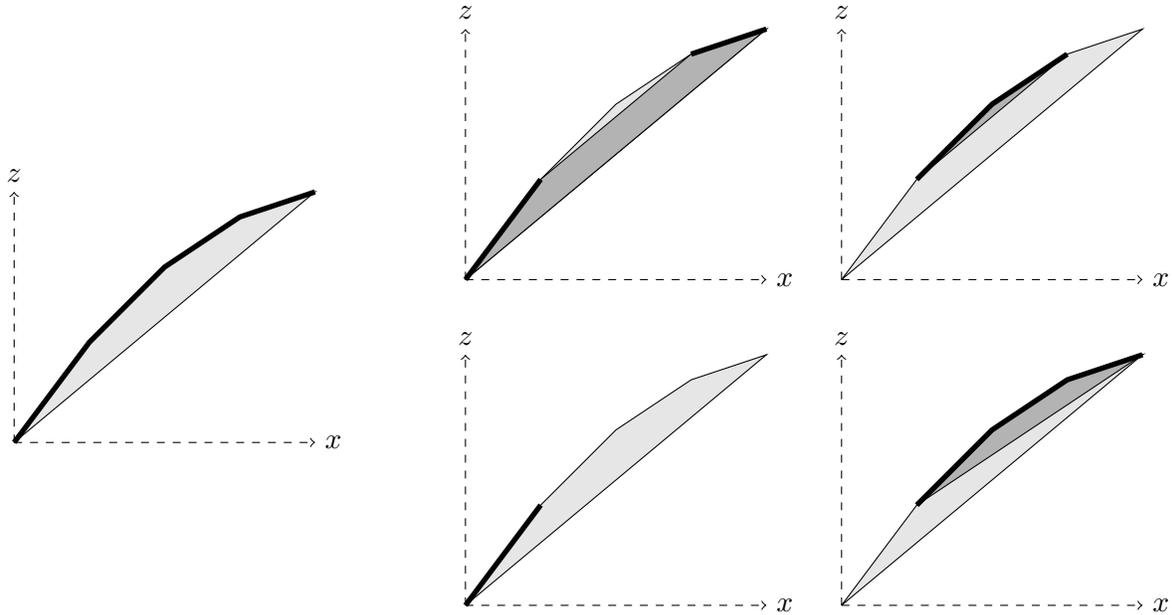
\begin{figure}[htpb]
\centering
\begin{tikzpicture}
\begin{scope}[shift={(0,-13/6)}]
    \draw [dashed,->] (0,0) -- (4,0);
    \draw [dashed,->] (0,0) -- (0,10/3);
    \node [right] at (4,0) {$x$};
    \node [above] at (0,10/3) {$z$};

    \draw [fill=gray!20] (0,0/3) -- (1,4/3) -- (2,7/3) -- (3,9/3) -- (4,10/3) -- (0,0);
    \draw [line width=2] (0,0/3) -- (1,4/3) -- (2,7/3) -- (3,9/3) -- (4,10/3);
\end{scope}
\begin{scope}[shift={(6,0)}]
    \draw [dashed,->] (0,0) -- (4,0);
    \draw [dashed,->] (0,0) -- (0,10/3);
    \node [right] at (4,0) {$x$};
    \node [above] at (0,10/3) {$z$};

    \draw [fill=gray!20] (0,0/3) -- (1,4/3) -- (2,7/3) -- (3,9/3) -- (4,10/3) -- (0,0);
    \draw [fill=gray!60] (0,0/3) -- (1,4/3) -- (3,9/3) -- (4,10/3) -- (0,0);
    \draw [line width=2] (0,0/3) -- (1,4/3);
    \draw [line width=2] (3,9/3) -- (4,10/3);
\end{scope}
\begin{scope}[shift={(11,0)}]
    \draw [dashed,->] (0,0) -- (4,0);
    \draw [dashed,->] (0,0) -- (0,10/3);
    \node [right] at (4,0) {$x$};
    \node [above] at (0,10/3) {$z$};

    \draw [fill=gray!20] (0,0/3) -- (1,4/3) -- (2,7/3) -- (3,9/3) -- (4,10/3) -- (0,0);
    \draw [fill=gray!60] (1,4/3) -- (2,7/3) -- (3,9/3) -- (1,4/3);
    \draw [line width=2] (1,4/3) -- (2,7/3) -- (3,9/3);
\end{scope}
\begin{scope}[shift={(6,-13/3)}]
    \draw [dashed,->] (0,0) -- (4,0);
    \draw [dashed,->] (0,0) -- (0,10/3);
    \node [right] at (4,0) {$x$};
    \node [above] at (0,10/3) {$z$};

    \draw [fill=gray!20] (0,0/3) -- (1,4/3) -- (2,7/3) -- (3,9/3) -- (4,10/3) -- (0,0);
    \draw [fill=gray!60,line width=2] (0,0) -- (1,4/3);
\end{scope}
\begin{scope}[shift={(11,-13/3)}]
    \draw [dashed,->] (0,0) -- (4,0);
    \draw [dashed,->] (0,0) -- (0,10/3);
    \node [right] at (4,0) {$x$};
    \node [above] at (0,10/3) {$z$};

    \draw [fill=gray!20] (0,0/3) -- (1,4/3) -- (2,7/3) -- (3,9/3) -- (4,10/3) -- (0,0);
    \draw [fill=gray!60] (1,4/3) -- (2,7/3) -- (3,9/3) -- (4,10/3) -- (1,4/3);
    \draw [line width=2] (1,4/3) -- (2,7/3) -- (3,9/3) -- (4,10/3);
\end{scope}
\end{tikzpicture}
    \caption{\textbf{(Left)} The LP relaxation of an ideal formulation (e.g. \Log{}) \eqref{eqn:log-example}} projected onto $(x,z)$-space. The \Log{} formulation after \textbf{(top center)} down-branching $y_1 \leq 0$ and \textbf{(top right)} up-branching $y_1 \geq 1$. The \Inc{} formulation after \textbf{(bottom center)} down-branching $y_1 \leq 0$ and \textbf{(bottom right)} up-branching $y_1 \geq 1$.
    \label{fig:log-branching}
\end{figure}

We turn to the \Log{} formulation for $d=4$ given in Example~\ref{eqn:log-example}. The mapping from the $\lambda$ variables to the original space is $(x,z) = (0,0)\lambda_1 + (1,4)\lambda_2 + (2,7) \lambda_3 + (3,9) \lambda_4 + (4,10) \lambda_5$. Qualitatively, in the top row of Figure~\ref{fig:log-branching} we see that the LP relaxation, projected onto the $(x,z)$-space of the graph $\gr(f)$, remains largely unchanged in the down-branching subproblem (i.e. when we branch $y_1 \leq 0$). This is undesirable, as it will not improve dual bounds in a branch-and-bound setting, which is crucial to ensuring fast convergence. Quantitatively, the strengthened proportion for this down-branch is 0, and so when minimizing $f$, the dual bound will be the same after branching as for the original LP relaxation (assuming both are feasible). From this, it is reasonable to infer that the high-performance of the \Log{} formulation is due to its strength and small size, and in spite of its poor branching behavior.

In contrast, the traditional SOS2 constraint branching~\citep{Beale:1970} induces much more balanced branches in $(x,z)$-space. Additionally, the incremental \Inc{} formulation~\citep{Dantzig:1960,Padberg:2000,Croxton:2003} is a MIP formulation that induces the same branching behavior, which we depict in the bottom row of Figure~\ref{fig:log-branching}. Quantitatively, both up- and down-branching result in subproblems with a strengthened proportion of 1, meaning that the formulation will always lead to strictly stronger dual bounds when minimizing $f$. In particular, we highlight the \emph{incremental branching} behavior of the \Inc{} formulation in the $(x,z)$-space: after selecting for branching a binary variable $y_k$ (\Inc{} has $d-1$ binary variables, so $k \in \llbracket d-1 \rrbracket$), the only points $(x,z) \in \gr(f)$ feasible for the down-branch (resp. up-branch) are those that lie on segments 1 to $k-1$, $\bigcup_{i=1}^{k-1} S^i$ (resp. $k$ to $d$, $\bigcup_{i=k}^d S^i$). Additionally, the \Inc{} formulation is \emph{hereditarily sharp}: each subproblem LP relaxation projects to exactly the convex hull (either $\Conv(\bigcup_{i=1}^{k-1} S^i)$ or $\Conv(\bigcup_{i=k}^d S^i)$) of the segments feasible for that subproblem~\citep{Jeroslow:1984,Jeroslow:1988c}. This combination has been observed to lead to very balanced branch-and-bound trees \citep{Yildiz:2013,Vielma:2015}, and the \Inc{} formulation has been observed to perform very well for small $d$, before its size (which scales linearly in $d$) becomes overwhelming (see the computational results in Section~\ref{ss:univariate-computational-results}).

\subsection{New zig-zag formulations for the SOS2 constraint}

We now present new embedding formulations for the SOS2 constraint that retain the size and strength of the \Log{} formulation, while repairing its degenerate branching behavior. For the remainder of the subsection, assume without loss of generality (w.l.o.g.) that $d$ is a power-of-two. Otherwise, construct the formulation for $\bar{d} = 2^{\lceil \log_2(d) \rceil}$ and fix the extraneous $\lambda_v$ variables to zero.

Take $K^r \in \mathcal{H}_r(d)$ as the BRGC for $d=2^r$ elements. Our first new encoding is the transformation of $K^r\in \{0,1\}^{d \times n}$ to $C^r \in \bbZ^{d \times r}$, where $C^r_{i,k} = \sum_{j=2}^i \left| K^r_{j,k}-K^r_{j-1,k} \right|$ for each $i \in \llbracket d \rrbracket$ and $k \in \llbracket r \rrbracket$. In words, $C^r_{i,k}$ is the number of times the sequence $(K^r_{1,k},\ldots,K^r_{i,k})$ changes value, and is monotonic nondecreasing in $i$. Our second encoding will be $Z^r \in \{0,1\}^{d \times r}$ with $Z^r_{i} = \scrA(C^r_i)$ for each $i \in \llbracket d \rrbracket$, where $\scrA : \bbR^r \to \bbR^r$ is the linear map  given by $\scrA(y)_k = y_k - \sum_{\ell=k+1}^r y_\ell$ for each component $k \in \llbracket r \rrbracket$. We show the encodings for $r=3$ in Figure~\ref{fig:encodings}, and include formal recursive definitions for them in Appendix~\ref{app:encodings}, where we additionally show that $C^r, Z^r \in \mathcal{H}_r(d)$. Applying Proposition~\ref{prop:sos2} with the new encodings gives two new small, strong formulations for the SOS2 constraint.
\begin{proposition} \label{prop:new-sos2}
Take $r = \lceil \log_2(d) \rceil$, along with $C^r_0 \equiv C^r_1$ and $C^r_{d+1} \equiv C^r_d$ for notational simplicity. Then two ideal formulations for the SOS2 constraint with $d$ segments are given by
  \begin{equation} \label{eqn:integer-zigzag-formulation}
     \sum\nolimits_{v=1}^{d+1} C^r_{v-1,k}\lambda_v \leq y_k \leq \sum\nolimits_{v=1}^{d+1} C^r_{v,k}\lambda_v \quad \forall k \in \llbracket r \rrbracket,\quad\quad (\lambda,y) \in \Delta^{d+1} \times \bbZ^r
    \end{equation}
    and
        \begin{equation} \label{eqn:binary-zigzag-formulation}
     \sum\nolimits_{v=1}^{d+1} C^r_{v-1,k}\lambda_v \leq y_k + \sum\nolimits_{\ell=k+1}^r 2^{\ell-k-1} y_{\ell} \leq
        \sum\nolimits_{v=1}^{d+1} C^r_{v,k}\lambda_v \quad \forall k \in \llbracket r \rrbracket,\quad\quad (\lambda,y) \in \Delta^{d+1} \times \{0,1\}^r.
    \end{equation}
\end{proposition}

We dub \eqref{eqn:binary-zigzag-formulation} the \emph{binary zig-zag} (\ZZB{}) formulation for the SOS2 constraint, as its associated binary encoding $Z^r$ ``zig-zags'' through the interior of the unit hypercube (See Figure~\ref{fig:encodings}). We will refer to formulation~\eqref{eqn:integer-zigzag-formulation} as the \emph{general integer zig-zag} (\ZZI{}) formulation because of its use of general integer encoding $C^r \in \bbZ^{d \times r}$. We emphasize that \ZZI{} and \ZZB{} are logarithmically-sized in $d$ and ideal: the same size and strength as the existing \Log{} formulation.

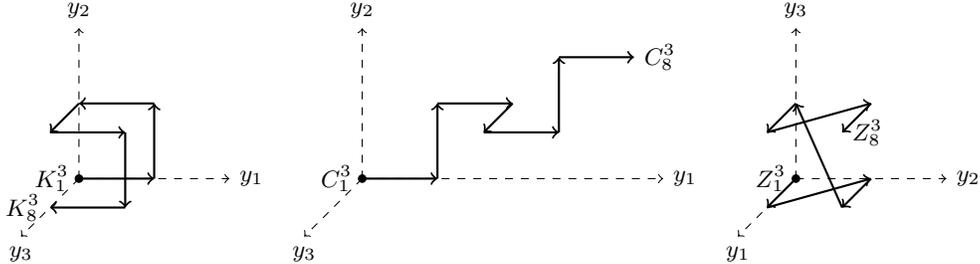
\begin{figure}[htpb]
    \centering
    \begin{tikzpicture}
        \draw [thick, ->] (0,0,0) -- (1,0,0);
        \draw [thick, ->] (1,0,0) -- (1,1,0);
        \draw [thick, ->] (1,1,0) -- (0,1,0);
        \draw [thick, ->] (0,1,0) -- (0,1,1);
        \draw [thick, ->] (0,1,1) -- (1,1,1);
        \draw [thick, ->] (1,1,1) -- (1,0,1);
        \draw [thick, ->] (1,0,1) -- (0,0,1);
        \begin{scope}[canvas is xy plane at z=0]
            [ultra thick] \node [right] at (2,0) {{\smaller $y_1$}};
            [ultra thick] \node [left] at (0,0) {{\smaller $K^3_1$}};
        \end{scope}
        \begin{scope}[canvas is xy plane at z=1]
            \node [left] at (0,0) {\smaller $K^3_8$};
        \end{scope}
        \begin{scope}[canvas is xz plane at y=0]
            [ultra thick] \node [below] at (0,2) {{\smaller $y_3$}};
        \end{scope}
        \begin{scope}[canvas is yz plane at x=0]
            [ultra thick] \node [above] at (2,0) {{\smaller $y_2$}};
        \end{scope}
        \begin{scope}[canvas is xy plane at z=0]
            \draw [fill] (0,0) circle [radius=.05];
        \end{scope}
        \draw [dashed, ->] (0,0,0) -- (2,0,0);
        \draw [dashed, ->] (0,0,0) -- (0,2,0);
        \draw [dashed, ->] (0,0,0) -- (0,0,2);
    \end{tikzpicture}
    \begin{tikzpicture}
        \draw [thick, ->] (0,0,0) -- (1,0,0);
        \draw [thick, ->] (1,0,0) -- (1,1,0);
        \draw [thick, ->] (1,1,0) -- (2,1,0);
        \draw [thick, ->] (2,1,0) -- (2,1,1);
        \draw [thick, ->] (2,1,1) -- (3,1,1);
        \draw [thick, ->] (3,1,1) -- (3,2,1);
        \draw [thick, ->] (3,2,1) -- (4,2,1);
        \begin{scope}[canvas is xy plane at z=0]
            [ultra thick] \node [right] at (4,0) {{\smaller $y_1$}};
            [ultra thick] \node [left] at (0,0) {{\smaller $C^3_1$}};
        \end{scope}
        \begin{scope}[canvas is xy plane at z=1]
            \node [right] at (4,2) {\smaller $C^3_8$};
        \end{scope}
        \begin{scope}[canvas is xz plane at y=0]
            [ultra thick] \node [below] at (0,2) {{\smaller $y_3$}};
        \end{scope}
        \begin{scope}[canvas is yz plane at x=0]
            [ultra thick] \node [above] at (2,0) {{\smaller $y_2$}};
        \end{scope}
        \begin{scope}[canvas is xy plane at z=0]
            \draw [fill] (0,0) circle [radius=.05];
        \end{scope}
        \draw [dashed, ->] (0,0,0) -- (4,0,0);
        \draw [dashed, ->] (0,0,0) -- (0,2,0);
        \draw [dashed, ->] (0,0,0) -- (0,0,2);
    \end{tikzpicture}
    \begin{tikzpicture}
        \draw [thick, ->] (0,0,0) -- (0,0,1);
        \draw [thick, ->] (0,0,1) -- (1,0,0);
        \draw [thick, ->] (1,0,0) -- (1,0,1);
        \draw [thick, ->] (1,0,1) -- (0,1,0);
        \draw [thick, ->] (0,1,0) -- (0,1,1);
        \draw [thick, ->] (0,1,1) -- (1,1,0);
        \draw [thick, ->] (1,1,0) -- (1,1,1);
        \begin{scope}[canvas is xy plane at z=0]
            [ultra thick] \node [right] at (2,0) {{\smaller $y_2$}};
            \node [left] at (0,0) {\smaller $Z^3_1$};
        \end{scope}
        \begin{scope}[canvas is xy plane at z=1]
            \node [right] at (1,1) {\smaller $Z^3_8$};
        \end{scope}
        \begin{scope}[canvas is xz plane at y=0]
            [ultra thick] \node [below] at (0,2) {{\smaller $y_1$}};
        \end{scope}
        \begin{scope}[canvas is yz plane at x=0]
            [ultra thick] \node [above] at (2,0) {{\smaller $y_3$}};
        \end{scope}
        \begin{scope}[canvas is xy plane at z=0]
            \draw [fill] (0,0) circle [radius=.05];
        \end{scope}
        \draw [dashed, ->] (0,0,0) -- (2,0,0);
        \draw [dashed, ->] (0,0,0) -- (0,2,0);
        \draw [dashed, ->] (0,0,0) -- (0,0,2);
    \end{tikzpicture}
    \caption{Depiction of $K^3$ (\textbf{Left}), $C^3$ (\textbf{Center}), and $Z^3$ (\textbf{Right}). The first row of each is marked with a dot, and the subsequent rows follow along the arrows. The axis orientation is different for $Z^3$ for visual clarity.}
    \label{fig:encodings}
\end{figure}



To study the branching behavior of the \ZZI{} formulation, we return to the SOS2 constraint with $d=4$ from Example~\ref{example:sos2}. The formulation consists of all $(\lambda,y) \in \Delta^5 \times \bbZ^2$ such that
\begin{equation} \label{eqn:zigzag-example}
    \lambda_3 + \lambda_4 + 2 \lambda_5 \leq y_1 \leq \lambda_2 + \lambda_3 + 2\lambda_4 + 2\lambda_5, \quad\quad \lambda_4 + \lambda_5 \leq y_2 \leq \lambda_3 + \lambda_4 + \lambda_5.
\end{equation}

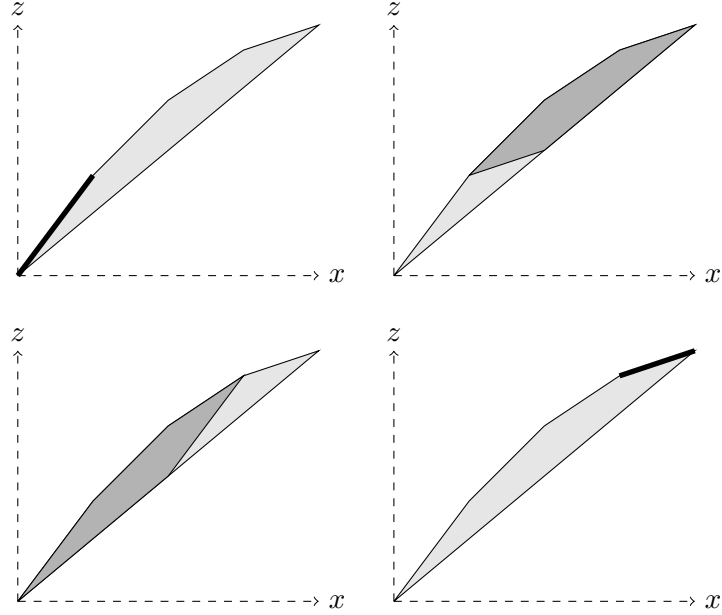
\begin{figure}[htpb]
    \centering
\begin{tikzpicture}
\begin{scope}[shift={(0,0)}]
    \draw [dashed,->] (0,0) -- (4,0);
    \draw [dashed,->] (0,0) -- (0,10/3);
    \node [right] at (4,0) {$x$};
    \node [above] at (0,10/3) {$z$};

    \draw [fill=gray!20] (0,0/3) -- (1,4/3) -- (2,7/3) -- (3,9/3) -- (4,10/3) -- (0,0);
    \draw [fill=gray!60, line width=2] (0,0/3) -- (1,4/3);
\end{scope}
\begin{scope}[shift={(5,0)}]
    \draw [dashed,->] (0,0) -- (4,0);
    \draw [dashed,->] (0,0) -- (0,10/3);
    \node [right] at (4,0) {$x$};
    \node [above] at (0,10/3) {$z$};

    \draw [fill=gray!20] (0,0/3) -- (1,4/3) -- (2,7/3) -- (3,9/3) -- (4,10/3) -- (0,0);
    \draw [fill=gray!60] (1,4/3) -- (2,7/3) -- (3,9/3) -- (4,10/3) -- (2,5/3) -- (1,4/3);
\end{scope}
\begin{scope}[shift={(0,-13/3)}]
    \draw [dashed,->] (0,0) -- (4,0);
    \draw [dashed,->] (0,0) -- (0,10/3);
    \node [right] at (4,0) {$x$};
    \node [above] at (0,10/3) {$z$};

    \draw [fill=gray!20] (0,0/3) -- (1,4/3) -- (2,7/3) -- (3,9/3) -- (4,10/3) -- (0,0);
    \draw [fill=gray!60] (0,0) -- (1,4/3) -- (2,7/3) -- (3,9/3) -- (2,5/3) -- (0,0);
\end{scope}
\begin{scope}[shift={(5,-13/3)}]
    \draw [dashed,->] (0,0) -- (4,0);
    \draw [dashed,->] (0,0) -- (0,10/3);
    \node [right] at (4,0) {$x$};
    \node [above] at (0,10/3) {$z$};

    \draw [fill=gray!20] (0,0/3) -- (1,4/3) -- (2,7/3) -- (3,9/3) -- (4,10/3) -- (0,0);
    \draw [fill=gray!60,line width=2] (3,9/3) -- (4,10/3);
\end{scope}
\end{tikzpicture}
    \caption{The LP relaxation of the \ZZI{} formulation \eqref{eqn:zigzag-example} projected onto $(x,z)$-space, after down-branching $y_1 \leq 0$ \textbf{(top center)}, up-branching $y_1 \geq 1$ \textbf{(bottom center)}, down-branching $y_1 \leq 1$ \textbf{(top right)}, and up-branching $y_1 \geq 2$ \textbf{(bottom right)}.}
    \label{fig:pwl-branching}
\end{figure}

We have two possibilities for branching on $y_1$, depicted in Figure~\ref{fig:pwl-branching}: down on $y_1 \leq 0$ and up on $y_1 \geq 1$, or down on $y_1 \leq 1$ and up on $y_1 \geq 2$. We note that after imposing either $y_1 \leq 0$ or $y_1 \geq 2$, the relaxation is then exact, i.e. the relaxation is equal to exactly one of the segments of the graph of $f$. Furthermore, when imposing either $y_1 \leq 1$ or $y_1 \geq 1$, we deduce a general inequality on the $\lambda$ variables that improves the strengthened proportion relative to \Log{}: either $\lambda_1 \leq \lambda_4 + \lambda_5$ or $\lambda_5 \leq \lambda_1 + \lambda_2$, respectively.

\begin{table}[htpb]
    \centering
    \smaller
    \begin{tabular}{c|c|cc|cc|cc:cc}
         Statistic & LP Relaxation & \Log{} $0\downarrow$ & \Log{} $1\uparrow$ & \Inc{} $0\downarrow$ & \Inc{} $1\uparrow$ & \ZZI{} $0\downarrow$ & \ZZI{} $1\uparrow$ & \ZZI{} $1\downarrow$ & \ZZI{} $2\uparrow$ \\ \hline
         Volume & 6 & 5.5 & 0.5 & 0 & 2 & 0 & 3.5 & 3.5 & 0 \\
         Strengthened Prop. & 0 & 0 & 1 & 1 & 1 & 1 & 0.5 & 0.5 & 1
    \end{tabular}
    \caption{Metrics for each possible branching decision on $z_1$ for \Log{}, \Inc{}, and \ZZI{} applied to \eqref{eqn:pwl-example}.}
    \label{tab:pwl-example-4}
\end{table}

As we see qualitatively in Figures~\ref{fig:log-branching} and \ref{fig:pwl-branching} and quantitatively in Table~\ref{tab:pwl-example-4}, the \ZZI{} formulation yields LP relaxations after branching that are stronger and more balanced than those of the \Log{} formulation. In Appendix~\ref{app:8-segment}, we offer a more complex example with an 8-segment concave piecewise linear function where this effect is even more pronounced. An instructive way to interpret the branching of \ZZI{} is that it emulates the SOS2 branching induced by the \Inc{} formulation. In particular, the \ZZI{} formulation also induces incremental branching, but has slightly weaker subproblem relaxations compared to the \Inc{} formulation as it does not maintain the hereditary sharpness property. In this way, the \ZZI{} formulation maintains the size and strength of the \Log{} formulation, while inducing branching behavior that is much closer to the \Inc{} formulation.

\subsection{Univariate computational experiments} \label{ss:univariate-computational-results}
To evaluate the new \ZZI{} and \ZZB{} formulations against the existing formulations for univariate piecewise linear functions, we reproduce a variant of the computational experiments of \cite{Vielma:2010}, with the addition of the \ZZB{} and \ZZI{} formulations. Although the \LogIB{} formulation outperformed the rest of the formulations considered in \cite{Vielma:2010}, it has also been observed that logarithmic formulations tends to suffer from a significant performance degradation when the number of segments $d$ of the piecewise linear functions is not a power-of-two~\citep{Vielma:2009a,Coppersmith:2005,Muldoon:2012,Muldoon:2013}. Therefore, we will focus on problems of this form in our computational experiments. This is precisely the setting in which \Log{} and \LogIB{} (which we will introduce more formally in Section~\ref{ssec:independent-branching}) are not equivalent, and so we include both variants in our experiments. Finally, we also include the previously mentioned \Inc{} formulation, the \MC{}, \CC{}, and \DLog{} formulations as described by \cite{Vielma:2010}, as well as the SOS2 native branching (\SOS{}) implementation of the corresponding MIP solver. 

We evaluate our formulations on single commodity transportation problems of the form
\begin{align*}
\min_{x \geq 0}\quad& \sum_{i \in S} \sum_{j \in D} f_{i,j}(x_{i,j}) \\
\text{s.t.}\quad& \sum_{i \in S} x_{i,j} = d_j \quad \forall j \in D, \quad \quad \sum_{j \in D} x_{i,j} = s_i \quad \forall i \in S,
\end{align*}
where we match supply from nodes $S$ with demand from nodes $D$, while minimizing the transportation costs given by the sum of continuous nondecreasing concave univariate piecewise linear functions $f_{i,j}$ for each arc pair in $S \times D$.

We perform a scaling analysis along two axes: the size of the network (i.e. the cardinality of $S$ and $D$), and the number of segments for each piecewise linear function $f_{i,j}$. Regarding the first axis, we study both \emph{small} networks ($|S|=|D|=10$) and \emph{large} networks ($|S|=|D|=20$). Regarding the second axis, we study families of instances where each piecewise linear function has $d \in \{6,13,28,59\}$ segments.

We use CPLEX v12.7.0 with the JuMP algebraic modeling library \citep{Dunning:2015a} in the Julia programming language \citep{Bezanson:2017} for all computational trials, here and for the remainder of this work. All such trials were performed on an Intel i7-3770 3.40GHz Linux workstation with 32GB of RAM. For each trial, we allow the solver to run for 30 minutes to prove optimality before timing out. For each formulation and each family ($d \in \{6,13,28,59\}$) of 100 instances, we report the average solve time, standard deviation in solve time, and the number of instances for which the formulation was either the fastest (Win), or was unable to prove to optimality in 30 minutes or less (Fail).

We start by studying the small network instances in Table~\ref{tab:univariate-pwl-non-power-of-two}.
We observe that the \Inc{} formulation is superior for smaller function instances (i.e. with functions with fewer segments). Additionally, the \Log{} and \LogIB{} formulations have similar performance on all families of instances. We observe that the new \ZZI{} and \ZZB{} formulations are the best performers for larger function instances, and one of the two is the fastest formulation for every instance in the largest function family with $d=59$. Additionally, \ZZI{} and \ZZB{} both offer roughly a 2x speed-up in average solve time over \Log{} and \LogIB{} for most families of instances ($d \in \{13,28,59\}$).


\begin{table}[htpb]
    \centering
    \smaller
    \begin{tabular}{r|c|rrrrrrr:rr}
$d$ & Metric & \texttt{MC} & \texttt{CC} & \texttt{SOS2} & \texttt{Inc} & \texttt{DLog} & \Log{} & \texttt{LogIB} & \texttt{ZZB} & \texttt{ZZI} \\ \hline
\multirow{4}{*}{6}
 & Mean (s)  & \textbf{0.6}  & 3.8  & 1.1  & \textbf{0.6}  & 1.1  & 1.4  & 2.6  & 1.1  & 0.9  \\
 & Std  & \textbf{0.3}  & 4.1  & 1.5  & \textbf{0.3}  & 1.0  & 1.2  & 2.4  & 0.9  & 0.5  \\
 & Win & 35 & 0 & 7 & \textbf{46} & 5 & 1 & 0 & 4 & 2  \\
 & Fail & \textbf{0} & \textbf{0} & \textbf{0} & \textbf{0} & \textbf{0} & \textbf{0} & \textbf{0} & \textbf{0} & \textbf{0}  \\
\hline
\multirow{4}{*}{13}
 & Mean (s)  & 3.0  & 71.2  & 4.5  & \textbf{1.7}  & 4.6  & 4.4  & 4.2  & 2.4  & 2.6  \\
 & Std  & 3.1  & 152.0  & 5.8  & \textbf{0.7}  & 3.5  & 3.4  & 3.0  & 1.8  & 1.7  \\
 & Win & 11 & 0 & 9 & \textbf{47} & 11 & 0 & 0 & 15 & 7  \\
 & Fail & \textbf{0} & \textbf{0} & \textbf{0} & \textbf{0} & \textbf{0} & \textbf{0} & \textbf{0} & \textbf{0} & \textbf{0}  \\
\hline
\multirow{4}{*}{28}
 & Mean (s)  & 18.4  & 178.9  & 87.4  & 5.5  & 11.1  & 8.8  & 8.9  & 5.1  & \textbf{4.6}  \\
 & Std  & 26.0  & 359.3  & 309.3  & 4.4  & 8.1  & 5.6  & 5.4  & 3.7  & \textbf{2.7}  \\
 & Win & 1 & 0 & 6 & 14 & 1 & 0 & 0 & 37 & \textbf{41}  \\
 & Fail & \textbf{0} & 3 & 3 & \textbf{0} & \textbf{0} & \textbf{0} & \textbf{0} & \textbf{0} & \textbf{0}  \\
\hline
\multirow{4}{*}{59}
 & Mean (s)  & 348.7  & 541.0  & 664.3  & 17.1  & 19.1  & 16.3  & 16.0  & 9.8  & \textbf{9.3}  \\
 & Std  & 523.7  & 610.3  & 746.4  & 14.9  & 11.3  & 10.3  & 9.3  & 6.1  & \textbf{5.0}  \\
 & Win & 0 & 0 & 0 & 0 & 0 & 0 & 0 & 41 & \textbf{59}  \\
 & Fail & 7 & 13 & 26 & \textbf{0} & \textbf{0} & \textbf{0} & \textbf{0} & \textbf{0} & \textbf{0}
    \end{tabular}
    \caption{Computational results \blue{with CPLEX} for univariate transportation problems on small networks.}
    \label{tab:univariate-pwl-non-power-of-two}
\vspace{1em}
    \begin{tabular}{r|c|rrrrrrr:rr}
$d$ & Metric & \texttt{MC} & \texttt{CC} & \texttt{SOS2} & \texttt{Inc} & \texttt{DLog} & \Log{} & \texttt{LogIB} & \texttt{ZZB} & \texttt{ZZI} \\ \hline
\multirow{4}{*}{28}
 & Mean (s)  & 828.0  & 1769.3  & 1498.6  & 196.9  & 242.1  & 332.9  & 295.8  & 147.4  & \textbf{98.0}  \\
 & Std  & 714.3  & 211.5  & 646.9  & 206.8  & 282.2  & 430.4  & 387.9  & 228.2  & 144.4  \\
 & Win & 0 & 0 & 11 & 6 & 1 & 1 & 5 & 10 & \textbf{66}  \\
 & Fail & 28 & 97 & 80 & \textbf{0} & 1 & 2 & 2 & 1 & \textbf{0}  \\
\hline
\multirow{4}{*}{59}
 & Mean (s)  & 1596.9  & 1800.0  & 1800.0  & 793.4  & 777.1  & 749.3  & 753.5  & 328.7  & \textbf{273.1}  \\
 & Std  & 475.7  & - & - & 557.7  & 593.5  & 593.3  & 591.3  & 383.0  & 341.6  \\
 & Win & 0 & 0 & 0 & 2 & 0 & 1 & 1 & 29 & \textbf{67}  \\
 & Fail & 82 & 100 & 100 & 11 & 15 & 16 & 17 & \textbf{2} & \textbf{2}
    \end{tabular}
    \caption{Computational results \blue{with CPLEX} for univariate transportation problems on large networks.}
    \label{tab:univariate-pwl-large-non-power-of-two}
    \vspace{1em}
            \begin{tabular}{c|rrrrrrr:rr}
Metric & \texttt{MC} & \texttt{CC} & \texttt{SOS2} & \texttt{Inc} & \texttt{DLog} & \Log{} & \texttt{LogIB} & \texttt{ZZB} & \texttt{ZZI} \\ \hline
 Mean (s)  & 1663.4  & 1800.0  & 1800.0  & 710.6  & 752.4  & 793.1  & 796.0  & 319.3  & \textbf{261.4}  \\
 Std  & 298.7  & -  & -  & 529.9  & 555.0  & 570.9  & 554.4  & 392.7  & 316.7  \\
 Win & 0 & 0 & 0 & 4 & 0 & 1 & 0 & 27 & \textbf{53}  \\
 Fail & 78 & 85 & 85 & 10 & 15 & 17 & 18 & 2 & \textbf{1}  \\
 Margin  & - & - & - & 207.0  & - & 5.6  & - & 320.1  & 348.9
    \end{tabular}
    \caption{Difficult univariate transportation problems on large networks.}
    \label{tab:univariate-large-hard-non-power-of-two}
\end{table}



In Table~\ref{tab:univariate-pwl-large-non-power-of-two} we present computational results for the large network instances. Here we observe a roughly 2-3x average speed-up on larger function instances for our new formulations over previous methods. Moreover, we highlight that the new formulations have lower variability in solve time, and time out on fewer instances than the existing methods. With $d=28$, the \SOS{} approach works very well for easier instances, winning on 11 of 100, though its variability is extremely high, timing out on 80 of 100 instances. The existing \Inc{}, \DLog{}, \Log{}, and \LogIB{} formulations all perform roughly comparably.

In Table~\ref{tab:univariate-large-hard-non-power-of-two}, we focus on those large network problems that are difficult (i.e. no approach is able to solve the instance in under 100 seconds) but still solvable (i.e. one formulation solves the instance in under 30 minutes). We see that the new zig-zag formulations are the fastest on 80 of 85 such instances. We also report the average margin: for those instances for which a given new (resp. existing) formulation is fastest, what is the absolute difference in solve time between it and the fastest existing (resp. new) formulation? In this way, we can measure the absolute improvement offered by our new formulation on an instance-by-instance basis. Here we see that the new formulations offer a substantial improvement on these difficult instances, with an absolute decrease of 5-6 minutes in average solve time over existing methods. Finally, we highlight that there are 5 instances that our new formulations can solve to optimality and for which all existing formulations are unable to solve in 30 minutes.

\blue{
We repeat the same experiments with the Gurobi v7.0.2 solver depicted in Tables~\ref{tab:univariate-pwl-non-power-of-two} and \ref{tab:univariate-pwl-large-non-power-of-two}, and include the results in Tables~\ref{tab:univariate-non-power-of-two-gurobi} and \ref{tab:univariate-large-network-non-power-of-two-gurobi}, respectively. On the whole, Gurobi is capable of solving these univariate instances much more efficiently than CPLEX; we omit an analogue of Table~\ref{tab:univariate-large-hard-non-power-of-two} as none of the instances satisfy the specified hardness criteria. Gurobi has a relatively superior implementation of native SOS2 branching that works very well for small and medium function instances. However, it performs very poorly on large function instances (timing out on 98 of 100 instances with $d=59$). We again observe on these larger instances that the \ZZI{} formulation offers a net improvement over the existing host of logarithmic formulations, and is the winner on a plurality of the largest instances in both families. However, as the solve time for all logarithmic formulations on the these largest univariate instances is relatively much lower with Gurobi than CPLEX, the average improvement of the new formulations is more muted than that which can be observed in Tables~\ref{tab:univariate-pwl-large-non-power-of-two} and \ref{tab:univariate-large-hard-non-power-of-two}.

\label{app:univariate-gurobi}
     \begin{table}[htpb]
         \centering
         \smaller
         \begin{tabular}{r|c|rrrrrrr:rr}
$d$ & Metric & \texttt{MC} & \texttt{CC} & \texttt{SOS2} & \texttt{Inc} & \texttt{DLog} & \Log{} & \texttt{LogIB} & \texttt{ZZB} & \texttt{ZZI} \\ \hline
\multirow{4}{*}{6}
 & Mean (s)  & 0.8  & 2.7  & \textbf{0.2}  & 0.5  & 0.7  & 0.7  & 0.7  & 1.0  & 0.7  \\
 & Std  & 0.4  & 3.4  & 0.2 & 0.2 & 0.8  & 0.7  & 0.8  & 0.8  & 0.6  \\
 & Win & 0 & 0 & \textbf{95} & 2 & 1 & 1 & 0 & 0 & 1  \\
 & Fail & \textbf{0} & \textbf{0} & \textbf{0} & \textbf{0} & \textbf{0} & \textbf{0} & \textbf{0} & \textbf{0} & \textbf{0}  \\
\hline
\multirow{4}{*}{13}
 & Mean (s)  & 4.2  & 13.4  & \textbf{0.9}  & 1.9  & 4.1  & 5.2  & 2.1  & 2.5  & 2.7  \\
 & Std  & 4.8  & 15.3  & 1.0  & 0.9  & 4.5  & 6.0  & 2.9  & 2.6  & 2.3  \\
 & Win & 0 & 0 & \textbf{90} & 4 & 0 & 0 & 1 & 2 & 3  \\
 & Fail & \textbf{0} & \textbf{0} & \textbf{0} & \textbf{0} & \textbf{0} & \textbf{0} & \textbf{0} & \textbf{0} & \textbf{0}  \\
\hline
\multirow{4}{*}{28}
 & Mean (s)  & 30.3  & 95.2  & 3.9  & 6.1  & 9.2  & 6.1  & \textbf{3.3}  & 4.4  & 4.4  \\
 & Std  & 43.0  & 261.3  & 8.1  & 5.2  & 8.7  & 10.2  & 2.7 & 4.6  & 3.7  \\
 & Win & 0 & 0 & \textbf{63} & 1 & 1 & 7 & 8 & 7 & 13  \\
 & Fail & \textbf{0} & 2 & \textbf{0} & \textbf{0} & \textbf{0} & \textbf{0} & \textbf{0} & \textbf{0} & \textbf{0}  \\
\hline
\multirow{4}{*}{59}
 & Mean (s)  & 265.5  & 372.3  & 1781.2  & 24.3  & 7.3  & 12.6  & 9.1  & 7.5  & \textbf{6.0}  \\
 & Std  & 409.5  & 530.0  & 134.7  & 23.1  & 6.7  & 12.5  & 9.3  & 7.2  & 5.3 \\
 & Win & 0 & 0 & 0 & 0 & 10 & 20 & 16 & 5 & \textbf{49}  \\
 & Fail & 2 & 8 & 98 & \textbf{0} & \textbf{0} & \textbf{0} & \textbf{0} & \textbf{0} & \textbf{0}
         \end{tabular}
         \caption{Computational results with Gurobi for univariate transportation problems on small networks.}
         \label{tab:univariate-non-power-of-two-gurobi}
         
\vspace{1em}
         
         \begin{tabular}{r|c|rrrrrrr:rr}
$d$ & Metric & \texttt{MC} & \texttt{CC} & \texttt{SOS2} & \texttt{Inc} & \texttt{DLog} & \Log{} & \texttt{LogIB} & \texttt{ZZB} & \texttt{ZZI} \\ \hline
\multirow4{*}{28}
 & Mean (s)  & 124.6  & 245.8  & 1784.8  & 31.5  & 27.1  & 19.8  & \textbf{16.3}  & 19.7  & 17.0  \\
 & Std  & 192.9  & 321.4  & 151.9  & 16.1  & 15.8  & 15.3  & 6.8 & 11.3  & 9.3  \\
 & Win & 0 & 0 & 0 & 0 & 5 & 16 & \textbf{38} & 11 & 30  \\
 & Fail & \textbf{0} & 2 & 99 & \textbf{0} & \textbf{0} & \textbf{0} & \textbf{0} & \textbf{0} & \textbf{0}  \\
\hline
\multirow{4}{*}{59}
 & Mean (s)  & 619.4  & 901.2  & 1800.0  & 87.3  & 23.9  & 27.4  & 26.3  & 24.7  & \textbf{20.9}  \\
 & Std  & 560.3  & 683.5  & - & 53.6  & 19.7  & 11.8 & 14.1  & 16.5  & 16.1  \\
 & Win & 0 & 0 & 0 & 0 & 10 & 9 & 20 & 7 & \textbf{54}  \\
 & Fail & 12 & 27 & 100 & \textbf{0} & \textbf{0} & \textbf{0} & \textbf{0} & \textbf{0} & \textbf{0}
         \end{tabular}
         \caption{Computational results with Gurobi for univariate transportation problems on large networks.}
         \label{tab:univariate-large-network-non-power-of-two-gurobi}
\end{table}
}

\section{Formulations for bivariate piecewise linear functions} \label{sec:bivariate}

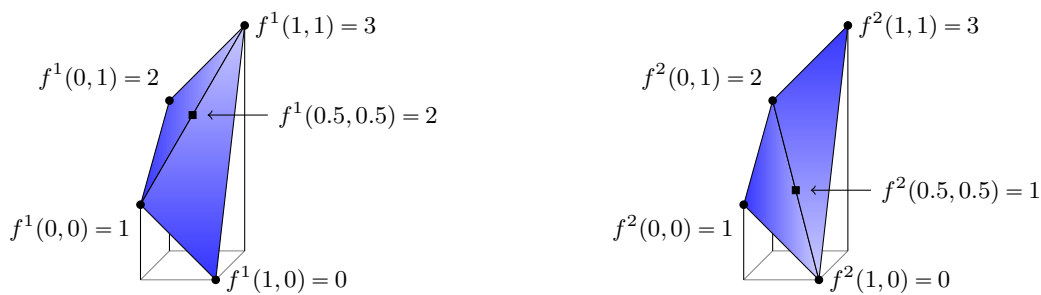
\begin{figure}[htpb]
    \centering
    \begin{tikzpicture}
        \begin{scope}[canvas is xz plane at y=0]
            \draw[help lines] (0,0) grid (1,1);
        \end{scope}

        \begin{scope}[canvas is xy plane at z=1]
            \draw (0,0) -- (0,1);
            \draw (1,0) -- (1,0);
        \end{scope}

        \begin{scope}[canvas is xy plane at z=0]
            \draw (0,0) -- (0,2);
            \draw (1,0) -- (1,3);
        \end{scope}

        \draw [left color=blue!80, right color=blue!20] (0,1,1) -- (0,2,0) -- (1,3,0) -- (0,1,1);
        \draw [bottom color=blue!80,top color=blue!20] (0,1,1) -- (1,0,1) -- (1,3,0) -- (0,1,1);

        \begin{scope}[canvas is xy plane at z=1]
            [ultra thick] \node [below left] at (0,1) {{\smaller $f^1(0,0)=1$}};
            [ultra thick] \node [right] at (1,0) {{\smaller $f^1(1,0)=0$}};
            \draw [fill] (0,1) circle [radius=.05];
            \draw [fill] (1,0) circle [radius=.05];
        \end{scope}

        \begin{scope}[canvas is xy plane at z=0]
            [ultra thick] \node [above left] at (0,2) {{\smaller $f^1(0,1)=2$}};
            [ultra thick] \node [right] at (1,3) {{\smaller $f^1(1,1)=3$}};
            \draw [fill] (0,2) circle [radius=.05];
            \draw [fill] (1,3) circle [radius=.05];
        \end{scope}

        [ultra thick] \node [right] at (1.5,2,0.5) {{\smaller $f^1(0.5,0.5)=2$}};
        \draw [fill] (0.5-0.075,2-0.075,0.5-0.075) rectangle (0.5+0.075,2+0.075,0.5+0.075);
        \draw [->] (1.5,2,0.5) -- (0.7,2,0.5);
    \end{tikzpicture}\hspace{5em}
    \begin{tikzpicture}
        \begin{scope}[canvas is xz plane at y=0]
            \draw[help lines] (0,0) grid (1,1);
        \end{scope}

        \begin{scope}[canvas is xy plane at z=1]
            \draw (0,0) -- (0,1);
            \draw (1,0) -- (1,0);
        \end{scope}

        \begin{scope}[canvas is xy plane at z=0]
            \draw (0,0) -- (0,2);
            \draw (1,0) -- (1,3);
        \end{scope}

        \draw [left color=blue!80, right color=blue!20] (0,1,1) -- (0,2,0) -- (1,0,1) -- (0,1,1);
        \draw [top color=blue!80, bottom color=blue!20] (1,0,1) -- (0,2,0) -- (1,3,0) -- (1,0,1);

        \begin{scope}[canvas is xy plane at z=1]
            [ultra thick] \node [below left] at (0,1) {{\smaller $f^2(0,0)=1$}};
            [ultra thick] \node [right] at (1,0) {{\smaller $f^2(1,0)=0$}};
            \draw [fill] (0,1) circle [radius=.05];
            \draw [fill] (1,0) circle [radius=.05];
        \end{scope}

        \begin{scope}[canvas is xy plane at z=0]
            [ultra thick] \node [above left] at (0,2) {{\smaller $f^2(0,1)=2$}};
            [ultra thick] \node [right] at (1,3) {{\smaller $f^2(1,1)=3$}};
            \draw [fill] (0,2) circle [radius=.05];
            \draw [fill] (1,3) circle [radius=.05];
        \end{scope}

        [ultra thick] \node [right] at (1.5,1,0.5) {{\smaller $f^2(0.5,0.5)=1$}};
        \draw [fill] (0.5-0.075,1-0.075,0.5-0.075) rectangle (0.5+0.075,1+0.075,0.5+0.075);
        \draw [->] (1.5,1,0.5) -- (0.7,1,0.5);
    \end{tikzpicture}
    \caption{Two bivariate functions over $D=[0,1]^2$ that match on the gridpoints, but differ on the interior of $D$.}
    \label{fig:triangulation-difference}
\end{figure}

Bivariate piecewise linear functions possess a much more complex structure than their univariate counterparts, which means that constructing logarithmic formulations for them is also correspondingly more difficult. This combinatorial structure is endowed by the pattern \blue{into} which the domain is \blue{decomposed}, the choice of which determines the values which the bivariate piecewise function takes (see Figure~\ref{fig:triangulation-difference} for an illustration). Although it is possible to extend the geometric construction of Proposition~\ref{prop:sos2} to the bivariate setting~\citep{Huchette:2019}, this technique requires us to compute the hyperplanes spanned by high-dimensional vectors  a la Proposition~\ref{prop:sos2}, which is, in general, very difficult. Instead, we turn to a combinatorial approach.

\blue{
For the remainder of the section, we will focus on bivariate functions with \emph{grid triangulation} domains, which we define formally as follows.
\begin{definition}
    Presume that $V = \llbracket d_1 + 1 \rrbracket \times \llbracket d_1 + 1 \rrbracket$. A \emph{grid triangulation} $\calT$ of $V$ is a family of sets $\calT$ where:
    \begin{itemize}
        \item Each $T \in \calT$ is a triangle: $|T|=3$.
        \item $\calT$ partitions the domain: $\bigcup_{T \in \calT} \Conv(T) = \Conv(V)$ and $\relint(\Conv(T)) \cap \relint(\Conv(T')) = \emptyset$ for each distinct $T, T' \in \calT$ (where $\relint(S)$ is the (relative) interior of set $S$).
        \item $\calT$ is on a regular grid: $||v-w||_\infty \leq 1$ for each $T \in \calT$ and $v,w \in T$.
    \end{itemize}
\end{definition}
}

\subsection{Independent branching formulations} \label{ssec:independent-branching}

The original logarithmic formulation \LogIB{} of~\cite{Vielma:2009a} for the SOS2 constraint is derived from the class of \emph{independent branching} formulations, which offers a combinatorial way of constructing formulations. \cite{Huchette:2016a} offer a complete characterization of its expressive power, as well as a graphical procedure to systematically construct independent branching formulations. \blue{
We start with two definitions.

\begin{definition}
Take a combinatorial disjunctive constraint given by the family of subsets $(T^i)_{i=1}^d$ over a ground set $V$.
\begin{itemize}
    \item The \emph{conflict graph} $G=(V,E)$ of the combinatorial disjunctive constraint is given by the edge set $E = \Set{\{u,v\} \in [V]^2 | \{u,v\} \not\subseteq T^i \text{ for all } i \in \llbracket d \rrbracket}$, where $[V]^2 \defeq \Set{\{u,v\} \in V \times V | u \neq v}$.
    \item A \emph{biclique} of some graph $G=(V,E)$ is a pair of sets $(A,B)$ such that $(V,A * B)$ is a subgraph of $G$ (i.e. $A * B \subseteq E)$, where $A * B \defeq \Set{ \{u,v\} \in [V]^2 | u \in A, v \in B}$.
    \item A \emph{biclique cover} of some graph $G = (V,E)$ is a family of bicliques $\{(A^k,B^k)\}_{k=1}^r$ such that $E = \bigcup_{k=1}^r (A^k * B^k)$. We say that such a biclique cover has $r$ \emph{levels}.
\end{itemize}
\end{definition}

Given a biclique cover for the conflict graph of a (suitably representable) combinatorial disjunctive constraint, the technique of \cite{Huchette:2016a} directly constructs a formulation as follows.}

\begin{proposition}[\cite{Huchette:2016a}]\label{bicliqueprop}
Let $\calT = (T^i \subseteq V)_{i=1}^d$ be the family of sets corresponding to either a univariate piecewise linear function, or a bivariate piecewise linear function with a grid triangulated domain. Take $E$ as the edge set for the conflict graph corresponding to $\calT$. \blue{If $\{(A^k,B^k)\}_{k=1}^r$ is a biclique cover for $(V,E)$,}
then an ideal independent branching formulation
  for $\bigcup_{i=1}^d P(T^i)$ is
  \begin{equation} \label{bicliqueform}
      \sum\nolimits_{v \in A^k} \lambda \leq y_k, \quad\quad \sum\nolimits_{v \in B^k} \lambda_v \leq 1 - y_k \quad \forall k \in \llbracket r \rrbracket,\quad
      (\lambda,y) \in \Delta^V \times \{0,1\}^r.
  \end{equation}
\end{proposition}
Intuitively, this formulation ensures that, for each level $k$, either $\lambda_v = 0$ for all $v \in A^k$, or $\lambda_v = 0$ for all $v \in B^k$.

As motivation, we return to Example~\ref{example-1} to construct the logarithmic independent branching formulation for the SOS2 constraint, \LogIB{}, as introduced by \cite{Vielma:2009a}.

\begin{example}\label{example:logIB}
Take the SOS2 constraint with $d=4$ (as seen in \eqref{eqn:pwl-example}). The edge set for the conflict graph is $E = \left\{\{1,3\},\{1,4\},\{1,5\},\{2,4\},\{2,5\},\{3,5\}\right\}$, which admits a biclique \blue{cover} with the sets $A^1 = \{3\}$, $B^1 = \{1,5\}$, $A^2 = \{4,5\}$, and $B^2 = \{1,2\}$. The corresponding \LogIB{} formulation is then
\begin{equation}\label{eqn:logIB-example}
    \lambda_3 \leq y_1, \quad\quad \lambda_1 + \lambda_5 \leq 1 - y_1, \quad\quad \lambda_4 + \lambda_5 \leq y_2, \quad\quad \lambda_1 + \lambda_2 \leq 1-y_2, \quad\quad (\lambda,y) \in \Delta^{V} \times \{0,1\}^2.
\end{equation}
See Figure~\ref{fig:biclique-example} for an illustration. As noted previously, \blue{the} \LogIB{} formulation \eqref{eqn:logIB-example} coincides with the \Log{} formulation \eqref{eqn:log-example} because $d$ is a power-of-two; see Appendix~\ref{app:log-logib} for an instance where this is not the case.

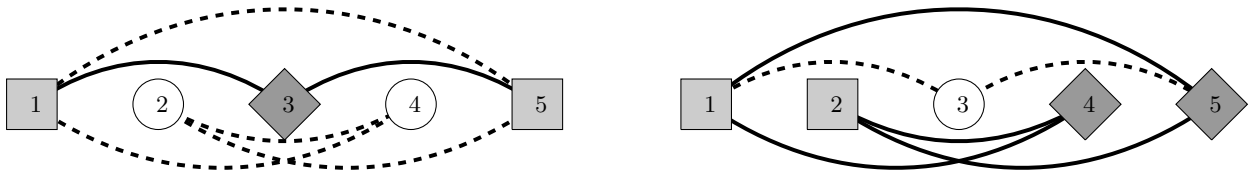
\begin{figure}
\centering
\begin{tikzpicture}[scale=0.84]
        \draw [dashed,line width=1.5] (0,0) to [out=40,in=140] (8,0);
        \draw [       line width=1.5] (0,0) to [out=35,in=145] (4,0);
        \draw [dashed,line width=1.5] (0,0) to [out=-35,in=-145] (6,0);
        \draw [       line width=1.5] (4,0) to [out=35,in=145] (8,0);
        \draw [dashed,line width=1.5] (2,0) to [out=-35,in=-145] (8,0);
        \draw [dashed,line width=1.5] (2,0) to [out=-30,in=-150] (6,0);

        \draw [fill=gray!40] (-0.4,-0.4) -- (0.4,-0.4) -- (0.4,0.4) -- (-0.4,0.4) -- (-0.4,-0.4);
        \draw [fill=white] (2,0) circle [radius=0.4];
        \draw [fill=gray!80] (4,0-0.5656854249492381) -- (4+0.5656854249492381,0) -- (4,0+0.5656854249492381) -- (4-0.5656854249492381,0) -- (4,0-0.5656854249492381);
        \draw [fill=white] (6,0) circle [radius=0.4];
        \draw [fill=gray!40] (8-0.4,-0.4) -- (8+0.4,-0.4) -- (8+0.4,0.4) -- (8-0.4,0.4) -- (8-0.4,-0.4);

        \node at (0,0) {\footnotesize $1$};
        \node at (2,0) {\footnotesize $2$};
        \node at (4,0) {\footnotesize $3$};
        \node at (6,0) {\footnotesize $4$};
        \node at (8,0) {\footnotesize $5$};
\end{tikzpicture}
\hfill
\begin{tikzpicture}[scale=0.84]
        \draw [       line width=1.5] (0,0) to [out=40,in=140] (8,0);
        \draw [dashed,line width=1.5] (0,0) to [out=35,in=145] (4,0);
        \draw [       line width=1.5] (0,0) to [out=-35,in=-145] (6,0);
        \draw [dashed,line width=1.5] (4,0) to [out=35,in=145] (8,0);
        \draw [       line width=1.5] (2,0) to [out=-35,in=-145] (8,0);
        \draw [       line width=1.5] (2,0) to [out=-30,in=-150] (6,0);

        \draw [fill=gray!40] (-0.4,-0.4) -- (0.4,-0.4) -- (0.4,0.4) -- (-0.4,0.4) -- (-0.4,-0.4);
        \draw [fill=gray!40] (2-0.4,-0.4) -- (2+0.4,-0.4) -- (2+0.4,0.4) -- (2-0.4,0.4) -- (2-0.4,-0.4);
        \draw [fill=white] (4,0) circle [radius=0.4];
		\draw [fill=gray!80] (6+0,0-0.5656854249492381) -- (6+0.5656854249492381,0) -- (6,0+0.5656854249492381) -- (6-0.5656854249492381,0) -- (6,0-0.5656854249492381);
		\draw [fill=gray!80] (8+0,0-0.5656854249492381) -- (8+0.5656854249492381,0) -- (8,0+0.5656854249492381) -- (8-0.5656854249492381,0) -- (8,0-0.5656854249492381);
        \node at (0,0) {\footnotesize $1$};
        \node at (2,0) {\footnotesize $2$};
        \node at (4,0) {\footnotesize $3$};
        \node at (6,0) {\footnotesize $4$};
        \node at (8,0) {\footnotesize $5$};
\end{tikzpicture}
\caption{The biclique cover for the conflict graph $G$ of the SOS2 constraint in Example~\ref{example:logIB}. \textbf{(Left)} The first level with $A^1$ and $B^1$ are diamonds and squares, respectively; and \textbf{(Right)} similarly for $A^2$ and $B^2$ in the second level. For each level, covered edges are solid and omitted edges are dashed.}
\label{fig:biclique-example}
\end{figure}

\end{example}



\subsection{Independent branching formulations for bivariate piecewise linear functions}\label{PWLTDSEC}

Recall that, using Proposition~\ref{bicliqueprop}, we can immediately construct a formulation for a bivariate function that is ideal and of size $\mathcal{O}(r)$ if we can find a biclique cover with $r$ levels for the corresponding conflict graph. \blue{A natural question then is: what is the smallest biclique cover that be constructed for a given grid triangulation? 

Computing a minimum cardinality biclique cover is NP-hard, even on bipartite graphs~\citep{Garey:1979}. \footnote{To the best of our knowledge, efficient algorithms for certain classes of structured non-bipartite graphs (e.g. the conflict graph of a grid triangulation), have not been investigated in the literature.}} \cite{Vielma:2009a} consider a highly structured grid triangulation known as the J1 or Union Jack~\citep{Todd:1977}, and (implicitly) present a biclique cover with $r = \lceil \log_2(d_1) \rceil + \lceil \log_2(d_2) \rceil + 1$ levels. More recently, \cite{Huchette:2016a} propose a \blue{construction} under a weaker structural condition involving the existence of a certain graph coloring that uses at most $r=\lceil \log_2(d_1) \rceil + \lceil \log_2(d_2) \rceil + 2$ levels, as well as a \blue{construction} for arbitrary grid triangulations with $r=\lceil \log_2(d_1) \rceil + \lceil \log_2(d_2) \rceil + 9$ levels. In this work, we present a new, even smaller \blue{construction} that is applicable for any grid triangulation.

\blue{

For the remainder of the subsection, consider a family of sets $\calT$ associated with a grid triangulation with $V = 
\llbracket d_1 + 1 \rrbracket \times \llbracket d_2 + 1 \rrbracket$, along with the edge set $E$ for its corresponding conflict graph. It will be useful to decompose these edges into three classes:
\begin{align*}
    E^\nearrow &\defeq \Set{ \{u,v\} \in E | \: ||u-v||_\infty = 1, \: |(u_1-v_1) + (u_2-v_2)| = 2 } \\
    E^\searrow &\defeq \Set{ \{u,v\} \in E | \: ||u-v||_\infty = 1, \: |(u_1-v_1) + (u_2-v_2)| = 0 } \\
    E^F &\defeq \Set{\{u,v\} \in E |\: ||u-v||_\infty > 1}.
\end{align*}
In words, $E^\nearrow$ and $E^\searrow$ are those ``nearby'' edges sharing a common subrectangle that are oriented along a diagonal (southwest to northeast) or an anti-diagonal (southeast to northwest), respectively. Contrastingly, $E^F$ are those edges that are ``far apart.'' It is straightforward from their definition to see that the three classes form a partition of $E$.

First, we show how to exactly cover the ``far apart'' edges by applying an aggregated SOS2 construction along each axis.

\begin{lemma} \label{lemma:aggregated-sos2}
    Take $\{(A^{1,k},B^{1,k})\}_{k=1}^{r_1}$ and $\{(A^{2,k},B^{2,k})\}_{k=1}^{r_2}$ as biclique covers for the conflict graphs associated with the SOS2 constraints on $d_1+1$ and $d_2+1$ breakpoints, respectively. Then
    \begin{equation} \label{eqn:far-away-edges}
        E^F = \bigcup_{k=1}^{r_1} \left((A^{1,k} \times \llbracket d_2 + 1 \rrbracket) * (B^{1,k} \times \llbracket d_2 + 1 \rrbracket)\right) \cup \bigcup_{k=1}^{r_2} \left((\llbracket d_1 + 1 \rrbracket \times A^{2,k}) * (\llbracket d_1 + 1 \rrbracket \times B^{2,k})\right).
    \end{equation}
\end{lemma}

\proof{Proof}
The conflict graph for the SOS2 constraint on $d+1$ breakpoints is given by the edge set $\Set{ \{u,v\} \in [d+1]^2 |\: |u-v| > 1}$. Therefore, if the families of sets are taken as described, we can infer that
\begin{subequations} \label{eqn:agg}
\begin{align}
    \Set{\{u,v\} \in [V]^2 |\: |u_1-v_1| > 1} &= \bigcup_{k=1}^{r_1} \left((A^{1,k} \times \llbracket d_2 + 1 \rrbracket) * (B^{1,k} \times \llbracket d_2 + 1 \rrbracket)\right) \\
    \Set{\{u,v\} \in [V]^2 |\: |u_2-v_2| > 1} &= \bigcup_{k=1}^{r_2} \left((\llbracket d_1 + 1 \rrbracket \times A^{2,k}) * (\llbracket d_1 + 1 \rrbracket \times B^{2,k})\right).
\end{align}
\end{subequations}
Observe that the union of the left-hand sides of both equations in \eqref{eqn:agg} is equal to $E^F$, while the union of the right-hand sides of both equations in \eqref{eqn:agg} is identical to the right-hand side of~\eqref{eqn:far-away-edges}, giving the result.
\Halmos\endproof

Next, we introduce an algorithm that constructs three bicliques. Each ``nearby'' edge in the conflict graph that is oriented along a diagonal will be contained in one of these 3 bicliques. However, these 3 bicliques also introduce additional edges that are not ``nearby'' edges oriented along a diagonal, and so the bicliques do not induce subgraphs of the graph induced by the diagonal edges, $(V,E^\nearrow)$. Fortunately, we will show that these extra edges are ``far apart'' edges contained in $E^F$, and so the 3 bicliques \emph{do} induce subgraphs of the conflict graph $(V,E)$ as required by Proposition~\ref{bicliqueprop}. Therefore, we can use these bicliques to construct a cover for the conflict graph; they will cover each ``nearby'' edge oriented along a diagonal, while not introducing any undesired additional ``nearby'' edges.


\begin{algorithm}
    \caption{Computing bicliques that cover diagonal triangle selection edges.}
    \begin{algorithmic}[1]
    \Require{Integer $\tau \in \{0,1,2\}$.}
    \Procedure{DiagonalBicliques}{$\tau$,E}
        \For{$\kappa \gets \tau, \tau+3, \ldots, d_1-1$}
            \Let{$\phi$}{\texttt{true}}
            \For{$i \gets (1+\kappa), \ldots, d_1$}
                \Let{$j$}{$i - \kappa$}
                \If{$1 \leq j \leq d_2$ and $\{(i,j),(i+1,j+1)\} \in E$}
                    \If{$\phi$}
                        \State Insert $(i,j) \to A$ and $(i+1,j+1) \to B$ 
                    \Else
                        \State Insert $(i+1,j+1) \to A$ and $(i,j) \to B$
                    \EndIf
                    \Let{$\phi$}{$\neg\phi$} \Comment{The unary operator $\neg$ denotes the negation of a boolean}
                \EndIf
            \EndFor
        \EndFor
        \For{$\kappa \gets 3-\tau, 6-\tau, \ldots, d_2-1$}
            \Let{$\phi$}{\texttt{true}}
            \For{$j \gets (1+\kappa), \ldots, d_2$}
                \Let{$i$}{$j - \kappa$}
                \If{$1 \leq i \leq d_1$ and $\{(i,j),(i+1,j+1)\} \in E$}
                    \If{$\phi$}
                        \State Insert $(i,j) \to A$ and $(i+1,j+1) \to B$ 
                    \Else
                        \State Insert $(i+1,j+1) \to A$ and $(i,j) \to B$ 
                    \EndIf
                    \Let{$\phi$}{$\neg\phi$}
                \EndIf
            \EndFor
        \EndFor
        \State \Return{$(A,B)$}
    \EndProcedure
    \end{algorithmic}
    \label{alg:diagonal}
\end{algorithm}

\begin{lemma} \label{lemma:diagonal}
    For each $\tau \in \{0,1,2\}$, construct the sets $(A^\tau,B^\tau) = \textsc{DiagonalBicliques}(\tau,E^\nearrow)$ according to the procedure listed in Algorithm~\ref{alg:diagonal}. Then $E^\nearrow \subseteq \bigcup_{\tau=0}^2 (A^\tau * B^\tau) \subseteq E$.
\end{lemma}
\proof{Proof}
\underline{The first inclusion: $E^\nearrow \subseteq \bigcup_{\tau=0}^2 (A^\tau * B^\tau)$} Consider some arbitrary $(i,j) \in \llbracket d_1 \rrbracket \times \llbracket d_2 \rrbracket$. It suffices to show that there exists some value for $\tau \in \{0,1,2\}$ and some $\kappa$ such that either: a) $\kappa \in \{\tau,\tau+3,\ldots,d_1-1\}$ and $j = i - \kappa$, or b) $\kappa \in \{3-\tau,6-\tau,\ldots,d_2-1\}$ and $i = j - \kappa$. If this is the case, then Algorithm~\ref{alg:diagonal} with reach either line 6 or 16, respectively, with the appropriate values for $(i,j)$, and so if $\{(i,j),(i+1,j+1)\} \in E^\nearrow$, by construction we will have build $(A,B)$ such that $\{(i,j),(i+1,j+1)\} \in A * B$.

To show that such a values exists, first consider the case where $i > j$, in which case $i - j \in \{0,\ldots,d_1-1\}$. It is straightforward to see that, if we take $\tau = (i-j)\mod 3$, then $\tau \in \{0,1,2\}$, and moreover one of the iterations of the for loop initiated in line 2 will have $\kappa = i-j$, giving the desired result. Similarly, if $i < j$, we can attain any value $\kappa \in \{1,\ldots,d_2-1\}$ in the loop initiated in line 12; choose $\kappa = j - i$ to give the result. Therefore, we conclude that $E^\nearrow \subseteq \bigcup_{\tau=0}^2 (A^\tau * B^\tau)$.

\underline{To second inclusion: $\bigcup_{\tau=0}^2 (A^\tau * B^\tau) \subseteq E$} We start by observing that, due to the for loop ranges and the explicit checks on the values of $j$ and $i$ in lines 6 and 16, respectively, any pair $\{u=(i,j),v=(i+1,j+1)\}$ that could possibly be inserted into $(A,B)$ at lines 8, 10, 18, or 20 will naturally satisfy $u,v \in V \equiv \llbracket d_1 + 1 \rrbracket \times \llbracket d_2 + 1 \rrbracket$.

Next, observe that by the definition of a grid triangulation, for each $u,v \in V$ with $||u-v||_\infty > 1$, necessarily $\{u,v\} \in E^F \subseteq E$. Therefore, the result follows if we can show that, for each $\tau \in \{0,1,2\}$, any $\{u,v\} \in (A^\tau * B^\tau) \backslash E^\nearrow$ satisfies $||u-v||_\infty > 1$, and therefore $\{u,v\} \in E^F$.

For the remainder of the proof fix $\tau \in \{0,1,2\}$ and presume some element $u=(i+\Delta,j+\Delta)$ was inserted into set $A$ on line 8 ($\Delta=0$) or on line 10 ($\Delta=1)$ with the iteration value $\kappa$ for the loop initiated on line 2. We will inspect possible values $v = (i'+\Delta',j'+\Delta)$ that can be inserted into set $B$ at lines 8 or 18 ($\Delta'=1$), or lines 10 or 20 ($\Delta'=0$) to verify that $\{u,v\} \in E^F$. All other possible cases will follow by symmetry.

First, consider possible insertions of $v$ to set $B$ on line 8 or 10 with the same iteration value $\kappa$ in the for loop initiated on line 2. Presume that $i' \neq i$ are distinct iteration values for the loop initiated on line 4. Take $\phi$ and $\phi'$ as the values for the boolean for each loop iteration, and presume w.l.o.g. that $\phi=\texttt{true}$ (and, therefore, that $\Delta=0$). If $i'=i+1$, then due to the negation on line 11, $\phi' = \texttt{false}$. Therefore, the two passes through the loop introduce exactly the elements $(i,j), (i+2,j+2) \to A$ and $(i+1,j+1) \to B$, and so no unnecessary edges are introduced. An analogous argument holds if $i'=i-1$. If, on the other hand, $|i-i'| > 1$, presume w.l.o.g. that $i' > i$, and observe that the elements added to the sets in the two loop iterations are $(i,j), (i'+\Delta',j'+\Delta') \to A$ and $(i+1,j+1), (i'+1-\Delta',j'+1-\Delta') \to B$, where $\Delta'=0$ if $\phi' = \texttt{true}$ and $\Delta'=1$ otherwise. Then $|i-(i'+1-\Delta')| = |(i'-i)+(1-\Delta')| > 1$, and so $\{(i,j),(i'+\Delta')\} \in E^F$. A similar argument holds for the other edges introduced, and for the case where $\phi=\texttt{false}$, and so, restricted to this single pass through the for loop, we have the result.

Next, consider possible insertions of $v$ to set $B$ with some distinct iteration value $\kappa' \neq \kappa$ in the for loop initiated on line 2. Since $\tau$ is fixed, we have that $|\kappa-\kappa'| \geq 3$. If $|(i+\Delta)-(i'+\Delta')| > 1$, we are done, so presume otherwise. In this case, since $|(i+\Delta)-(i'+\Delta')| \leq 1$ and $|\kappa-\kappa'| \geq 3$, we conclude that 
\begin{align*}
    |(j+\Delta)-(j'+\Delta')| &= |(i+\Delta-\kappa) - (i'+\Delta'-\kappa')| \\
    &= |((i+\Delta)-(i'+\Delta')) - (\kappa-\kappa')| \\
    &\geq \left||(i+\Delta)-(i'+\Delta')| - |\kappa-\kappa'|\right| \\
    &\geq |1 - 3| = 2.
\end{align*}
Therefore, $\{u,v\} \in E^F$.

Finally, consider the case where $v$ was inserted into set $B$ in the for loop initiated on line 12. Define the quantity $\gamma = (i+\Delta) - (j' + \Delta')$. Using the identities $j=i-\kappa$ and $i'=j'-\kappa'$ given by lines 5 and 15, respectively, we can write
\begin{subequations}
\begin{align}
    |(i+\Delta)-(i'+\Delta')| &= |(i+\Delta) - (j'+\Delta'-\kappa')| = |\gamma + \kappa'| \label{eqn:first-component} \\
    |(j+\Delta)-(j'+\Delta')| &= |(i+\Delta - \kappa) - (j' + \Delta')| = |\gamma - \kappa| \label{eqn:second-component}
\end{align}
\end{subequations}
Therefore, in order for $|(i+\Delta)-(i'+\Delta')| \leq 1$ (\emph{condition 1}), we must have $\gamma + \kappa' \in \{-1,0,1\}$, i.e. $\gamma \in \{-\kappa'-1,-\kappa',-\kappa'+1\}$. Since $\kappa' \geq 3-\tau$ from the loop iteration definition of line 12, we can then infer that $\gamma \leq \tau - 2$ if this condition holds. Similarly, in order for $|(j+\Delta)-(j'+\Delta')| \leq 1$ (\emph{condition 2}), we must have $\gamma - \kappa \in \{-1,0,1\}$, i.e. $\gamma \in \{\kappa-1,\kappa,\kappa+1\}$. Since $\kappa \geq \tau$ from the loop iteration definition of line 2, we can infer that $\gamma \geq \tau-1$ for this condition to hold. As there does not exist a value of $\gamma$ such that $\gamma \leq \tau - 2$ and $\gamma \geq \tau - 1$, we can infer that both condition 1 and condition 2 cannot hold at the same time. Equivalently, either $|u_1-v_1| > 1$ or $|u_2-v_2| > 1$, which implies that $||u-v||_\infty > 1$. Therefore, $\{u,v\} \in E^F$, completing the proof.
\Halmos\endproof

It is straightforward to adapt the construction from Lemma~\ref{lemma:diagonal} to separate the nearby edges $E^\searrow$ along the anti-diagonals, by, for example, reflecting the ground set $V$ along $x_1$ direction via the invertible mapping $M: V \to V$ where $M(u,v) = (d_1+2-u,v)$, and then applying Algorithm~\ref{alg:diagonal} to the edge set transformed using this mapping. Additionally, in Appendix~\ref{app:anti-diagonal} we give an explicit statement from first principles.

\begin{lemma} \label{lemma:antidiagonal}
    There exists a family of bicliques $\{(A^k,B^k)\}_{k=0}^2$ such that $E^\searrow \subseteq \bigcup_{\tau=0}^2 (A^\tau * B^\tau) \subseteq E$.
\end{lemma}

Combining these results gives an explicit procedure to construct a small formulation for arbitrary bivariate grid triangulations.

\begin{theorem}
    There exists an independent branching formulation for a bivariate grid triangulation over $V = \llbracket d_1 + 1 \rrbracket \times \llbracket d_2 + 1 \rrbracket$ of depth $\lceil \log_2(d_1) \rceil + \lceil \log_2(d_2) \rceil + 6$. 
\end{theorem}
\proof{Proof}
The construction follows by applying Proposition~\ref{bicliqueprop} to the biclique cover given by the union of all pairs of sets as defined in Lemmas~\ref{lemma:aggregated-sos2}, \ref{lemma:diagonal}, and \ref{lemma:antidiagonal}. The constructions from Lemmas~\ref{lemma:diagonal} and \ref{lemma:antidiagonal} each introduce 3 bicliques. Furthermore, as noted in Section~\ref{sec:univariate}, \cite{Vielma:2009a} presented an independent branching formulation for SOS2 on $d+1$ breakpoints that requires $\lceil \log_2(d) \rceil$ level, meaning that we can adopt the construction of Lemma~\ref{lemma:aggregated-sos2} using $\lceil \log_2(d_1) \rceil + \lceil \log_2(d_2) \rceil$ bicliques, giving the result.
\Halmos\endproof

We can show the construction pictorially in Figures~\ref{fig:aggregated-sos2} and \ref{fig:6-stencil} on one particular grid triangulation with $d_1 = d_2 = 8$. In Figure~\ref{fig:aggregated-sos2}, we see the bicliques derived in Lemma~\ref{lemma:aggregated-sos2}, where a logarithmically-sized SOS2 biclique is ``aggregated'' vertically or horizontally (top and bottom rows, respectively). In the top row of Figure~\ref{fig:6-stencil}, we see the construction derived in Lemma~\ref{lemma:diagonal}, which covers all diagonal edges $E^\nearrow$. Note the three panels, each of which aggregates diagonal lines that are sufficiently far apart, starting with an offset of $\tau \in \{0,1,2\}$ (from left to right). In particular, note that all of the edges introduced that are not in $E^\nearrow$ are sufficiently far apart that they are contained in $E^F$, meaning that no undesired edges are introduced. In the second row of Figure~\ref{fig:6-stencil}, we see the analogous construction that covers the antidiagonal lines presented in Lemma~\ref{lemma:antidiagonal}.
}

\begin{figure}[htpb]
    \centering
	\includegraphics[width=.25\linewidth]{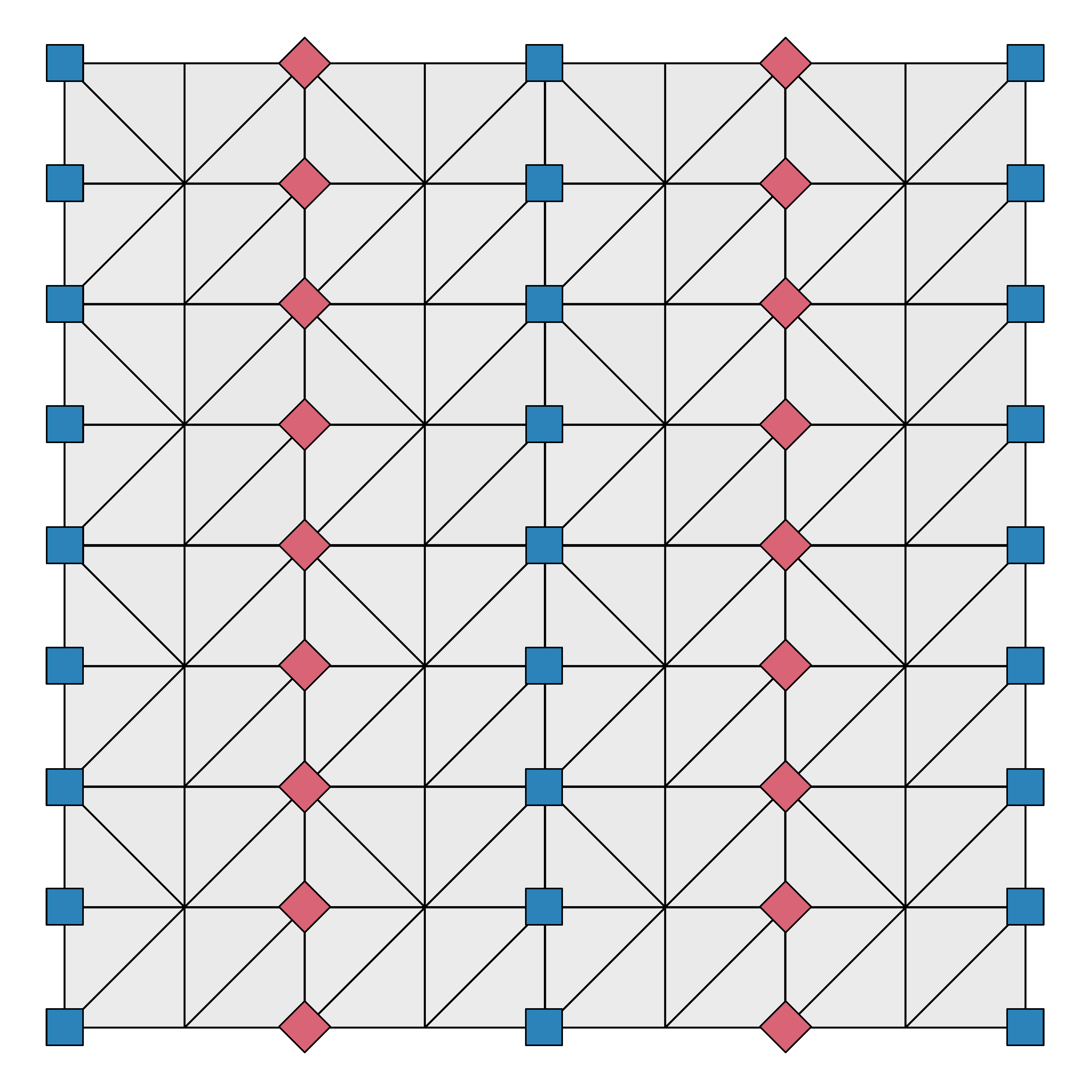}
    \includegraphics[width=.25\linewidth]{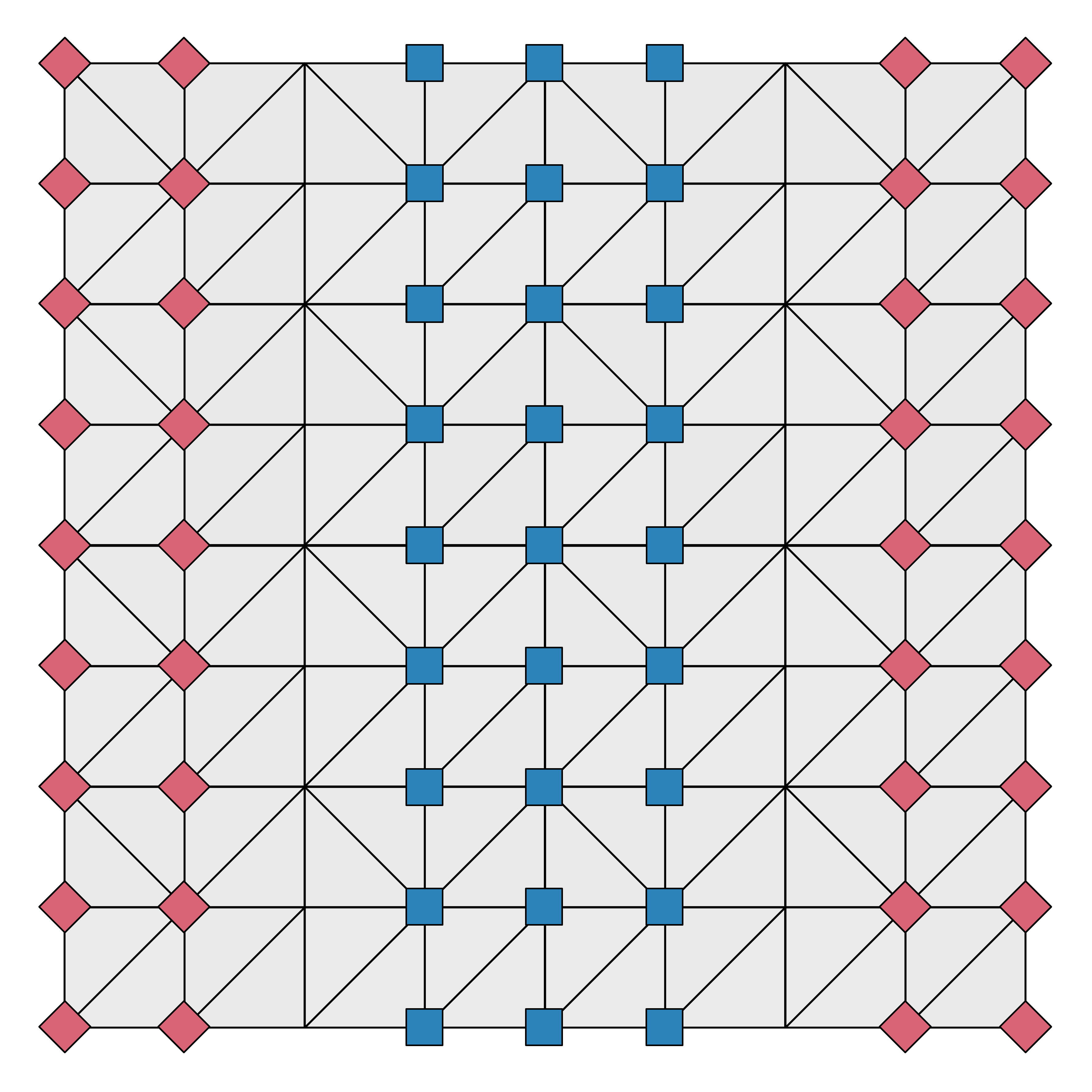}
    \includegraphics[width=.25\linewidth]{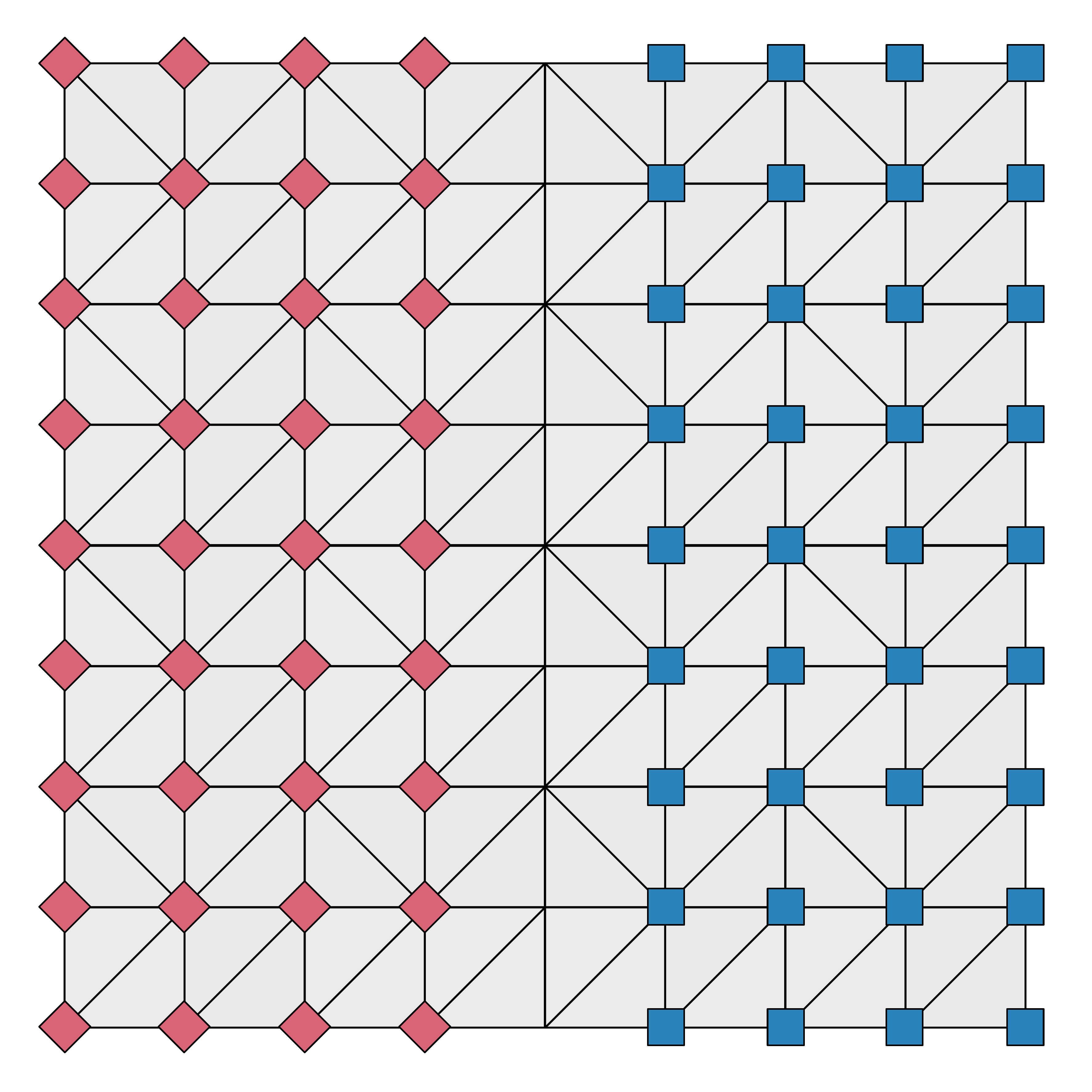} \\
    \includegraphics[width=.25\linewidth]{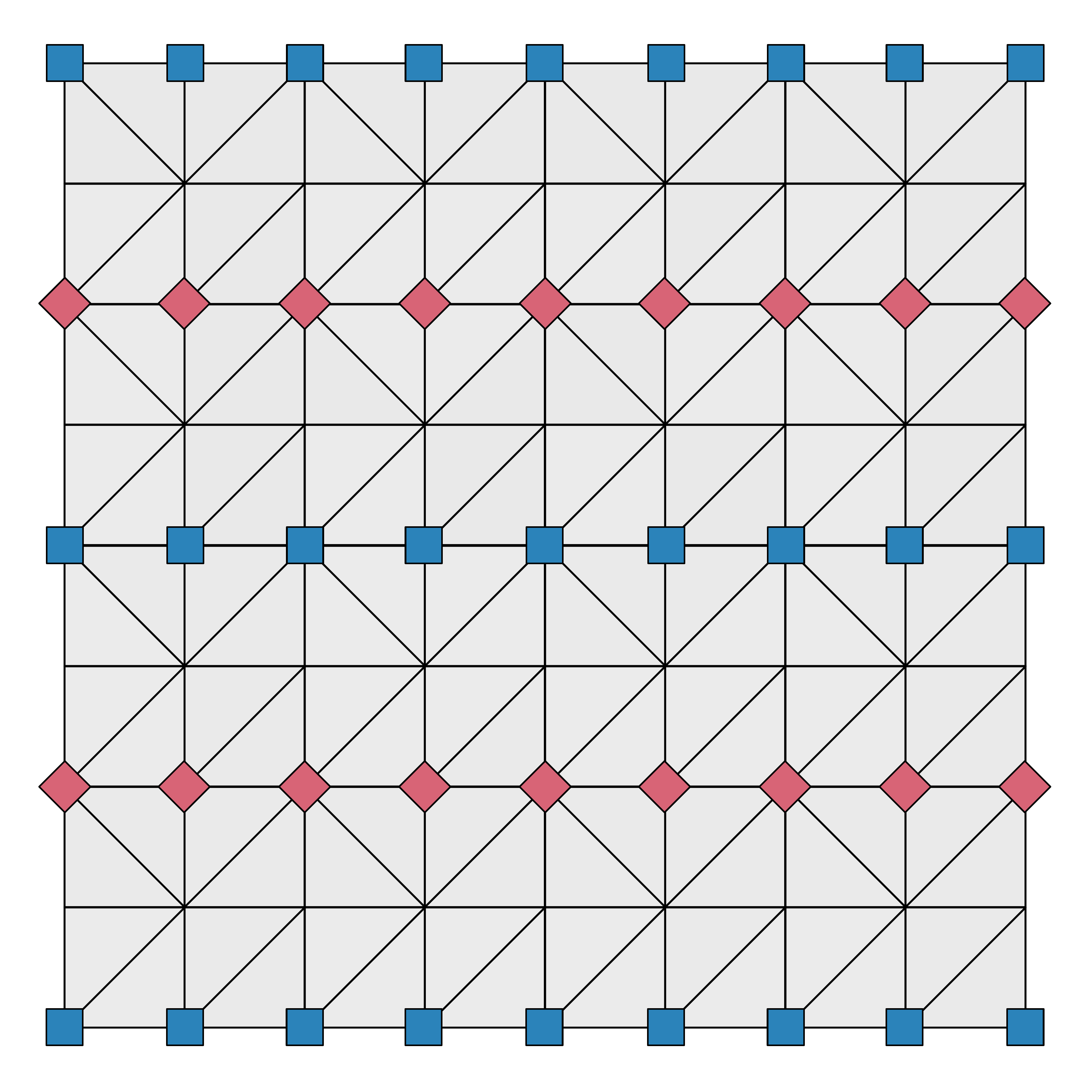}
    \includegraphics[width=.25\linewidth]{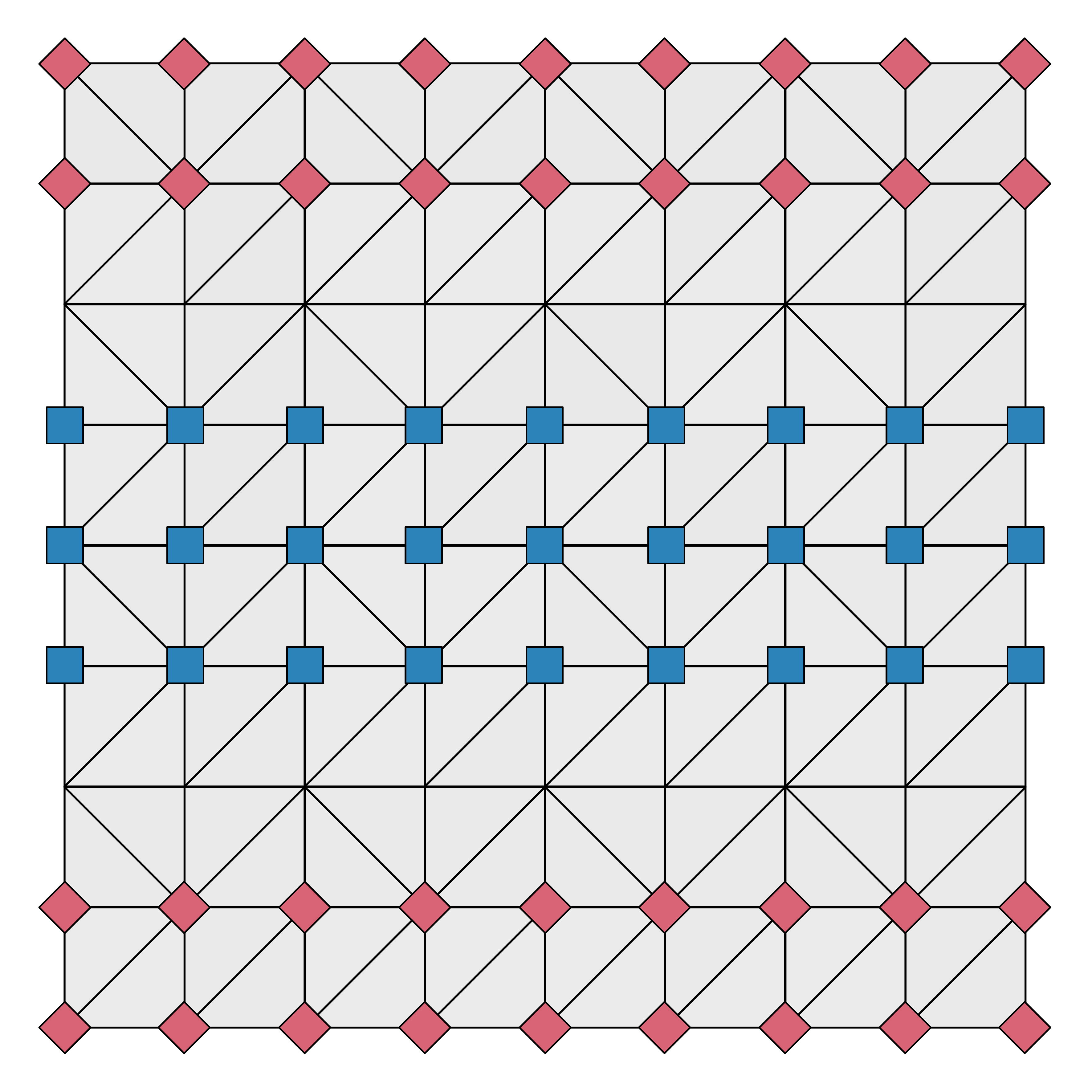}
    \includegraphics[width=.25\linewidth]{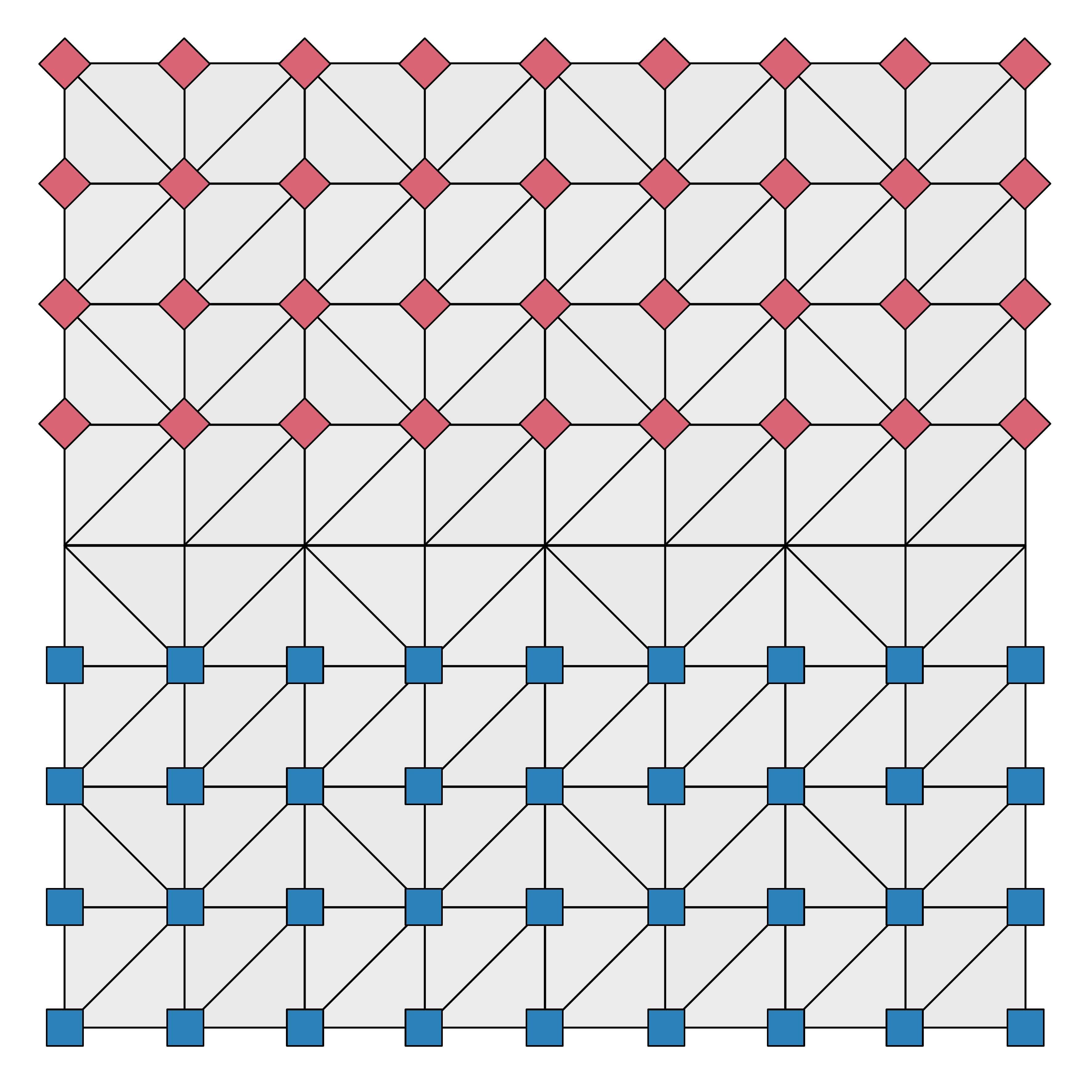}
    \caption{The ``aggregated SOS2'' biclique construction from Lemma~\ref{lemma:aggregated-sos2}. The first row depicts the sets $A^{1,k} \times \llbracket d_2 + 1 \rrbracket$ and $B^{1,k} \times \llbracket d_2 + 1 \rrbracket$ as squares and diamonds, respectively, while the second row depicts the sets $\llbracket d_1 + 1 \rrbracket \times A^{2,k}$ and $\llbracket d_1 + 1 \rrbracket \times B^{2,k}$ as squares and diamonds, respectively.}
    \label{fig:aggregated-sos2}
    \includegraphics[width=.25\linewidth]{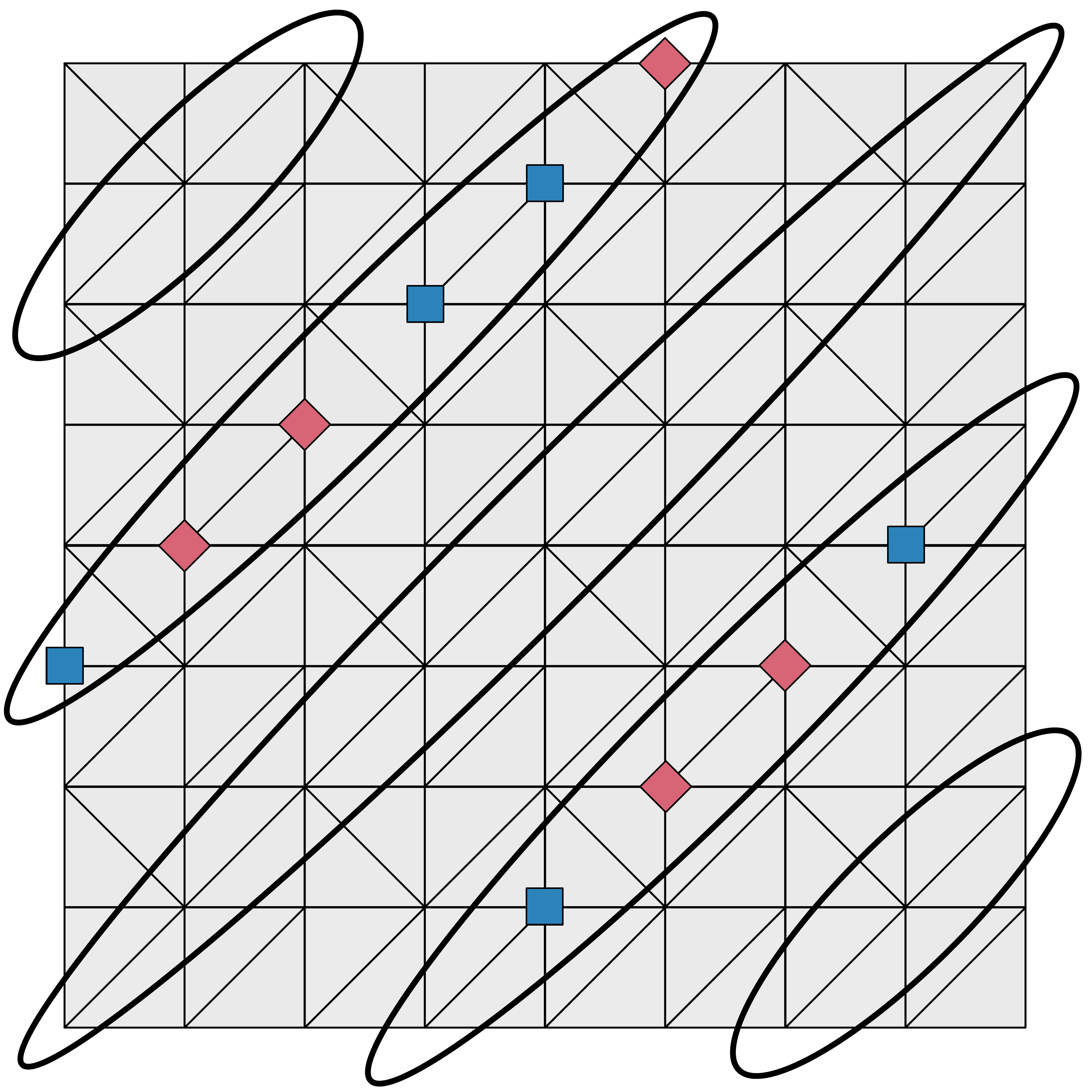}
    \includegraphics[width=.25\linewidth]{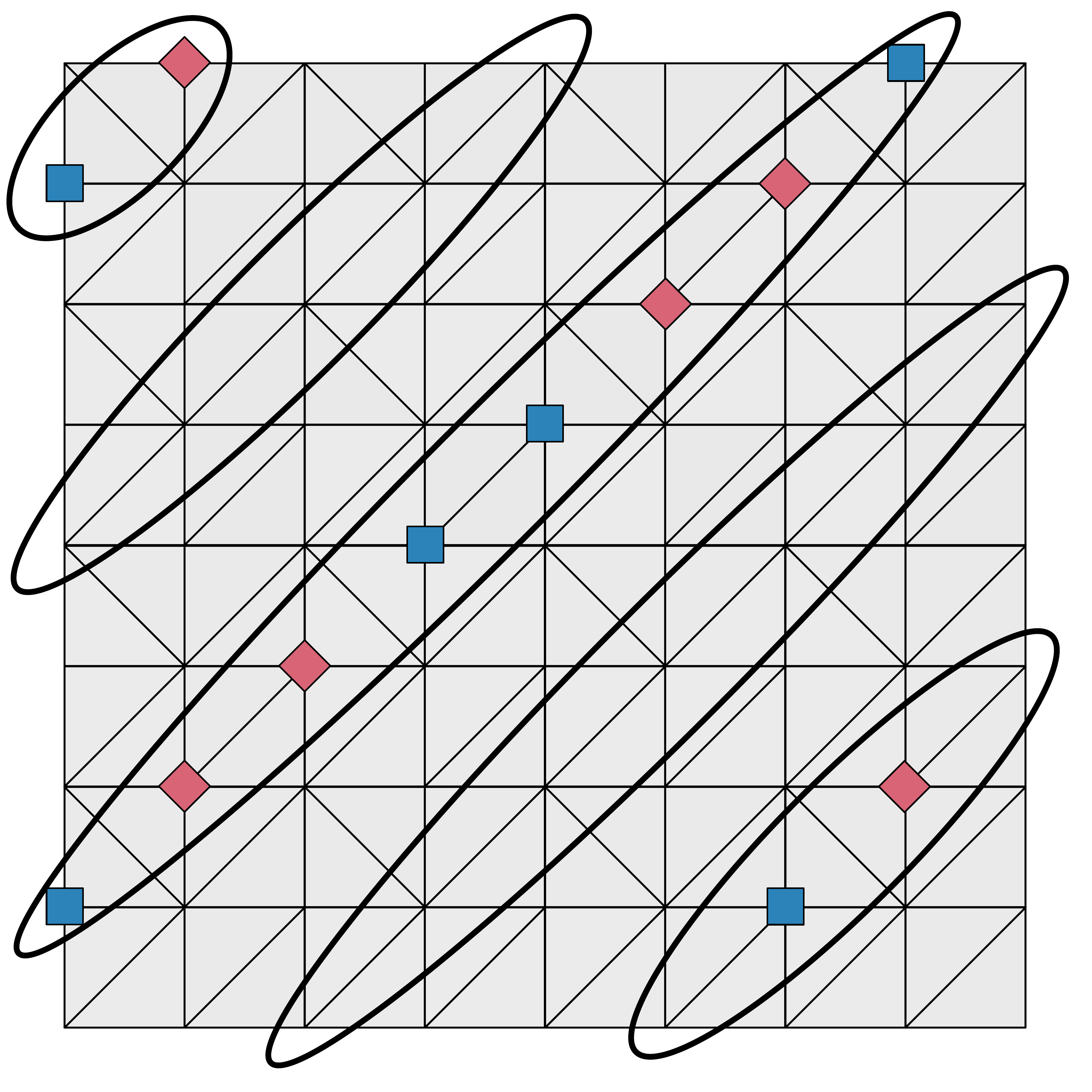}
    \includegraphics[width=.25\linewidth]{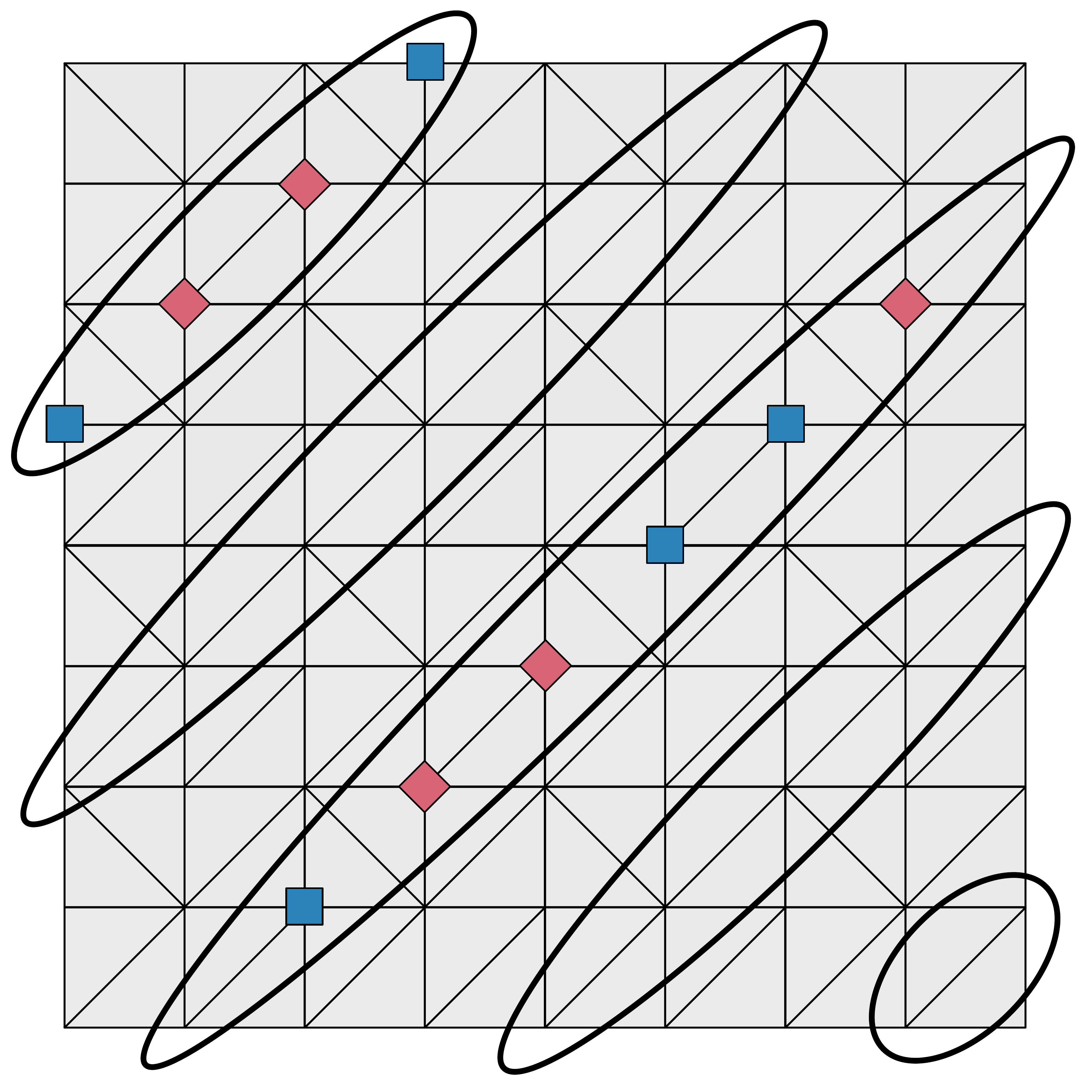} \\
    \includegraphics[width=.25\linewidth]{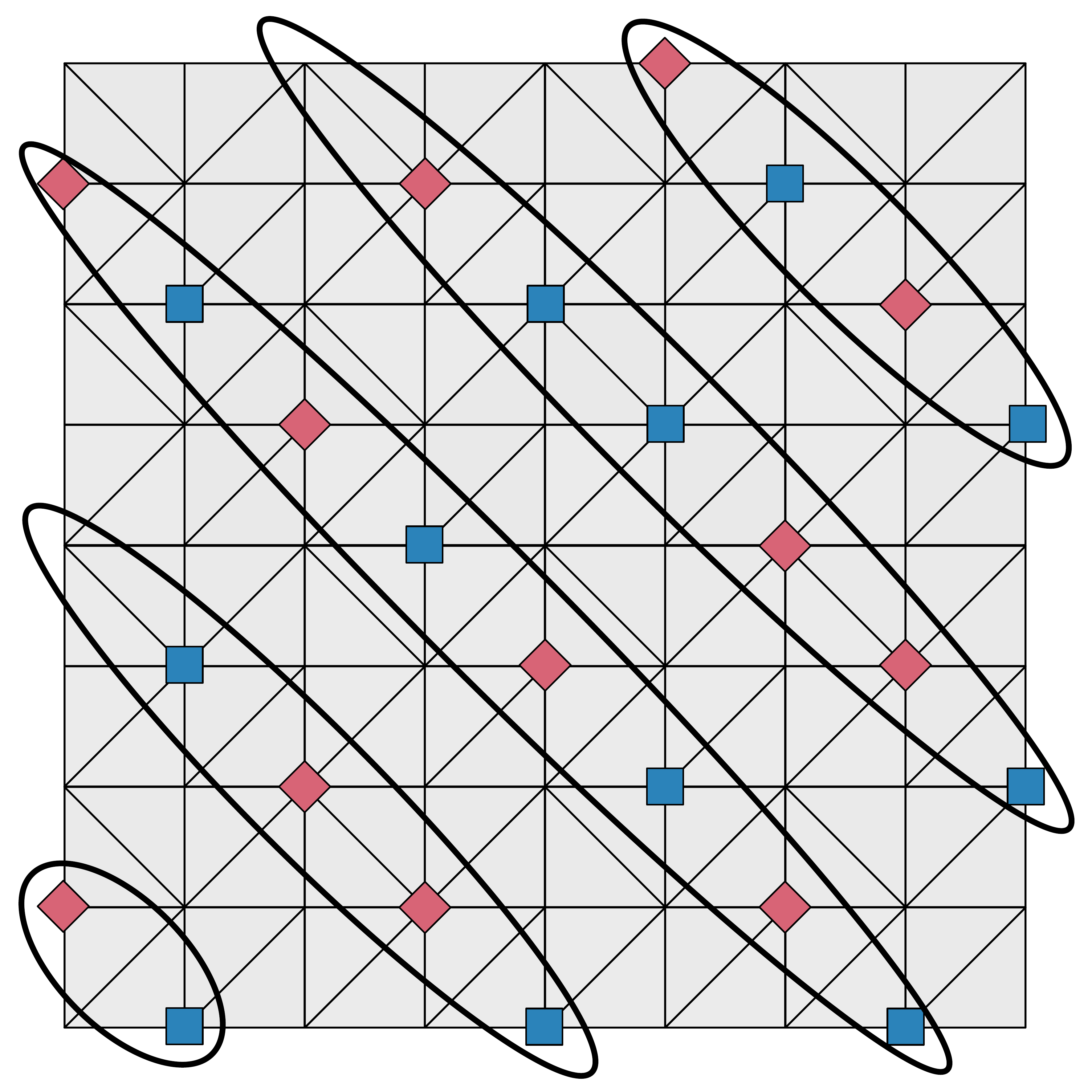}
    \includegraphics[width=.25\linewidth]{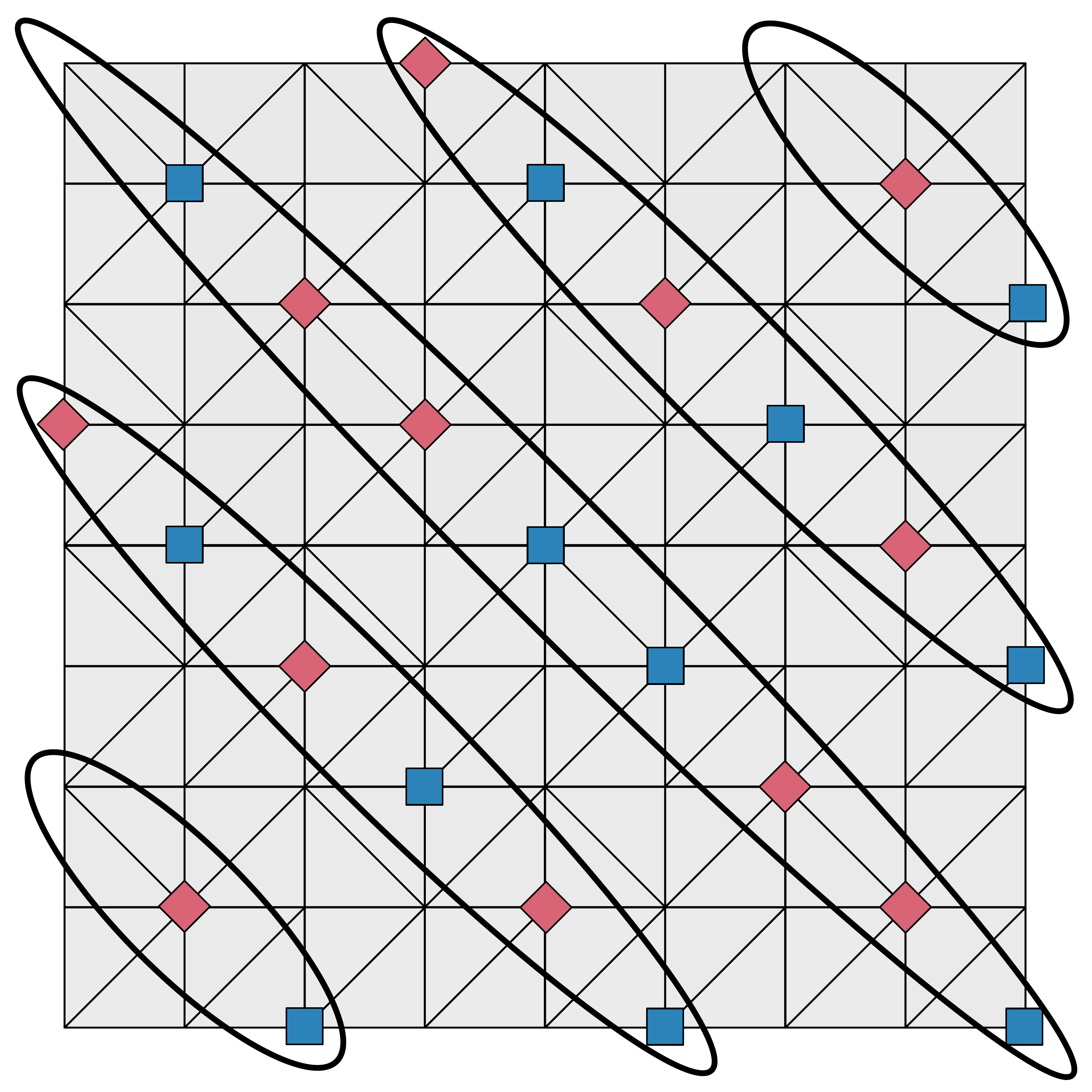}
    \includegraphics[width=.25\linewidth]{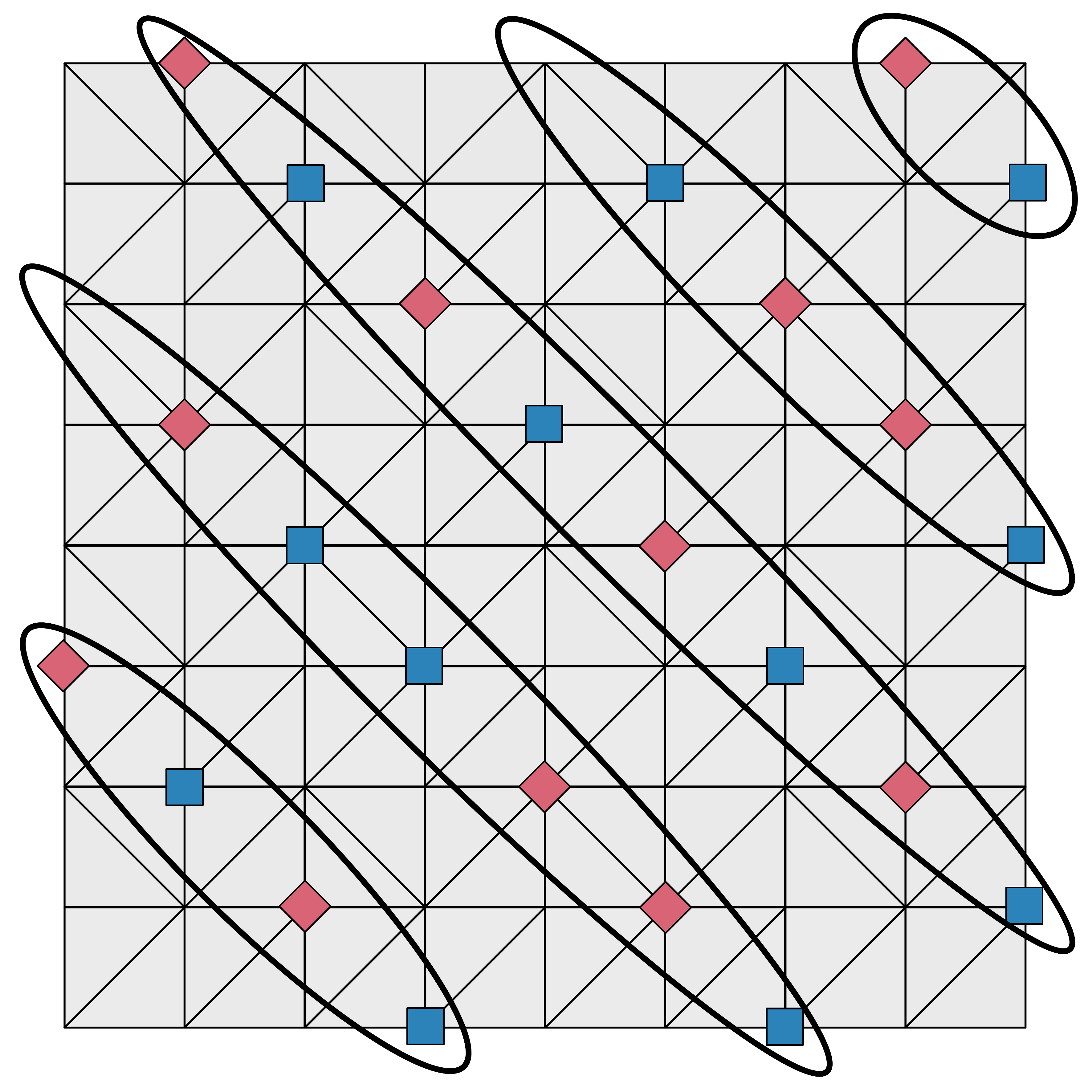}
    \caption{The diagonal and anti-diagonal biclique construction from Lemmas~\ref{lemma:diagonal} and \ref{lemma:antidiagonal}, respectively in the top and bottom rows. In each row, the sets $A^\tau$ and $B^\tau$ are depicted in squares and diamonds, respectively. As a visual aid, the diagonal/antidiagonal lines which are covered in each cell are circled.}
    \label{fig:6-stencil}
\end{figure}

\subsection{Combination of formulations}
Since our formulations for bivariate piecewise linear functions are comprised of two (aggregated) SOS2 constraints and a biclique \blue{cover} for the ``triangle selection'', we could potentially replace the independent branching formulations for the two SOS2 constraints with \emph{any} SOS2 formulation and maintain validity. For example, we can construct a \emph{hybrid} formulation for bivariate functions over a grid triangulation by applying the \ZZI{} formulation for the aggregated SOS2 constraint along the $x_1$ and the $x_2$ dimension, and the 6-stencil independent branching formulation to enforce triangle selection. However, in general the intersection of ideal formulations will not be ideal, with independent branching formulations being a notable exception. Fortunately, the following proposition (proven in Appendix~\ref{app:combining-formulations}) shows that this preservation of strength is not restricted to independent branching formulations, but holds for any intersection of ideal formulations of combinatorial disjunctive constraints.


\begin{theorem} \label{thm:combining-formulations}
    Fix $m \in \bbN$ and take:
    \begin{itemize}
        \item $U^t = \bigcup_{i=1}^{s_t} P(T^{i,t})$, where $\bigcup_{i=1}^{s_t} T^{i,t} = V$, and
        \item $\Pi^t\subseteq \mathbb{R}^{V} \times \bbR^{r_t}$ such that $\Set{(\lambda,z^t)\in \Pi^t | z^t \in \bbZ^{r_t}}$ is an ideal formulation of $U^t$
    \end{itemize}
    for each $t \in \llbracket m \rrbracket$. Then, an ideal formulation for $\bigcap_{t=1}^m U^t$ is
    \begin{equation}\label{prop:combining-formulations:form}
        \Set{
            (\lambda,z^1,\ldots,z^m) | \begin{array}{cl}
                (\lambda,z^t) \in \Pi^t \quad &\forall t \in \llbracket m \rrbracket \\
                z^t \in \bbZ^{r_t} \quad &\forall t \in \llbracket m \rrbracket
            \end{array}
        }.
        \end{equation}
\end{theorem}

\subsection{Computational experiments with bivariate piecewise linear functions} \label{ss:bivariate-computations}
To study the computational efficacy of the 6-stencil approach, we perform a computational study on a series of bicommodity transportation problems studied in Section 5.2 of \cite{Vielma:2010}. The network for each instance is fixed with 5 supply nodes and 5 demand nodes, and the objective functions are the sum of 25 concave, nondecreasing bivariate piecewise linear functions over grid triangulations with $d_1=d_2=N$ for $N \in \{4,8,16,32\}$. The triangulation of each bivariate function is generated randomly, which is the only difference from \citep{Vielma:2010}, where the Union Jack triangulation was used. To handle generic triangulations, we apply the 6-stencil formulation for triangle selection, coupled with either the \Log{}, \ZZB{}, or \ZZI{} formulation for the SOS2 constraints, taking advantage of Theorem~\ref{thm:combining-formulations} (recall that \Log{} and \LogIB{} coincide when $d$ is a power-of-two). We compare these new formulations against the \CC{}, \MC{}, and \DLog{} formulations, which readily generalize to bivariate functions. We note in passing that the \Inc{} formulation approach also generalizes to bivariate piecewise linear functions, but requires the computation of a Hamiltonian cycle~\citep{Wilson:1998}, a nontrivial task which may not be practically viable for unstructured triangulations.

\begin{table}[htpb]
    \centering
    \smaller
    \begin{tabular}{r|c|rrr:rrr}
    \multicolumn{5}{c}{} & \multicolumn{3}{c}{6-Stencil +} \\
$N$ & Metric & \texttt{MC} & \texttt{CC} & \texttt{DLog} & \Log{} & \texttt{ZZB} & \texttt{ZZI} \\ \hline
\multirow{4}{*}{4}
 & Mean (s)  & 1.4  & 1.5  & 0.9  & \textbf{0.4}  & \textbf{0.4}  & \textbf{0.4}  \\
 & Std  & 1.3  & 1.5  & 0.6  & 0.2 & 0.2 & 0.2  \\
 & Win & 0 & 0 & 0 & 29 & 31 & \textbf{40}  \\
 & Fail & \textbf{0} & \textbf{0} & \textbf{0} & \textbf{0} & \textbf{0} & \textbf{0}  \\
\hline
\multirow{4}{*}{8}
 & Mean (s)  & 39.3  & 97.2  & 12.6  & \textbf{2.7}  & 3.0  & 3.0  \\
 & Std  & 75.0  & 179.6  & 9.8  & 2.2  & 2.4  & 2.9  \\
 & Win & 0 & 0 & 0 & \textbf{51} & 17 & 32  \\
 & Fail & \textbf{0} & \textbf{0} & \textbf{0} & \textbf{0} & \textbf{0} & \textbf{0}  \\
\hline
\multirow{4}{*}{16}
 & Mean (s)  & 1370.9  & 1648.1  & 352.8  & \textbf{24.6}  & 26.5  & 35.2  \\
 & Std  & 670.4  & 360.8  & 499.4  & 24.5  & 27.4  & 40.4  \\
 & Win & 0 & 0 & 0 & \textbf{43} & 31 & 6  \\
 & Fail & 53 & 66 & 6 & \textbf{0} & \textbf{0} & \textbf{0}  \\
\hline
\multirow{4}{*}{32}
 & Mean (s)  & 1800.0  & 1800.0  & 1499.6  & \textbf{133.5}  & 167.6  & 246.5  \\
 & Std  & -  & -  & 475.2  & 162.7 & 226.7  & 306.6  \\
 & Win & 0 & 0 & 0 & \textbf{63} & 15 & 2  \\
 & Fail & 80 & 80 & 50 & \textbf{0} & \textbf{0} & 1
    \end{tabular}
    \caption{Computational results \blue{with CPLEX} for bivariate transportation problems on grids of size $N=d_1=d_2$.}
    \label{tab:bivariate-pwl}
\end{table}

In Table~\ref{tab:bivariate-pwl}, we see that the new formulations are the fastest on every instance in our test bed. For $N\in\{16,32\}$, we see an average speed-up of over an order of magnitude as compared to the \DLog{} formulation, the best of the existing approaches from the literature. We see that the \Log{} 6-stencil formulation wins a plurality or majority of instances for $N \in \{8,16,32\}$, and that the hybrid \ZZI{} 6-stencil formulation is outperformed by the hybrid \ZZB{}  6-stencil formulation by a non-trivial amount on larger instances. In particular, we highlight the largest family of instances ($N=32$), where existing methods are unable to solve 50 of 80 instances in 30 minutes or less, whereas our new formulations can solve all in a matter of minutes, on average.

For completeness, we also perform bivariate computational experiments where $N$ is not a power-of-two, now adding the \LogIB{} 6-stencil formulation as an option for the SOS2 constraints. We present the results in Appendix~\ref{app:bivariate-pwl-non-power-of-two}. Qualitatively the results are quite similar to those in Table~\ref{tab:bivariate-pwl}, although the hybrid \ZZB{} and \ZZI{} 6-stencil formulations perform slightly better on these instances, relative to the \Log{}/\LogIB{} formulations, as compared to when $N$ is a power-of-two. There is no significant difference between the \Log{} and \LogIB{} 6-stencil formulations.

\blue{
\begin{table}[htpb]
         \centering
         \smaller
         \begin{tabular}{r|c|rrr:rrr}
             \multicolumn{5}{c}{} & \multicolumn{3}{c}{6-Stencil +} \\
     $N$ & Metric & \texttt{MC} & \texttt{CC} & \texttt{DLog} & \Log{} & \texttt{ZZB} & \texttt{ZZI} \\ \hline
\multirow{4}{*}{4}
 & Mean (s)  & 1.1  & 1.8  & 0.7  & \textbf{0.3}  & \textbf{0.3}  & \textbf{0.3}  \\
 & Std  & 0.8  & 1.6  & 0.6  & 0.1 & 0.1 & 0.1 \\
 & Win & 0 & 0 & 0 & \textbf{43} & 20 & 37  \\
 & Fail & \textbf{0} & \textbf{0} & \textbf{0} & \textbf{0} & \textbf{0} & \textbf{0}  \\
\hline
\multirow{4}{*}{8}
 & Mean (s)  & 13.0  & 54.9  & 12.4  & \textbf{2.1}  & 2.3  & \textbf{2.1}  \\
 & Std  & 12.5  & 79.9  & 14.8  & 2.2  & 2.1  & 1.9 \\
 & Win & 0 & 0 & 0 & \textbf{52} & 19 & 29  \\
 & Fail & \textbf{0} & \textbf{0} & \textbf{0} & \textbf{0} & \textbf{0} & \textbf{0}  \\
\hline
\multirow{4}{*}{16}
 & Mean (s)  & 440.8  & 1154.9  & 266.7  & \textbf{16.0}  & 18.7  & 16.2  \\
 & Std  & 560.9  & 724.3  & 438.3  & 21.1  & 20.6  & 18.8 \\
 & Win & 0 & 0 & 0 & \textbf{45} & 12 & 23  \\
 & Fail & 6 & 39 & 3 & \textbf{0} & \textbf{0} & \textbf{0}  \\
\hline
\multirow{4}{*}{32}
 & Mean (s)  & 1521.6  & 1799.0  & 1291.1  & \textbf{111.6}  & 129.0  & 121.0  \\
 & Std  & 515.6  & - & 599.8  & 145.8 & 156.6  & 163.6  \\
 & Win & 0 & 0 & 0 & \textbf{48} & 10 & 22  \\
 & Fail & 56 & 79 & 38 & \textbf{0} & \textbf{0} & \textbf{0}
         \end{tabular}
         \caption{Computational results with Gurobi for bivariate transportation problems on grids of size $N=d_1=d_2$.}
         \label{tab:bivariate-gurobi}
\end{table}

We also reproduce the computational results from Table~\ref{tab:bivariate-pwl} using the Gurobi v7.0.2 solver, and include the results in Table~\ref{tab:bivariate-gurobi}. The takeaway remains the same, as the new stencil formulations are the fastest performers on every instance in the test bed. Additionally, we observe that, as in the univariate case, Gurobi is generally more efficient than CPLEX on these instances, although the difference between solvers is not nearly as dramatic as on the univariate instances. Additionally, on the largest instances ($N=32$), we observe relatively uniform behavior among the three stencil formulation variants, as opposed to CPLEX, for which we observe a substantial relative degradation of the \ZZI{} variant with respect to the other two stencil formulations.
}

\section{Computational tools for piecewise linear modeling: \texttt{PiecewiseLinearOpt}} \label{sec:computational-tools}
Throughout this work, we have investigated a number of possible formulations for optimization problems containing piecewise linear functions. The performance of these formulations can be highly dependent on latent structure in the function, and there are potentially a number of formulations one may want to try on a given instance. However, these formulations can seem quite complex and daunting to a practitioner, especially one unfamiliar with the idiosyncrasies of MIP modeling. Anecdotally, we have observed that the complexity of these formulations has driven potential users to simpler but less performant models, or to abandon MIP approaches altogether.

This gap between high-performance and accessibility is fundamental throughout optimization. One essential tool to help close the gap is the modeling language, which allows the user to express an optimization problem in a user-friendly, pseudo-mathematical style, and obviates the need to interact with the underlying optimization solver directly. Because they offer a much more welcoming experience for the modeler, algebraic modeling languages have been widely used for decades, with AMPL~\citep{Fourer:1989} and GAMS~\citep{Rosenthal:2014} being two particularly storied and successful commercial examples. JuMP~\citep{Dunning:2015a} is a recently developed open-source algebraic modeling language in the Julia programming language~\citep{Bezanson:2017} which offers state-of-the-art performance and advanced functionality, and is readily extensible.

To accompany this work, we have created \texttt{PiecewiseLinearOpt}, a Julia package that extends JuMP to offer all the formulation options discussed herein through a simple, high-level modeling interface. The package supports continuous univariate piecewise linear functions, and bivariate piecewise linear functions over grid triangulations. It supports all the formulations used in the computational experiments in this work, and can handle the construction and formulation of both structured or unstructured grid triangulations. All this complexity is hidden from the user, who can embed piecewise linear functions in their optimization problem in a single line of code with the \texttt{piecewiselinear} function.

In Figure~\ref{fig:1D-pwl-code}, we see sample code for adding piecewise linear functions to JuMP models.
\begin{figure}
\footnotesize
\begin{lstlisting}
using JuMP, PiecewiseLinearOpt, CPLEX
model = Model(solver=CplexSolver())
@variable(model, 0 <= x[1:2] <= 4)
xval = [0,1,2,3,4]
fval = [0,4,7,9,10]
z1 = piecewiselinear(model, x, xval, fval, method=:Log)
g(u,v) = 2*(u-1/3)^2 + 3*(v-4/7)^4
dx = dy = linspace(0, 1, 9)
z2 = piecewiselinear(model, x[1], x[2], dx, dy, g, method=:ZZI)
@objective(model, Min, z1 + z2)
\end{lstlisting}
\caption{\texttt{PiecewiseLinearOpt} code to set the univariate function \eqref{eqn:pwl-example} as the objective, using the \Log{} formulation.}
\label{fig:1D-pwl-code}
\end{figure}
After loading the required packages, we define the \texttt{Model} object, and add the \texttt{x} variables to it. We add the univariate function \eqref{eqn:pwl-example} to our model, specifying it in terms of the breakpoints \texttt{xval} of the domain, and the corresponding function values \texttt{fval} at these breakpoints. We call the \texttt{piecewiselinear} function, while using the \Log{} formulation. It returns a JuMP variable $\texttt{z1}$ which is constrained to be equal to $f(\texttt{x})$, and can then used anywhere in the model, e.g. in the objective function. After this, we add a bivariate piecewise linear function to our model by approximating a nonlinear function \texttt{g} on the box domain $[0,1]^2$. We use the \ZZI{} formulation along each axis $x_1$ and $x_2$; it will automatically choose the triangulation that best approximates the function values at the centerpoint of each subrectangle in the grid, and then use the 6-stencil triangle selection portion of the formulation, as the triangulation is unstructured.

To showcase the \texttt{PiecewiseLinearOpt} package in a more practical setting, we consider a share-of-choice product design problem arising in marketing (e.g. see \citep{Bertsimas:2015b,Camm:2006,Wang:2009}). We are given a product design space $x \in [0,1]^\eta$, along with with $\nu$ customer types, each with a $\lambda_i \in [0,1]$ share of the market and a \emph{partworth} (i.e. preference vector) $\beta^i \in \bbR^\eta$. For each customer type $i$, the probability of purchase is $p_i(x) = \frac{1}{1 + \exp(u_i-\beta^i \cdot x)}$, where $u_i$ is a minimum ``utility hurdle'' given by existing good products. 

Given that the true preference vectors $\beta^i$ are typically unknown, we may consider a stochastic optimization version of our problem. For each scenario $s \in \llbracket S \rrbracket$, we observe a realized preference vector $\beta^{i,s}$. Our objective is to select the product specification $x$ in order to maximize the expected number of purchases, while ensuring the product performance on each individual realized scenario is not too poor. Mathematically, we may write the optimization problem as
\begin{subequations}
\begin{alignat}{2}
    \max_{x,\mu,\bar{\mu},p,\bar{p}}\quad& \sum_{i=1}^\nu \lambda_i \bar{p}_i \\
    \text{s.t.}\quad& \bar{\mu}_i = \frac{1}{S}\sum_{s=1}^S\beta^{i,s} \cdot x \quad &\forall i \in \llbracket \nu \rrbracket \\
    & \bar{p}_i = \frac{1}{1 + \exp(u_i-\bar{\mu}_i)} \quad &\forall i \in \llbracket \nu \rrbracket \label{eqn:mean-nonlinearity} \\
    & \mu^s_i = \beta^{i,s} \cdot x \quad &\forall s \in \llbracket S \rrbracket, i \in \llbracket \nu \rrbracket \\
    & p^s_i = \frac{1}{1 + \exp(u_i-\mu^s_i)} \quad &\forall s \in \llbracket S \rrbracket, i \in \llbracket \nu \rrbracket \label{eqn:scenario-nonlinearity} \\
    & \sum_{i=1}^\nu \lambda_i p^s_i \geq C \sum_{i=1}^\nu \lambda_i \bar{p}_i \quad &\forall s \in \llbracket S \rrbracket \label{eqn:scenario-balancing} \\
    & 0 \leq x_j \leq 1 \quad &\forall j \in \llbracket \eta \rrbracket
\end{alignat}
\end{subequations}
Here $C$ is some nonnegative scaling constant, and \eqref{eqn:scenario-balancing} ensures that the expected number of purchases in a given scenario is not significantly reduced from the overall expected purchases. Our solution approach is to apply a piecewise linear approximation to the nonlinearities arising in \eqref{eqn:mean-nonlinearity} and \eqref{eqn:scenario-nonlinearity}. This can be easily accomplished with the \texttt{PiecewiseLinearOpt} package, as the code in Figure~\ref{code:share-of-choice} illustrates.

In Table~\ref{tab:share-of-choice} we report the computational performance of high-performing formulations for 18 randomly generated instances of the share-of-choice problem with a $\eta=15$ dimensional product design space, $\nu=20$ customer types, $S=12$ scenarios, scaling constant $C=0.2$, and $N=50$ pieces for each piecewise linear discretization. We observe that the \ZZI{} formulation is the best performer on the majority of instances, and substantially outperforms the \Log{} formulation, which is unable to solve any instance to optimality in 30 minutes or less. Note that for this problem the piecewise linear function will appear in both the objective and the constraints of the optimization problem.
\begin{table}
\centering
\smaller
\begin{tabular}{c|rr:rr}
Metric & \texttt{Inc} & \Log{} & \texttt{ZZB} & \texttt{ZZI} \\ \hline
Mean (s) & 880.9 & 3600.0 & 3525.5 & \textbf{776.2} \\
Std & 1202.9 & - & 316.1 & 1037.1 \\
Win & 5 & 0 & 0 & \textbf{13} \\
Fail & 2 & 18 & 17 & \textbf{1}
    \end{tabular}
    \caption{Aggregate statistics for share-of-choice problems with 50 piece discretizations.}
    \label{tab:share-of-choice}
\end{table}

\begin{figure}
\footnotesize
\begin{lstlisting}
using JuMP, Distributions,PiecewiseLinearOpt
model = Model()
@variable(model, 0 <= x[1:eta] <= 1)
@variable(model, mu[1:nu, 1:S])
@variable(model, mu_bar[1:nu])
@variable(model, p[1:nu, 1:S])
@variable(model, p_bar[1:nu])
for i in 1:nu
    @constraint(model, mu_bar[i] == 1/S * sum(dot(beta[i,s], x) for s in S))
    f(t) = 1 / (1 + exp(u[i] - t))
    @constraint(model, p_bar[i] == piecewiselinear(model, mu_bar[i], prob_min[i], prob_max[i], f)
    for s in 1:S
        @constraint(model, mu[i,s] == dot(beta[i,s], x))
        @constraint(model, p[i,s] == piecewiselinear(model, mu[i,s], scen_prob_min[i,s], scen_prob_max[i,s], f))
    end
end
for s in 1:S
    @constraint(model, sum(lambda[i]*p[i,s] for i in 1:nu) >= C * sum(lambda[i]*p_bar[i] for i in 1:nu))
end
@objective(model, Max, sum(lambda[i]*p_bar[i] for i in 1:nu))
\end{lstlisting}
\caption{\texttt{PiecewiseLinearOpt} code to solve a stochastic share-of-choice problem.}
\label{code:share-of-choice}
\end{figure}



We believe that this exemplifies the value of \texttt{PiecewiseLinearOpt}, and modeling languages more generally: it allows a user to quickly and easily write their problem as code, and then iterate as-needed to solve more quickly or to add complexity. For example, we can alter the breakpoint values in the code in Figure~\ref{code:share-of-choice} to modify the model to produce feasible solutions and upper bounds on the optimal solution. We hope that this simple computational tool will make the advanced formulations available for modeling piecewise linear functions more broadly accessible to researchers and practitioners.

%

\ACKNOWLEDGMENT{This material is based upon work supported by the National Science Foundation under Grant CMMI-1351619.}

\theendnotes

\bibliographystyle{ormsv080}
\bibliography{master.bib}

\newpage
\renewcommand{\theHsection}{A\arabic{section}}
\begin{APPENDICES}


\section{Binary reflected Gray codes, related encodings, and proof of Proposition~\ref{prop:new-sos2}} \label{app:encodings}
     \newcommand{\rev}{\operatorname{rev}}
    \newcommand{\rep}{\operatorname{rep}}
The following straightforward lemma gives a recursive construction for $K^r$, $C^r$, and $Z^r$.
    \begin{lemma} \label{lem:encodings}
    $K^1=C^1=Z^1\defeq (0,1)^T$, and for $r \in \bbN$ (and $d=2^r$):
    \[
        \setlength\arraycolsep{5pt}K^{r+1}\defeq \begin{pmatrix} K^r & \mathbf{0}^{d} \\ \rev(K^r) & \mathbf{1}^{d} \end{pmatrix},\quad\quad C^{r+1}\defeq \begin{pmatrix}C^r & \mathbf{0}^{d} \\C^r + \mathbf{1}^{d}\otimes C^r_{d} & \mathbf{1}^{d} \end{pmatrix},\quad\text{and}\quad Z^{r+1}\defeq \begin{pmatrix}Z^r & \mathbf{0}^{d} \\Z^r  & \mathbf{1}^{d} \end{pmatrix},
    \]
    where $\mathbf{0}^{r}, \mathbf{1}^{r} \in \bbR^r$ are the vectors with all components equal to $0$ or $1$, respectively, $u\otimes v=uv^T\in\mathbb{R}^{m\times n}$ for any $u\in \mathbb{R}^m$ and $v\in \mathbb{R}^n$, and $\rev(A)$ reverses the rows of the matrix $A$.
    \end{lemma}
    \proof{\textbf{Proof of Proposition~\ref{prop:new-sos2}}}
        First, we observe that $K^r, Z^r \in \{0,1\}^{d \times r}$ and that $\scrA$ is an invertible linear map. Therefore, for each $r \in \bbN$, $K^r$, $C^r$, and $Z^r$ are in convex position. Additionally, as $K^r$ and $Z^r$ are binary matrices, they are trivially hole-free. Additionally, the hole-free property is inherited by $C^r$ from $Z^r$ since $\scrA$ is invertible and linear, and both  $\scrA$ and  $\scrA^{-1}$ are unimodular ($\scrA(w) \in \bbZ^r$ if and only if $w \in \bbZ^r$).        

        Now the result is direct from Proposition~\ref{prop:sos2}, as $\{c^i \equiv C^r_{i+1}-C^r_i\}_{i=1}^{d-1} = \{{\bf e}^k\}_{k=1}^r$, where ${\bf e}^k$ is the canonical unit vector with support on component $k$, and the inverse of $\scrA$ is $\scrA^{-1}(y)_k = y_k + \sum_{\ell=k+1}^r 2^{\ell-k-1}y_\ell$ for each $k \in \llbracket r \rrbracket$. Formulations \eqref{eqn:integer-zigzag-formulation} and \eqref{eqn:binary-zigzag-formulation} correspond to encodings $C^r$ and $Z^r$, respectively.
    \Halmos\endproof
    
\section{An example where \Log{} and \LogIB{} do not coincide} \label{app:log-logib}

Consider the SOS2 instance with $d=3$ segments. The \Log{} formulation is all $(\lambda,y) \in \Delta^4 \times \{0,1\}^2$ such that
\begin{subequations}\label{eqn:log-non-power-of-two-example}
\begin{alignat}{2}
    \lambda_3 + \lambda_4 &\leq y_1, \quad\quad \lambda_2 + \lambda_3 + \lambda_4 &\geq y_1 \label{eqn:log-non-power-of-two-example-1} \\
    \lambda_4 &\leq y_2, \quad\quad \lambda_3 + \lambda_4 &\geq y_2.
\end{alignat}
\end{subequations}
This follows from Proposition~\ref{prop:sos2}, after observing that the spanning hyperplanes needed are given by the directions $b^1=(1,0)$ and $b^2=(0,1)$.

The \LogIB{} formulation is all $(\lambda,y) \in \Delta^4 \times \{0,1\}^2$ such that
\begin{subequations}\label{eqn:logib-non-power-of-two-example}
\begin{alignat}{2}
    \lambda_3 &\leq y_1, \quad\quad \lambda_2 + \lambda_3 + \lambda_4 &\geq y_1 \label{eqn:logib-non-power-of-two-example-1} \\
    \lambda_4 &\leq y_2, \quad\quad \lambda_3 + \lambda_4 &\geq y_2.
\end{alignat}
\end{subequations}
This follows from Proposition~\ref{bicliqueprop}, after observing that a biclique cover for the conflict graph edge set $E = \{\{1,3\},\{1,4\},\{2,4\}\}$ is $A^1=\{3\}$, $B^1=\{1\}$, $A^2=\{4\}$, and $B^2=\{1,2\}$. We then transform the formulation using the equation $\lambda_1+\lambda_2+\lambda_3+\lambda_4=1$ to present the \LogIB{} formulation in a way analogous to \eqref{eqn:log-non-power-of-two-example}, where we can observe that the first inequality in \eqref{eqn:log-non-power-of-two-example-1} differs from the first inequality in \eqref{eqn:logib-non-power-of-two-example-1}.

\section{8-segment piecewise linear function formulation branching} \label{app:8-segment}
Consider the univariate piecewise linear function $f : [0,8] \to \bbR$ given by
\begin{equation} \label{eqn:pwl-example-8}
    f(x) = \begin{cases}
        8x      & 0 \leq x \leq 1 \\
        7x +  1 & 1 \leq x \leq 2 \\
        6x +  3 & 2 \leq x \leq 3 \\
        5x +  6 & 3 \leq x \leq 4 \\
        4x + 10 & 4 \leq x \leq 5 \\
        3x + 15 & 5 \leq x \leq 6 \\
        2x + 21 & 6 \leq x \leq 7 \\
         x + 28 & 7 \leq x \leq 8.
    \end{cases}
\end{equation}
The corresponding \LogIB{}/\Log{} formulation is
\begin{subequations} \label{eqn:log-8seg}
\begin{gather}
    x = \lambda_2 + 2\lambda_3 + 3\lambda_4 + 4\lambda_5 + 5\lambda_6 + 6\lambda_7 + 7\lambda_8 + 8\lambda_9 \\
    z = 8\lambda_2 + 15\lambda_3 + 21\lambda_4 + 26\lambda_5 + 30\lambda_6 + 33\lambda_7 + 35\lambda_8 + 36\lambda_9 \\
    \lambda_3 + \lambda_7 \leq y_1 \leq \lambda_2 + \lambda_3 + \lambda_4 + \lambda_6 + \lambda_7 + \lambda_8 \\
    \lambda_4 + \lambda_5 + \lambda_6 \leq y_2 \leq \lambda_3 + \lambda_4 + \lambda_5 + \lambda_6 + \lambda_7 \\
    \lambda_6 + \lambda_7 + \lambda_8 + \lambda_9 \leq y_3 \leq \lambda_5 + \lambda_6 + \lambda_7 + \lambda_8 + \lambda_9 \\
    (\lambda,y) \in \Delta^9 \times \{0,1\}^3,
\end{gather}
\end{subequations}
and the corresponding \ZZI{} formulation is
\begin{subequations} \label{eqn:zigzag-8seg}
\begin{gather}
    x = \lambda_2 + 2\lambda_3 + 3\lambda_4 + 4\lambda_5 + 5\lambda_6 + 6\lambda_7 + 7\lambda_8 + 8\lambda_9 \\
    z = 8\lambda_2 + 15\lambda_3 + 21\lambda_4 + 26\lambda_5 + 30\lambda_6 + 33\lambda_7 + 35\lambda_8 + 36\lambda_9 \\
    \lambda_3 + \lambda_4 + 2\lambda_5 + 2\lambda_6 + 3\lambda_7 + 3\lambda_8 + 4\lambda_9 \leq y_1 \leq \lambda_2 + \lambda_3 + 2\lambda_4 + 2\lambda_5 + 3\lambda_6 + 3\lambda_7 + 4\lambda_8 + 4\lambda_9 \\
    \lambda_4 + \lambda_5 + \lambda_6 + \lambda_7 + 2\lambda_8 + 2\lambda_9 \leq y_2 \leq \lambda_3 + \lambda_4 + \lambda_5 + \lambda_6 + 2\lambda_7 + 2\lambda_8 + 2\lambda_9 \\
    \lambda_6 + \lambda_7 + \lambda_8 + \lambda_9 \leq y_3 \leq \lambda_5 + \lambda_6 + \lambda_7 + \lambda_8 + \lambda_9 \\
    (\lambda,y) \in \Delta^9 \times \bbZ^3
\end{gather}
\end{subequations}
In Table~\ref{tab:pwl-example-8}, we show statistics for the relaxations of the both. We observe that the \ZZI{} formulation yields more balanced branching.
\begin{table}[htpb]
    \centering
    \smaller
    \begin{tabular}{c|cc||cc|cc|cc|cc}
         Statistic & \Log{} $0\downarrow$ & \Log{} $1\uparrow$ & \ZZI{} $0\downarrow$ & \ZZI{} $1\uparrow$ & \ZZI{} $1\downarrow$ & \ZZI{} $2\uparrow$ & \ZZI{} $2\downarrow$ & \ZZI{} $3\uparrow$ & \ZZI{} $3\downarrow$ & \ZZI{} $4\uparrow$ \\ \hline
         Volume & 41 & 17 & 0 & 38.5 & 11.5 & 27 & 27 & 11.5 & 38.5 & 0 \\
         Strengthened Prop. & 0 & 1 & 1 & 0.25 & 0.75 & 0.5 & 0.5 & 0.75 & 0.25 & 1
    \end{tabular}
    \caption{Metrics for each possible branching decision on $z_1$ for \Log{} and \ZZI{} applied to \eqref{eqn:pwl-example-8}.}
    \label{tab:pwl-example-8}
\end{table}

\begin{figure}[htpb]
    \centering

    \includegraphics[width=.24\linewidth]{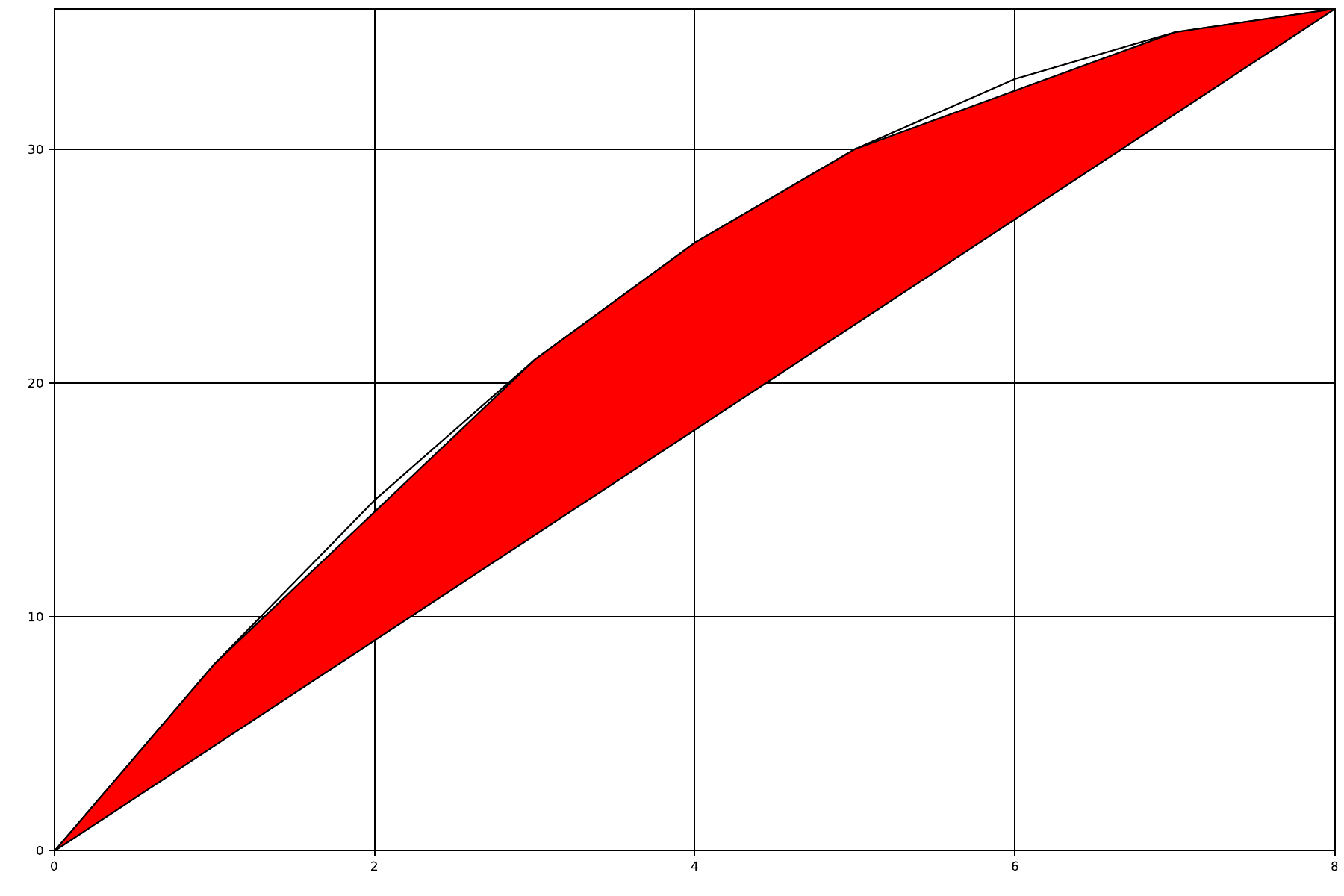}
    \includegraphics[width=.24\linewidth]{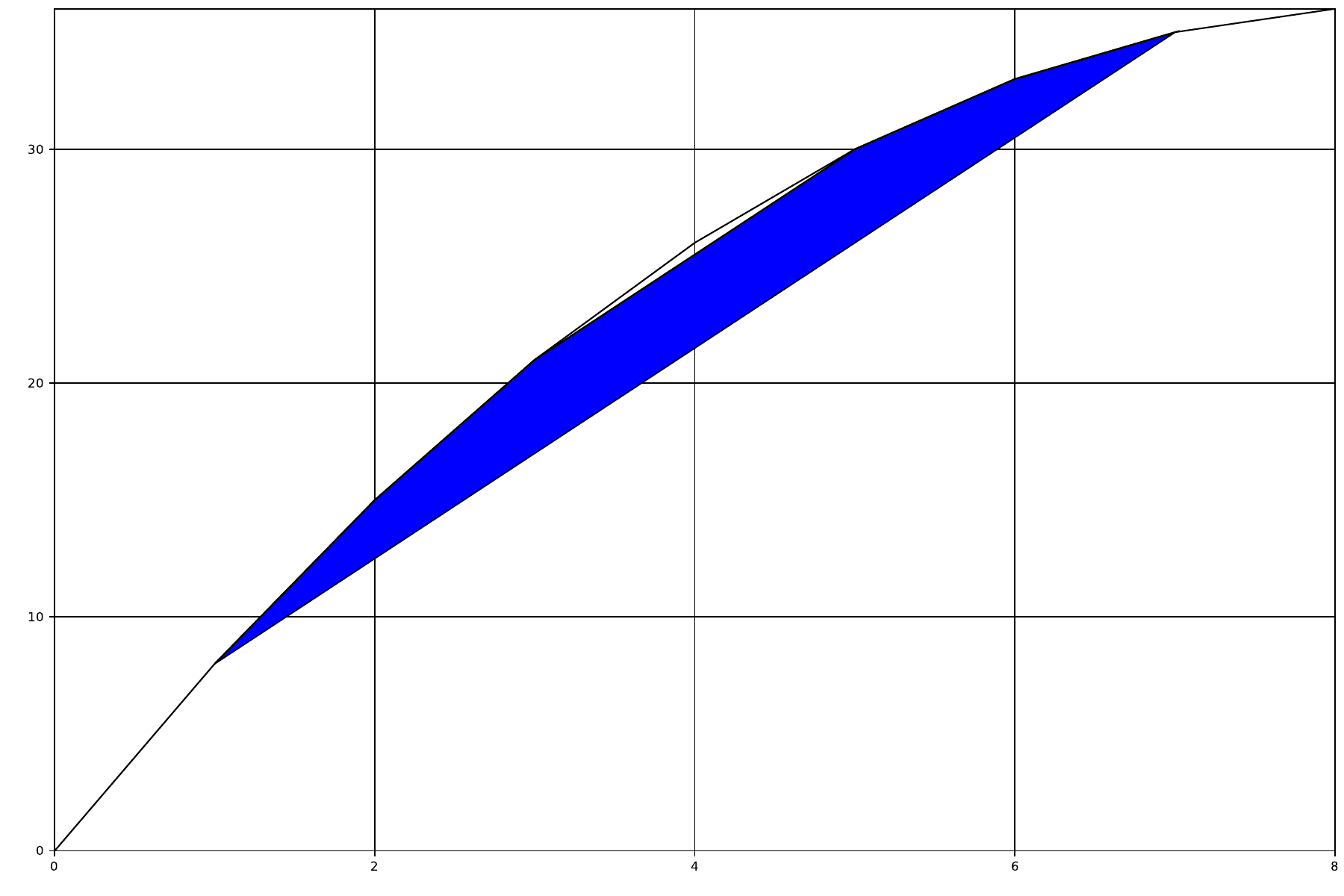}

    \caption{Feasible region in the $(x,z)$-space for the \Log{} formulation \eqref{eqn:log-8seg} after: down-branching $y_1 \leq 0$ \textbf{(left)}, and up-branching $y_1 \geq 1$ \textbf{(right)}.}

    \includegraphics[width=.24\linewidth]{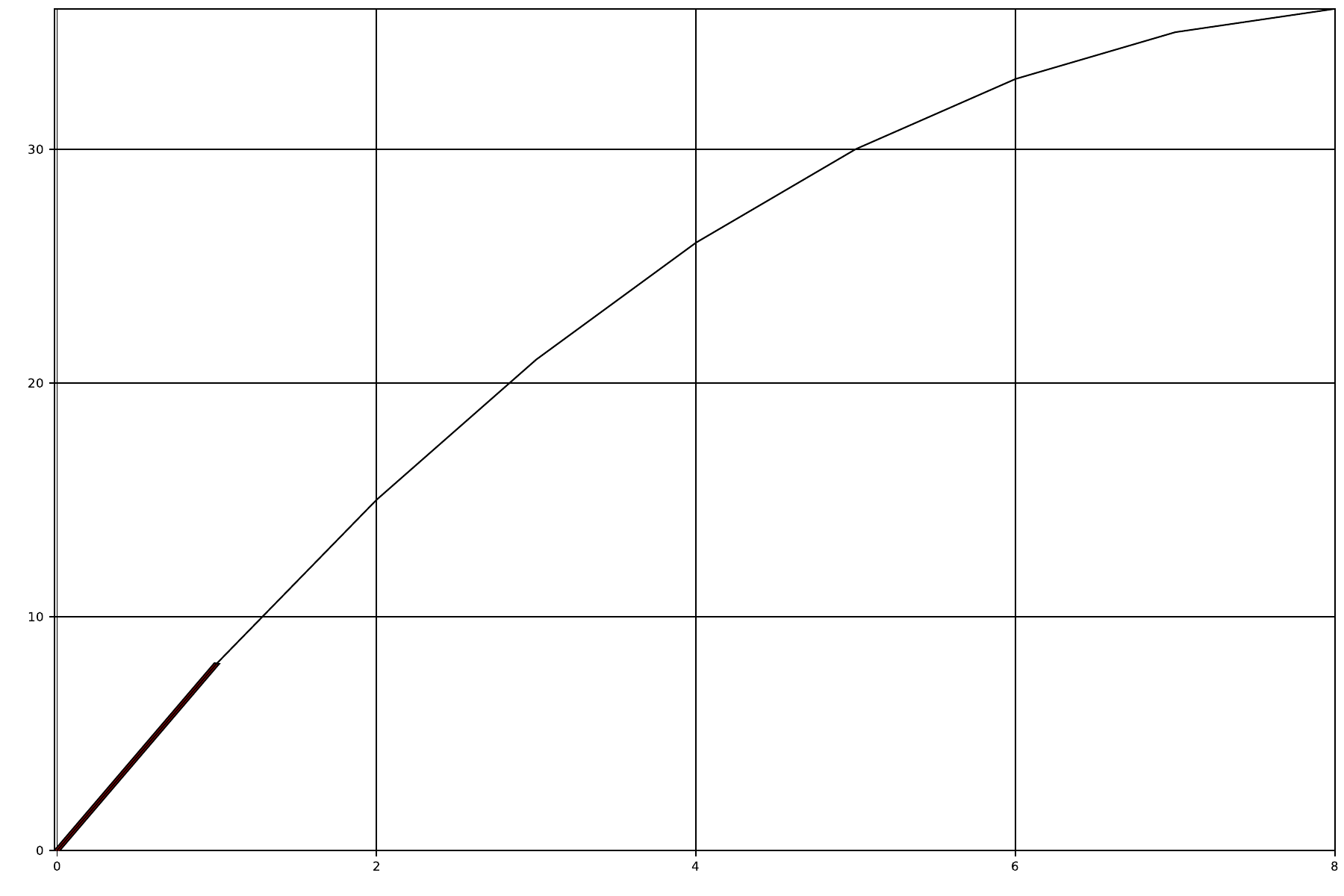}
    \includegraphics[width=.24\linewidth]{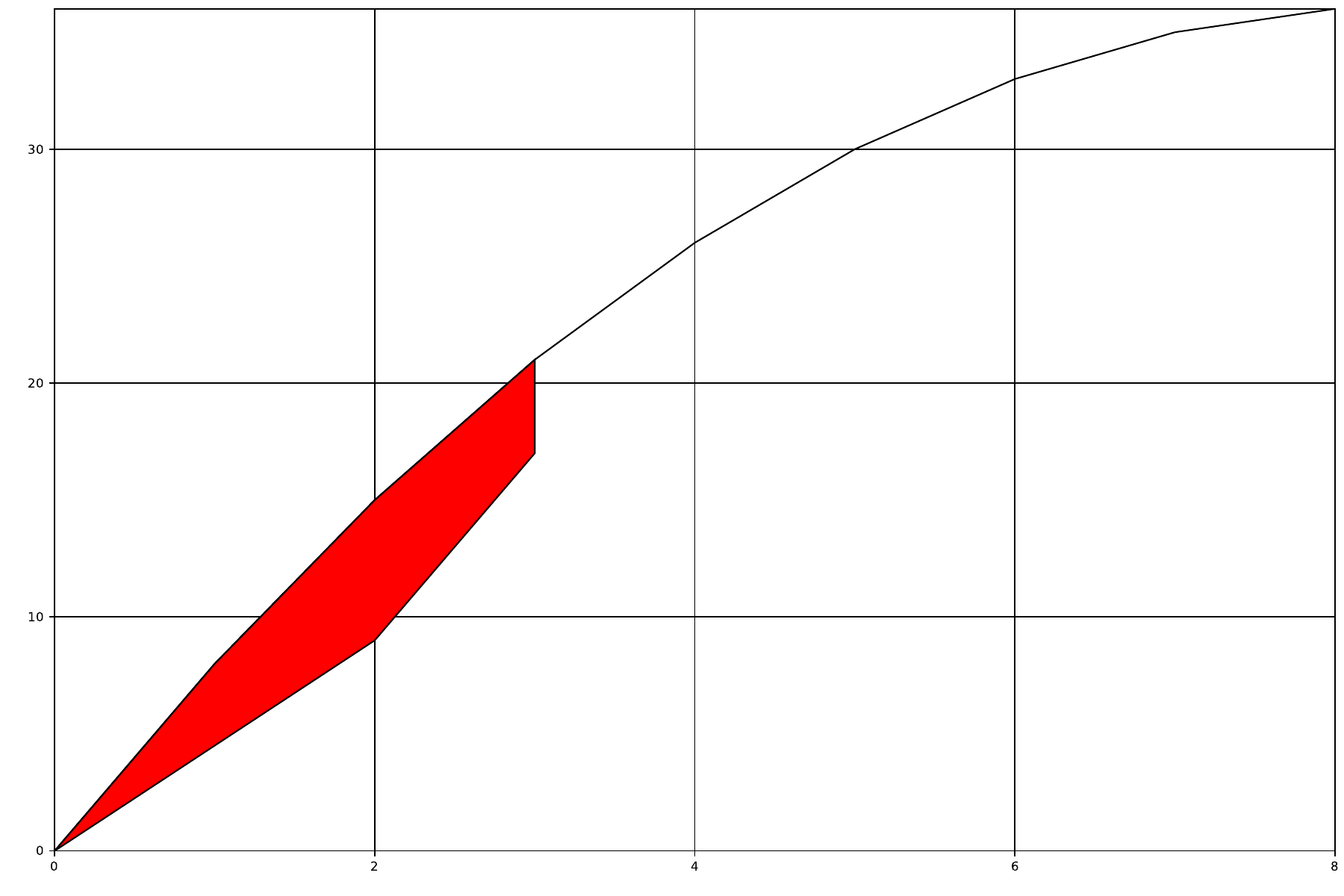}
    \includegraphics[width=.24\linewidth]{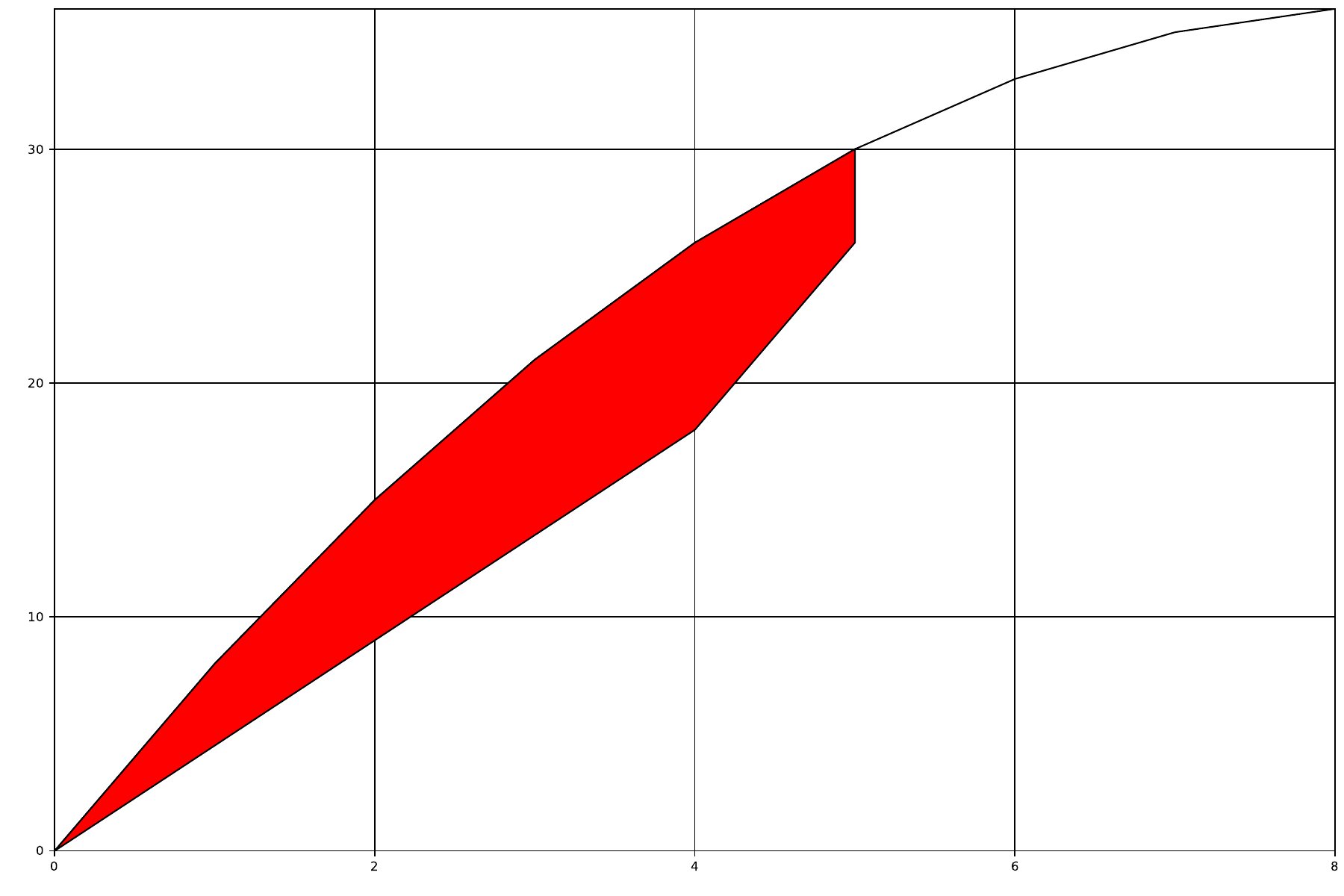}
    \includegraphics[width=.24\linewidth]{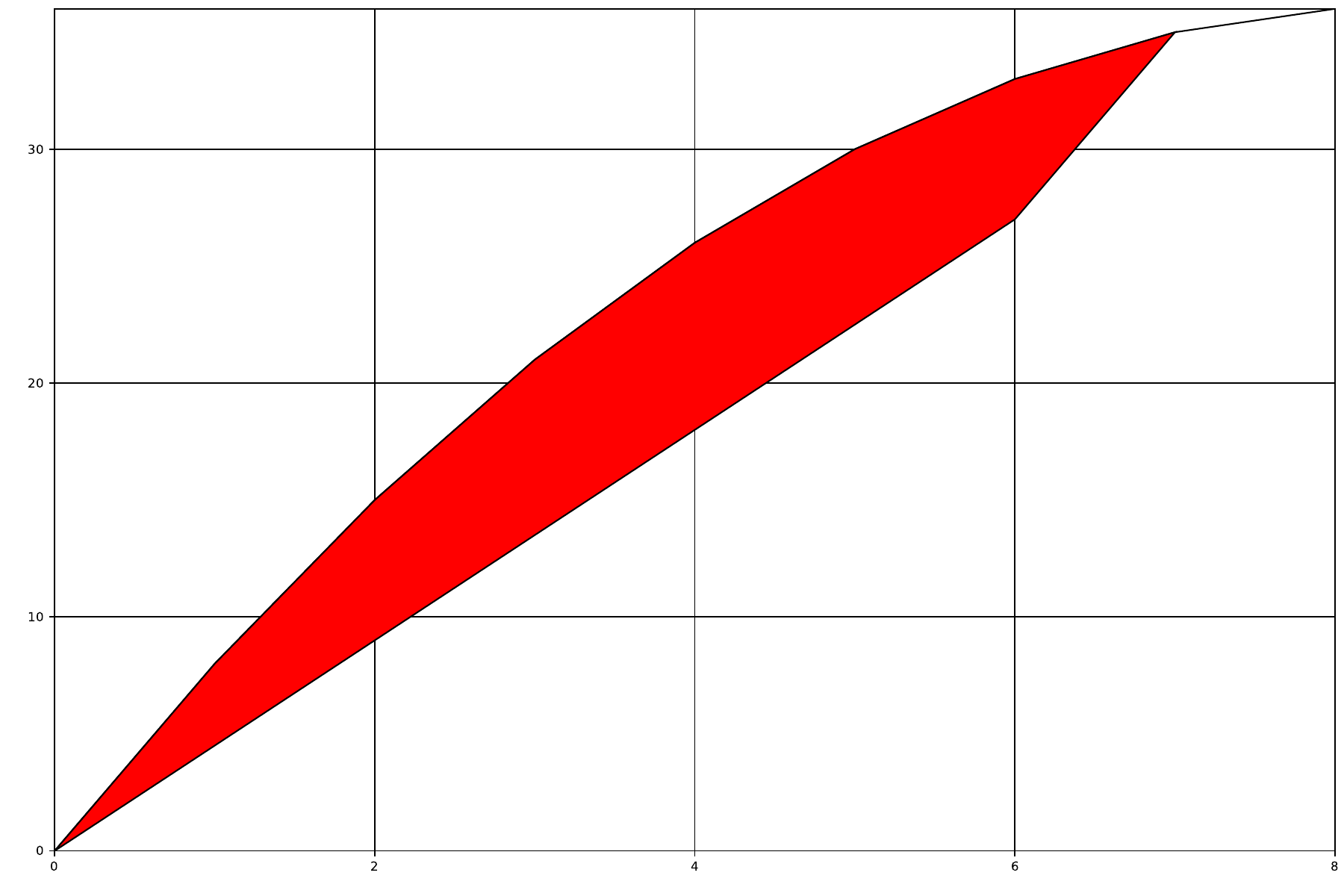} \\
    \includegraphics[width=.24\linewidth]{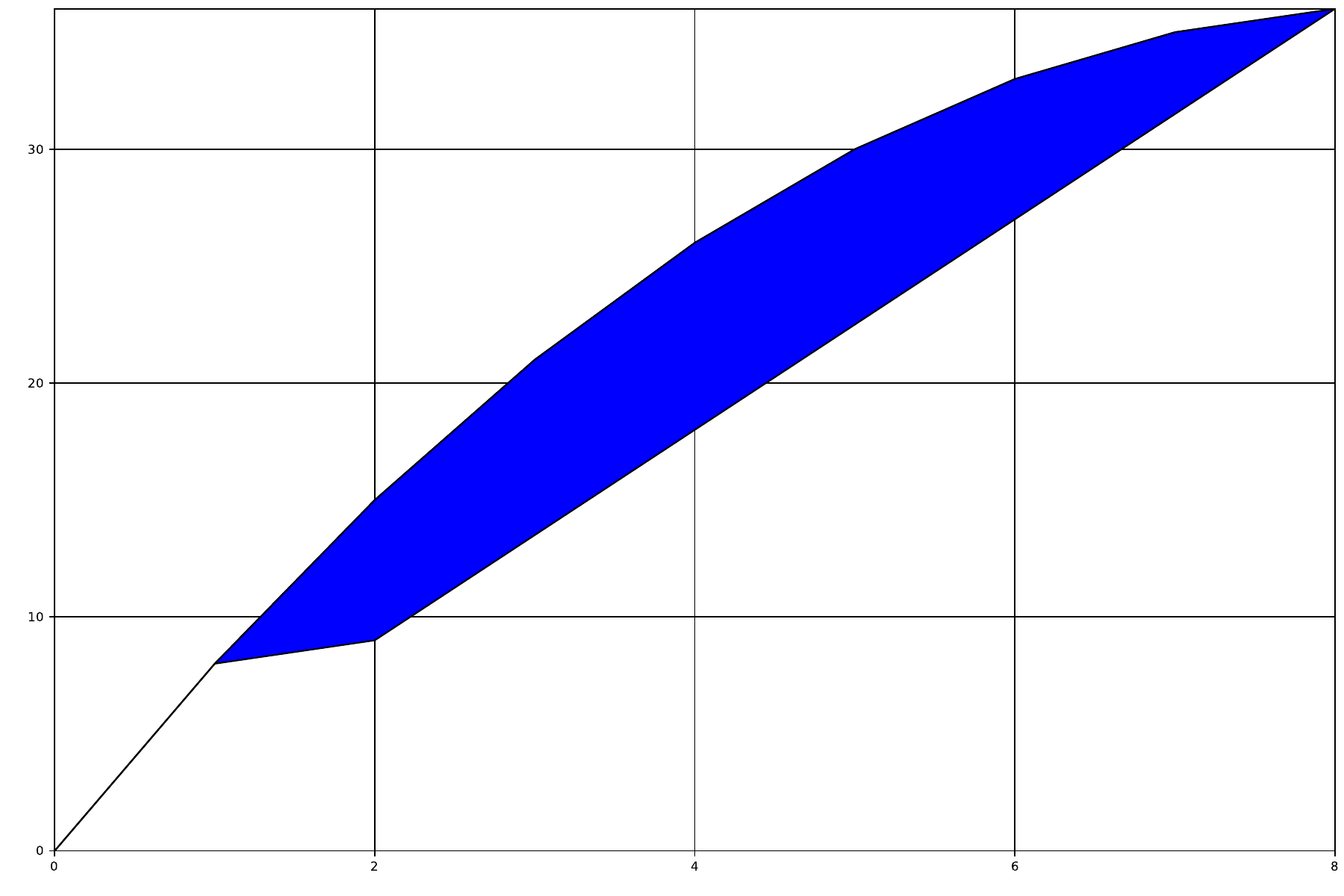}
    \includegraphics[width=.24\linewidth]{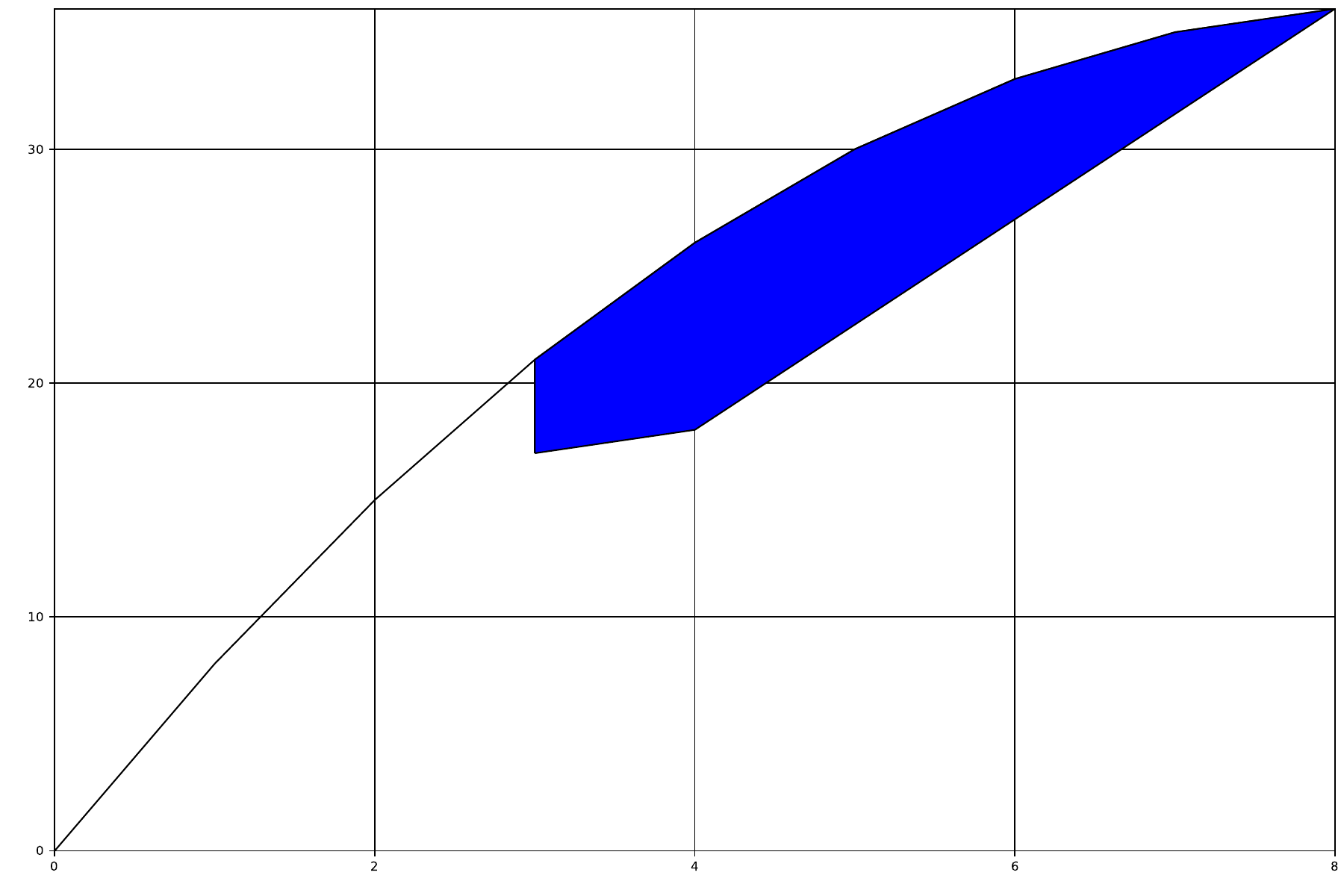}
    \includegraphics[width=.24\linewidth]{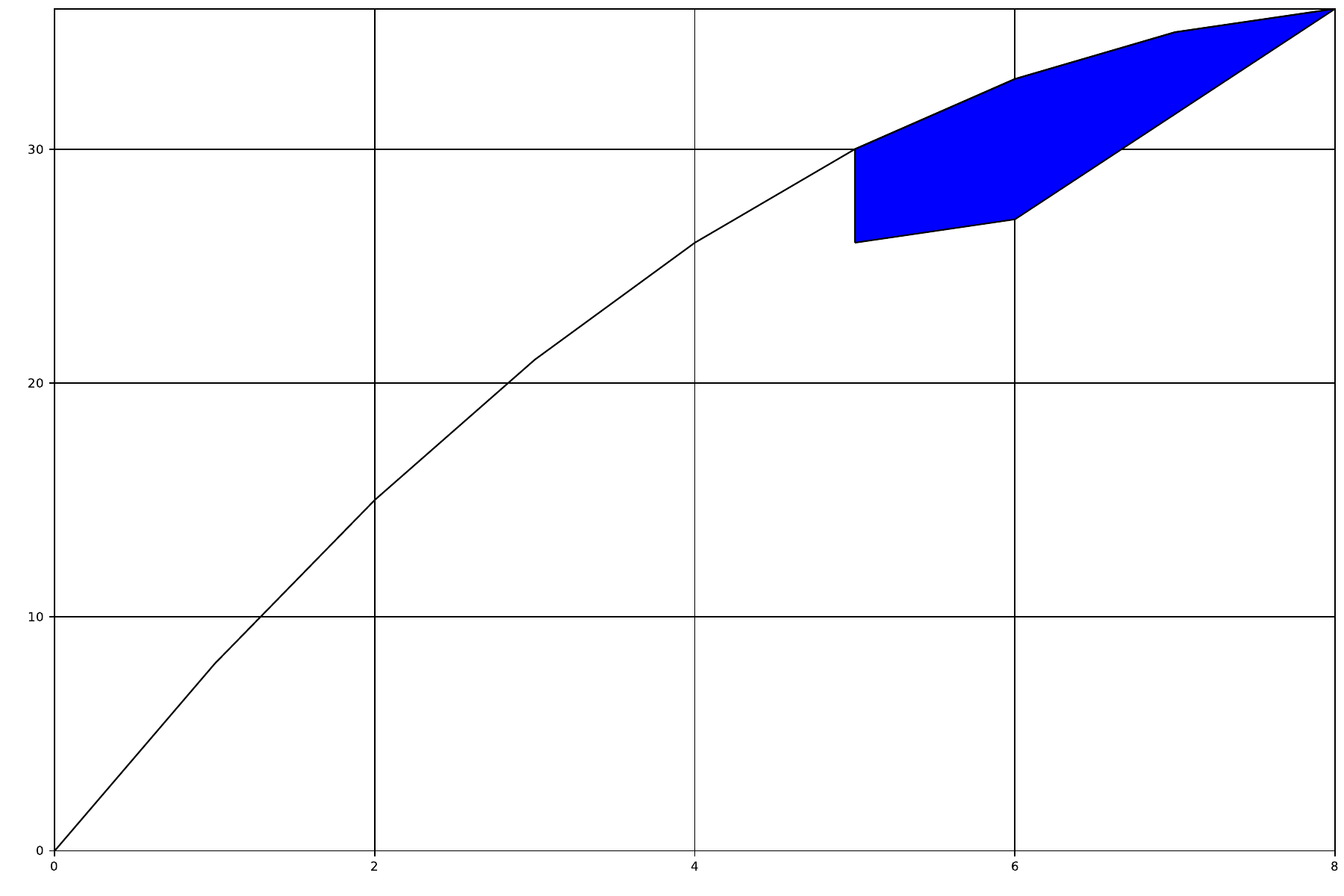}
    \includegraphics[width=.24\linewidth]{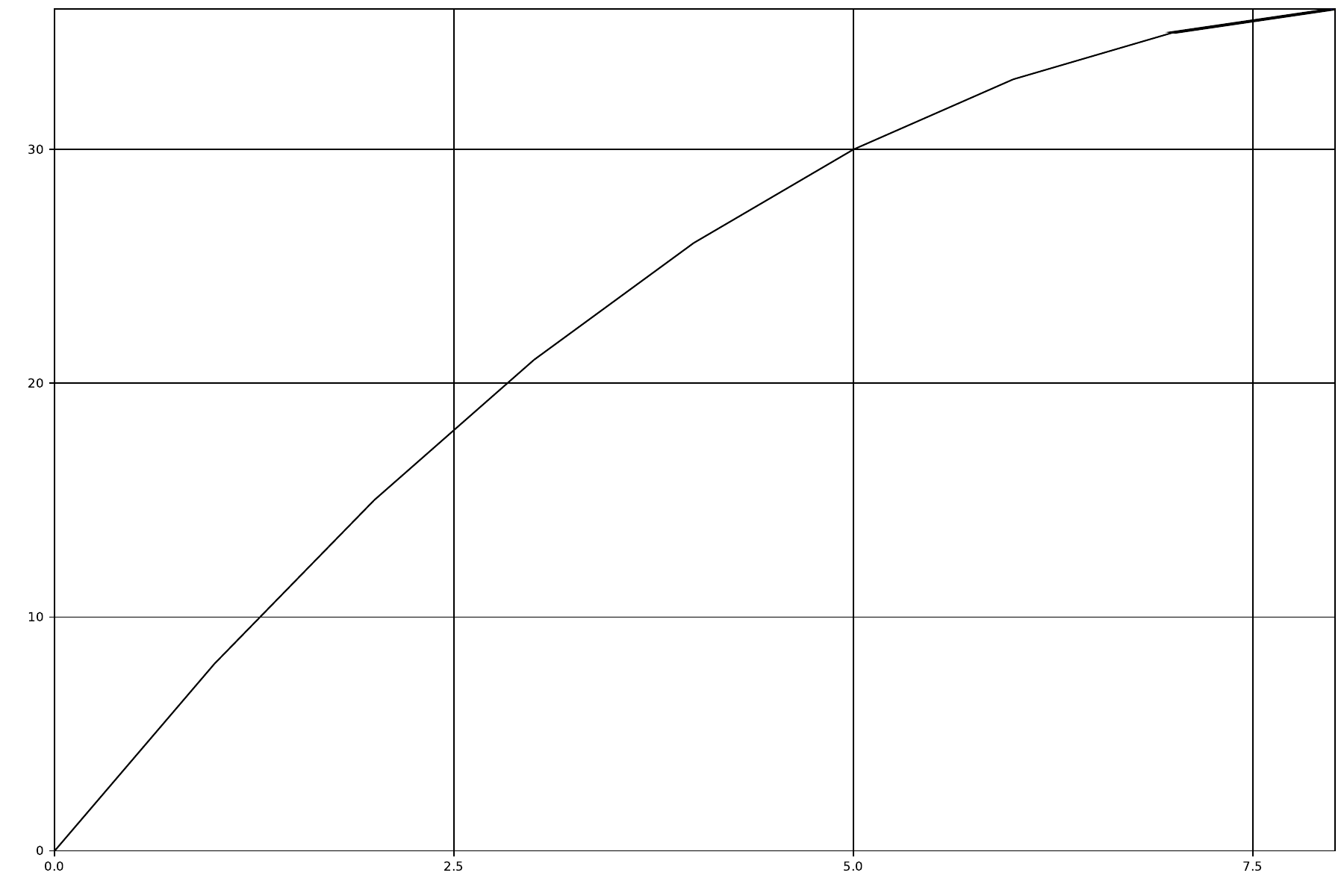}
    \caption{Feasible region in the $(x,z)$-space for the \ZZI{} formulation \eqref{eqn:zigzag-8seg} after: \textbf{(Top first column)} down-branching on $y_1 \leq 0$, \textbf{(Bottom first column)} up-branching on $y_1 \geq 1$; \textbf{(Top second column)} down-branching on $y_1 \leq 1$, \textbf{(Bottom second column)} up-branching on $y_1 \geq 2$; \textbf{(Top third column)} down-branching on $y_1 \leq 2$, \textbf{(Bottom third column)} up-branching on $y_1 \geq 3$; \textbf{(Top fourth column)} down-branching on $y_1 \leq 3$, and \textbf{(Bottom fourth column)} up-branching on $y_1 \geq 4$.}
    \label{fig:pwl-8seg-branching}
\end{figure}

\blue{
\section{Algorithm to compute bicliques for anti-diagonal nearby edges} \label{app:anti-diagonal}

Algorithm~\ref{alg:anti-diagonal} presents the analog procedure to Algorithm~\ref{alg:diagonal} for the anti-diagonal case. 
\begin{algorithm}[htpb]
    \caption{Computing bicliques that cover anti-diagonal triangle selection edges.}
    \begin{algorithmic}[1]
    \Require{Integer $\tau \in \{0,1,2\}$.}
    \Procedure{AntiDiagonalBicliques}{$\tau$,E}
        \For{$\kappa \gets \tau, \tau+3, \ldots, d_1-1$}
            \Let{$\phi$}{\texttt{true}}
            \For{$j \gets 1, \ldots, d_2$}
                \Let{$i$}{$d_1 + 1 - j - \kappa$}
                \If{$1 \leq i \leq d_1$ and $\{(i,j),(i+1,j+1)\} \in E$}
                    \If{$\phi$}
                        \Let{$A$}{$A \cup (i+1,j)$}
                        \Let{$B$}{$B \cup (i,j+1)$}
                    \Else
                        \Let{$A$}{$A \cup (i,j+1)$}
                        \Let{$B$}{$B \cup (i+1,j)$}
                    \EndIf
                    \Let{$\phi$}{$\neg\phi$}
                \EndIf
            \EndFor
        \EndFor
        \For{$\kappa \gets 3-\tau, 6-\tau, \ldots, d_2-1$}
            \Let{$\phi$}{\texttt{true}}
            \For{$j \gets (1+\kappa), \ldots, d_2$}
                \Let{$i$}{$d_1 + 1 - j + \kappa$}
                \If{$1 \leq i \leq d_1$ and $\{(i,j),(i+1,j+1)\} \in E$}
                    \If{$\phi$}
                        \Let{$A$}{$A \cup (i+1,j)$}
                        \Let{$B$}{$B \cup (i,j+1)$}
                    \Else
                        \Let{$A$}{$A \cup (i,j+1)$}
                        \Let{$B$}{$B \cup (i+1,j)$}
                    \EndIf
                    \Let{$\phi$}{$\neg\phi$}
                \EndIf
            \EndFor
        \EndFor
        \State \Return{$(A,B)$}
    \EndProcedure
    \end{algorithmic}
    \label{alg:anti-diagonal}
\end{algorithm}
}

\section{Proof of Theorem~\ref{thm:combining-formulations}} \label{app:combining-formulations}
     \proof{Proof of Theorem~\ref{thm:combining-formulations}}

        For simplicity, assume w.l.o.g. that $V = \llbracket n \rrbracket$. Let \[
            \Pi = \Set{(\lambda,z^1,\ldots,z^m)\in \mathbb{R}^{n+\sum_{i=1}^m r_i} | (\lambda,z^t) \in \Pi^t \:\: \forall t \in \llbracket m \rrbracket}
        \]
       be the LP relaxation of \eqref{prop:combining-formulations:form}. Because the original formulations are ideal (and therefore also sharp), we have
        \[
            \Proj_\lambda(\Pi) = \bigcap_{t=1}^m \Proj_\lambda(\Pi^t) = \bigcap_{t=1}^m \Conv(U^t) \subseteq \Delta^n = \Conv\left(\bigcap_{t=1}^m U^t \right),
        \]
        and hence \eqref{prop:combining-formulations:form} is sharp, as $\Proj_\lambda(\Pi)= \Delta^n$.

        To show \eqref{prop:combining-formulations:form} is also ideal, consider any point $(\hat{\lambda},\hat{z}^1,\ldots,\hat{z}^m) \in \Pi$. First, we show that if this point is extreme, then $\hat{\lambda} = {\bf e}^v$ for some $v \in \llbracket n \rrbracket$. Consider some point where $\hat{\lambda}$ is fractional; w.l.o.g., presume that $0 < \hat{\lambda}_1, \hat{\lambda}_2 < 1$. Define $\lambda^+ \defeq \hat{\lambda} + \epsilon {\bf e}^1 - \epsilon {\bf e}^2$ and $\lambda^- \defeq \hat{\lambda} - \epsilon {\bf e}^1 + \epsilon {\bf e}^2$ for sufficiently small $\epsilon > 0$; clearly $\hat{\lambda} = \frac{1}{2}\lambda^+ + \frac{1}{2}\lambda^-$. We would like to construct points $z^{t,+}$ and $z^{t,-}$ for each $t \in \llbracket m \rrbracket$ such that $\hat{z}^t = \frac{1}{2}z^{t,+} + \frac{1}{2}z^{t,-}$, and such that $(\lambda^+,z^{t,+}), (\lambda^-,z^{t,-}) \in \Pi^t$. Then $(\hat{\lambda},\hat{z}^1,\ldots,\hat{z}^m) = \frac{1}{2}(\lambda^+,\hat{z}^{1,+},\ldots,\hat{z}^{m,+}) + \frac{1}{2}(\lambda^-,\hat{z}^{1,-},\ldots,\hat{z}^{m,-})$ is the convex combination of two other feasible points for $\Pi$, and so is not extreme.

        For a given $t \in \llbracket m \rrbracket$, define $E^t = \Set{(k,h) | ({\bf e}^k,h) \in \ext(\Pi^t) }$, which is equivalent to the set of all extreme points of $\Pi^t$.
        As $(\hat{\lambda},\hat{z}^t) \in \Pi^t$, there must exist some $\gamma^t \in \Delta^{E^t}$ where $(\hat{\lambda},\hat{z}^t) = \sum_{(k,z) \in E^t} \gamma^t_{(k,z)} ({\bf e}^k,h)$. As $1,2 \in \supp(\hat{\lambda})$, there must exist some $\tilde{h}^t$ and $\grave{h}^t$ wherein $(1,\tilde{h}^t), (2,\grave{h}^t) \in E^t$ and $0 < \gamma^t_{(1,\tilde{h}^t)}, \gamma^t_{(2,\grave{h}^t)} < 1$.
        Now define
        \[
            \gamma^{t,\pm}_{(k,h)} = \begin{cases} \gamma^t_{(k,h)} \pm \epsilon & k = 1, h = \tilde{h}^t \\ \gamma^t_{(k,h)} \mp \epsilon & k = 2, h = \grave{h}^t \\ \gamma^t_{(k,h)} & \text{o.w.} \end{cases}
        \]
        Note that, as $\gamma^t \in \Delta^{E^t}$, so is $\gamma^{t,\pm} \in \Delta^{E^t}$. Therefore, we may take
        \begin{alignat*}{3}
            z^{t,+} &\defeq \sum_{(k,h) \in E^t} \gamma^{t,+}_{(k,h)} h =&  &\epsilon \tilde{h}^t - \epsilon \grave{h}^t + \sum_{(k,h) \in E^t} \gamma^t_{(k,h)} h \\
            z^{t,-} &\defeq \sum_{(k,h) \in E^t} \gamma^{t,-}_{(k,h)} h =& -&\epsilon \tilde{h}^t + \epsilon \grave{h}^t + \sum_{(k,h) \in E^t} \gamma^t_{(k,h)} h.
        \end{alignat*}
        Then we may observe that $z^{t,+},z^{t,-} \in \Pi^t$, and that $\hat{z}^t = \frac{1}{2} z^{t,+} + \frac{1}{2} z^{t,-}$. Now see that
        \[
            \lambda^\pm = \sum_{(k,h) \in E^t} \gamma^{t,\pm}_{(k,h)} {\bf e}^k = \sum_{(k,h) \in E^t} \gamma^t_{(k,h)} {\bf e}^k \pm \epsilon {\bf e}^1 \mp \epsilon {\bf e}^2 = \hat{\lambda} \pm \epsilon {\bf e}^1 \mp \epsilon {\bf e}^2
        \]
        Therefore, for each $t \in \llbracket m \rrbracket$, we have that $(\lambda^+,z^{t,+}),(\lambda^-,z^{t,-}) \in \Pi^t$, and that $(\hat{\lambda},\hat{z}^t) = \frac{1}{2}(\lambda^+,z^{t,+}) + \frac{1}{2}(\lambda^-,z^{t,-})$. This implies that $(\lambda^+,h^{1,+},\ldots,h^{m,+}),(\lambda^+,h^{1,-},\ldots,h^{m,-}) \in \Pi$ and that $(\hat{\lambda},\hat{z}^1,\ldots,\hat{z}^m) = \frac{1}{2}(\lambda^+,h^{1,+},\ldots,h^{m,+}) + \frac{1}{2}(\lambda^-,h^{1,-},\ldots,h^{m,-})$.
        Therefore, as our original point is a convex combination of two distinct points also feasible for $\Pi$, it cannot be extreme. Therefore, we must have that $\lambda = {\bf e}^v$ for some $v \in \llbracket n \rrbracket$ for any extreme point of $\Pi$.

        Now, assume for contradiction that $\Pi$ has a fractional extreme point. Using property of extreme points just stated, we may assume without loss of generality that this fractional extreme point is of the form $({\bf e}^1,\hat{z}^1,\ldots,\hat{z}^m)$ with $\hat{z}^1 \not\in \bbZ^{r_1}$. As $({\bf e}^1,\hat{z}^1) \in \Pi^1$, then $({\bf e}^1,\hat{z}^1) = \sum_{(v,h) \in E^1} \gamma_{(v,h)} ({\bf e}^v,h)$ for some $\gamma \in \Delta^{E^1}$. Also, as $\Pi^1$ is ideal and $\hat{z}^1$ is fractional, $({\bf e}^1,\hat{z}^1) \notin \ext(\Conv(\Pi^1))$, and so $\gamma$ must  have at least two non-zero components. But then
        \[
             (\hat{\lambda},\hat{z}^1,\hat{z}^2,\ldots,\hat{z}^m) = \sum_{(v,h) \in E^1} \gamma_{(v,h)} ({\bf e}^1,h,\hat{z}^2,\ldots,\hat{z}^m),
        \]
        a contradiction of the points extremality. Therefore, $\Pi$ is ideal.\Halmos
     \endproof

\section{Non-power-of-two bivariate computational results}\label{app:bivariate-pwl-non-power-of-two}
     See Table~\ref{tab:bivariate-pwl-non-power-of-two}.
     \begin{table}[htpb]
         \centering
         \smaller
         \begin{tabular}{r|c|rrr:rrrr}
     $N$ & Metric & \texttt{MC} & \texttt{CC} & \texttt{DLog} & \Log{} & \texttt{LogIB} & \texttt{ZZB} & \texttt{ZZI} \\ \hline
\multirow{4}{*}{6}
 & Mean (s)  & 9.2  & 20.8  & 4.7  & 1.2  & 1.5  & 1.5  & \textbf{1.1}  \\
 & Std  & 12.0  & 33.0  & 3.4  & 0.7  & 1.1  & 1.2  & 0.6 \\
 & Win & 0 & 0 & 0 & 31 & 9 & 12 & \textbf{48}  \\
 & Fail & \textbf{0} & \textbf{0} & \textbf{0} & \textbf{0} & \textbf{0} & \textbf{0} & \textbf{0}  \\
\hline
\multirow{4}{*}{13}
 & Mean (s)  & 1092.9  & 1507.9  & 320.3  & 16.8  & \textbf{16.5}  & 17.3  & 18.1  \\
 & Std  & 729.7  & 535.4  & 478.7  & 18.6  & 15.7 & 18.6  & 19.3  \\
 & Win & 0 & 0 & 0 & 16 & \textbf{26} & 23 & 15  \\
 & Fail & 37 & 58 & 4 & \textbf{0} & \textbf{0} & \textbf{0} & \textbf{0}  \\
\hline
\multirow{4}{*}{28}
 & Mean (s)  & 1768.1  & 1800.0  & 1426.2  & 127.3  & 131.2  & \textbf{113.4}  & 192.7  \\
 & Std  & 139.6  & -  & 513.5  & 174.5  & 188.7  & 129.7 & 254.9  \\
 & Win & 0 & 0 & 0 & 20 & 26 & \textbf{31} & 3  \\
 & Fail & 75 & 80 & 46 & \textbf{0} & \textbf{0} & \textbf{0} & \textbf{0}
         \end{tabular}
         \caption{Computational results for transportation problems whose objective function is the sum of bivariate piecewise linear objective functions on grids of size $N=d_1=d_2$, when $N$ is not a power-of-two.}
         \label{tab:bivariate-pwl-non-power-of-two}
     \end{table}

\end{APPENDICES}

\end{document}